\documentclass[12pt]{amsart}

\usepackage[]{graphicx} \usepackage{mathrsfs} \usepackage{amsmath}
\usepackage{amsfonts} \usepackage{mathtools} \usepackage{cite}
\usepackage{amsthm} \usepackage{enumerate} \usepackage[shortlabels]{enumitem}
\usepackage{mathtools} \usepackage{amssymb} \usepackage{appendix}
\usepackage{systeme} 
\usepackage{subfig}
\usepackage{multirow}

\usepackage{xcolor}

\usepackage{algpseudocode,algorithm,algorithmicx} %

\newcommand*\Let[2]{\State #1 $\gets$ #2}
\algrenewcommand\algorithmicrequire{\textbf{Input:}}
\algrenewcommand\algorithmicensure{\textbf{Output:}}

\theoremstyle{definition} \newtheorem{remark}{Remark}
\theoremstyle{definition} \newtheorem{definition}{Definition}
\theoremstyle{theorem} 

\def\etalchar{}

= -.75truein \footskip = 18pt

\relax 
\bibstyle{alpha}
\bibdata{numerical_ref}
\gdef \@abspage@last{42}

\begin{document}

\title[]{COMPUTING THE INVARIANT CIRCLE AND THE FOLIATION BY STABLE MANIFOLDS FOR A 2-D MAP BY THE PARAMETERIZATION METHOD: \\NUMERICAL IMPLEMENTATION AND RESULTS}

\author{Yian Yao}
\address{School of Mathematics\\ 
  Georgia Institute of Technology \\
  686  Cherry St. \\ 
  Atlanta GA 30332-160}
\email{yyao93@gatech.edu}
\author{Rafael de la Llave} 
\address{School of Mathematics\\ 
  Georgia Institute of Technology \\
  686  Cherry St. \\ 
  Atlanta GA 30332-160}
\email{rafael.delallave@gatech.edu} 
\thanks{Supported in part by NSF DMS-1800241}

\maketitle

\begin{abstract}

We present and implement an algorithm for computing the invariant
circle and the corresponding stable manifolds for 2-dimensional
maps. The algorithm is based on the parameterization method, and it is
backed up by an a-posteriori theorem established in \cite{YaoL21a}.

The algorithm works irrespective of whether the internal dynamics in
the invariant circle is a rotation or it is phase-locked. The
algorithm converges quadratically and the number of operations and
memory requirements for each step of the iteration are linear with
respect to the size of the discretization.

We also report on the result of running the implementation in some
standard models to uncover new phenomena. In particular, we explored a
``bundle merging'' scenario in which the invariant circle loses hyperbolicity
because the angle between the stable directions and the tangent becomes 
zero even if the rates of contraction are separated. 

We also discuss and implement a generalization of the algorithm to 3
dimensions, and implement it on the 3-dimensional Fattened Arnold
Family (3D-FAF) map with non-resonant eigenvalues and present numerical
results. 

\end{abstract}

\textbf{keywords:} 

\section{Introduction} \label{sec_introduction}

In \cite{YaoL21a}, we have developed an algorithm for computing invariant circle
for 2-dimensional maps and the stable manifolds of its points.  
The algorithm is quadratically convergent as a Newton 
method, but the storage requirements and the number of operations per step is only 
proportional to the number of variables used in the discretization. 
(The regular Newton methods need matrix operations that, of course, 
require a storage square the number of variables, and even larger number of 
operations.)
In paper \cite{YaoL21a}, we presented rigorous results on 
the convergence of the algorithm. The results of \cite{YaoL21a} were
in a-posteriori format. They identified  some condition numbers
and showed that if the initial residual is small enough (with respect to
the condition numbers) the method 
converges. 

In this paper, we present details on the numerical implementation of 
the algorithm and present the details on running it in some examples. 

Even if \cite{YaoL21a} and this paper consider the same problem
and they serve as inspiration for each other, the methodology
and the problems considered are very different, so that the overlap is 
rather minimal. 

The basic idea, following the parameterization method 
\cite{CFdlL05, H16}, is to formulate a functional equation 
(See Secction~\ref{subsec_invariance} and \eqref{invariance}) 
for a parameterization of the invariant circle, the stable foliation 
and the dynamics on them.  The paper \cite{YaoL21a} develops a 
very efficient quasi-Newton method (summarized 
here in Section~\ref{subsec_quasiNewton} ).

One of the advantages of the parameterization approach comparing to the
graph-based methods is that the parameterization method can follow the turns
and oscillations of the stable manifolds very easily.
On the other hand, 
because the stable manifolds often turn and become complicated, graph-based methods
is not applicable when the stable manifolds fail to be a graph. 
On top of that, the parameterization method also gives information on the tangents and high order
derivatives, which may be useful to compute intersections. 

Note that the unknowns for 
our invariance equation are functions, so an important task
for the numerical implementations is to decide on discretizations 
that allow efffective treatment and estimate the truncation error, etc. 
For practical purposes, one also needs to discuss storage strategies, 
operation counts and possible paralellization strategies. 

It is remarkable that the same cancellations 
that lead to good estimates in \cite{YaoL21a} 
lower the operation count and the storage requirements
in the implementation presented here. Conversely, the obstructions 
to regularity reappear here as obstructions to use efficient algorithms. 

We have implemented the algorithm and run it into some standard examples 
with the goal of uncovering some new mathematical phenomena 
that happen in the boundary of the existence of the invariant circles. Besides 
the mathematical interest of these problems, this exploration provides 
a good testing ground for our algorithm.
Having the algorithm  backed up by ``a-posteriori'' theorems, allows 
us to be confident in the results even close to breakdown. This
assurance is invaluable when exploring unexpected phenomena. 
 We note however that, in some cases, 
we have carried out numerical explorations when some of the regularity 
etc. hypothesis of the rigorous results are not satisfied. 

In the case considered in this paper, the notion of breakdown is 
rather  subtle since the regularity of the circle decreases
continuously with parameters.
In normal hyperbolicity theory,  it is customary to 
follow \cite{Mane78} and consider the breakdown when 
normal hyperbolicity is lost and the invariant object ceases 
to be a $C^1$ manifold. In some of our examples, it is natural
to continue the invariant  objects even when they are continuous curves
(with 
a H\"older exponent).  To put our results in context, we present 
a small review of known results in Section~\ref{sec_breakdown}.

In the final Section~\ref{sec_example_3d}, we show that our methods adapt quite 
straightforwardly to higher dimensions,  and we present some numerical results for 
an example. Even if the methods apply to many cases in higher dimensions, we 
show  that  new phenomena may appear: resonances 
among  normal eigenvalues, which will require some adaptation. 
We postpone this study.

\subsection{Other results in the literature} 

The numerical study of limit cycles and their properties 
is very vast and we cannot hope to make a systematic survey of all
the field, but in this subsection, we will try to present some results that are 
somewhat close in methodology and compare the differences.  The study of 
the stable foliations is much more recent.

The paper \cite{HdlL13} also considers limit cycles and the stable manifold 
of their points for differential equations in the plane. The method of 
\cite{HdlL13} is also based in studying an invariance equation, developing
efficient quasi-Newton method and implementing it in 2-D autonomous ODEs.

There are, however very substantial differences
between differential equations and maps of the plane. In the case of 
maps, the circle may be \emph{``phase locked''}. That is, it may contain stable and unstable
periodic orbits. One of the strong point of our algorithm is that our algorithm can deal 
in the same footing with the phase locked circles and the circles on 
which the dynamics is conjugate to a rotation or other cases.   
The phase locked circles studied here do not appear in 2D autonomous ODE's in the plane. 
On the other hand, as indicated in Section~\ref{sec_example_3d}, new phenomena 
appear studying 3-dimensional maps.  

From the numerical point of view, the paper \cite{HdlL13} presents a general 
theory as well as  theory optimized for Fourier-Taylor discretizations. 
The actual implementation is done using only Fourier-Taylor series. A shortcoming 
of Fourier methods is that they are not very adaptative (the grid of 
points is uniform)  and 
they become inefficient for spiking phenomena (when important changes 
happen in small scale).  In this 
paper, 
we use a discretization by splines, which is very adaptative
and allows to consider finer grids in the regions where the invariant circle oscillates more. 

A different approach to compute invariant foliations for ODE's in the plane is developed
in \cite{LangfieldKO14}. The method is based  on the study of  two-point boundary value problems 
by collocation methods. It indeed allows to concentrate  the efforts in places where complicated
phenomena are taking place. 

Another paper which studies numerically invariant foliations is 
 \cite{ChungJ15} which develops
numerical algorithms for invariant foliations and indeed makes 
available the package {\tt FOLI8PAK}.  This paper is also based in 
solving a functional equation for the invariant foliation, but the algorithms 
are not based on a Newton method.

For the main example we consider here,the dissipative standard map \eqref{DST},
there is a literature of computing 
invariant circles  such that the dynamics on them is conjugated to 
a rotation. This can be achived by adjusting one parameter in \eqref{DST}. 
Rigorous results are developed in \cite{CallejaCL13,CanadellH17b} and numerical 
explorations are carried out in \cite{CC10,CallejaF12,CanadellH17a,Rand92a,Rand92b,Rand92c}.

In this paper, we take advantage of the fact that our algorithm works 
equally well for phase-locked and for rotational tori, and explore the breakdown in 
families whose rotation number changes from rational -- presumably phase locked -- and 
conjugate to a rotation. 

\subsection{Organization of the Paper}

The paper is organized as follows:

In Section~\ref{sec_algorithm}, we recall
briefly the main algorithm developed in \cite{YaoL21a}, and present the
iterative steps in an algorithmic form (Algorithm~\ref{algorithm}). In
Section~\ref{sec_implementation}, we discuss some implementation details. 
We specify the methods of discretizing the functions and indicate some
algorithms to accomplish some of the steps. Notably, in Subsection~\ref{subsec_algorithm_coho},
we present
Algorithm~\ref{algorithm_coho} that accelerates the solution of 
cohomology equations \eqref{coho} which are of the basic building blocks
in the iterative step. 

In Section~\ref{sec_continuation}, we consider parameter dependent
problems.  The a-posteriori format of the results in \cite{YaoL21a}
shows that, given enough computer resources, a continuation algorithm
(e.g. Algorithm~\ref{algorithm_continuation}) will reach parameter
values for which the non-degeneracy assumptions of the theorem fails. These
assumptions are basically quantitative versions of hyperbolicity, so
that a continuation method can compute close to the breakdown of 
hyperbolicity (again we recall that the concept of breakdown
depends on the regularity we allow on the invariant circle). 

In Section~\ref{sec_example_2d}, we present the results of running
Algorithm~\ref{algorithm} for the dissipative standard map \eqref{DST}
in computers (with large but not unlimited resources).

 In Section~\ref{sec_breakdown}, we briefly discuss the breakdown of
 the invariant circle when the perturbation is large.

In Section~\ref{sec_example_3d}, we generalize
Algorithm~\ref{algorithm} to 3-dimensional
maps and present some numerical results regarding the 3D-FAF map for
non-resonant eigenvalues. We uncover a new dynamical phenomenon 
that appears when we increase the dimension, but, we postone a full 
numerical treatment.

Finally, in the Appendix \ref{appendix}, we present a brief discussion on
a propose of a parallel implementation of our algorithm.

\section{The Numerical Algorithm} \label{sec_algorithm} In this section, we
concisely discuss the setup of the problem, the invariance equation for the
invariant circle and
the corresponding stable manifolds, and
recall the main algorithm developed in \cite{YaoL21a} to solve it. 
To avoid repetitions with \cite{YaoL21a}, we shortened significantly
the presentation. More details are available in \cite{YaoL21a}.

\subsection{The Invariance Equation} \label{subsec_invariance} Given a smooth
diffeomorphism $f : \mathbb{T} \times \mathbb{R} \rightarrow \mathbb{T} \times
\mathbb{R}$ that generates a discrete dissipative dynamical system in
$\mathbb{T} \times \mathbb{R}$ with an invariant circle, the goal 
is to develop an algorithm  for computing both the invariant circle and
corresponding isochrons.

\begin{remark}
  Following \cite{YaoL21a}, we remark that in the phase-locking phenomenon, the
  isochrons (\cite{W74}) are
  different from the foliation by the stable manifolds. In fact,	
  both the tangent and normal direction of the periodic orbit on the invariant	
  circle are attractive. Nonetheless, in the paper, we abuse the term ``isochron''	
  a little bit and refer ``isochron'' as the leaves of the foliation by the	
  stable manifolds. 
\end{remark}

Following the idea of the parameterization method, we seek an injective immersion
$W: \mathbb{T} \times \mathbb{R} \rightarrow \mathbb{T} \times \mathbb{R}$,
which parameterizes the neighborhood of the invariant circle,
a homeomorphism $a: \mathbb{T} \rightarrow \mathbb{T}$, which describes the
internal dynamics on
the invariant circle, and $\lambda: \mathbb{T} \times \mathbb{R} \rightarrow
\mathbb{R}$, which describes the dynamics on the isochrons, such that
\begin{equation} \label{invariance_origin} f \circ W(\theta, s) - W(a(\theta),
\lambda(\theta, s)) = 0.
\end{equation}

We emphasize that equation \eqref{invariance_origin} is a functional equation. 
The unknowns are the functions $W,a,\lambda$. Hence, 
to treat it in a computer, we will need to specify how to discretize it
(as discussed in Section~\ref{function_representation}). 
  Note that the 
unknowns appear composed with each other.

Taking advantage of the underdetermination  of  the invariance equation,
\eqref{invariance_origin}, it  is equivalent to
\begin{equation} \label{invariance} f \circ W(\theta, s) - W(a(\theta),
\lambda(\theta) s) = 0,
\end{equation}
where $\lambda(\theta, s)$ in equation \eqref{invariance_origin} is reduced to
be a linear function w.r.t. $s$: $\lambda(\theta)s$ (see \cite{YaoL21a}). For the rest
of the paper, the function $\lambda$ is always refering to $\lambda(\theta)$.

\begin{remark}
  \label{undertermincy}
  It can been shown that there are two sources of underdetermincy from equation
  \eqref{invariance}:
  \begin{itemize}
  \item Conjugacy on $\theta$: $\widetilde{W}(\theta, s) = W(h(\theta), s), \ \widetilde{a}(\theta) =
    h^{-1} \circ a \circ h(\theta), \text{ and } \widetilde{\lambda}(\theta) =
    \lambda(h(\theta))$ where $h: \mathbb{T} \rightarrow \mathbb{T}$ is a
    diffeomorphism.
  \item Conjugacy on $s$: $\widehat{W}(\theta, s) = W(\theta, \phi_{\theta}(s)), \ 
    \widehat{a}(\theta) = a(\theta), \text{ and } \widehat{\lambda}(\theta)s =
    \phi^{-1}_{a(\theta)} \circ (\lambda(\theta)\phi_{\theta}(s))$, where
    $\phi_{\theta}(s) = \phi(\theta, s)$ is strictly increasing for any given
    $\theta \in \mathbb{T}$, and $\phi_{\theta}(0) = 0$.
  \end{itemize}
  We remark that all the underdetermincies are up to the conjugacy of some
  homeomorphism. Because of this, properties such as the rotation number of
  $a(\theta)$ are preserved.
\end{remark}

\begin{remark}
  \label{remark_constant_lambda}
  As discussed in \cite{YaoL21a}, taking advantage of the underdeterminancy, when the
  internal dynamics $a(\theta)$
  is conjugate to an irrational rotation $\theta + \omega$, where $\omega$ is
  Diophantine, it can
  be shown that the dynamics on
  the isochron 
  $\lambda(\theta)$ can be reduced to a constant $\lambda$, and thus the
  invariance equation is reduced to
  \begin{equation*}
    f \circ W(\theta, s) - W(\theta + \omega, \lambda s) = 0,
  \end{equation*}
  where $\omega$ is the rotation number of $a(\theta)$, and one need to adjust
  parameters of the map to maintain the rotation number.
  These algorithms to find invariant curves and invariant circles 
  when the inner dyamics is a fixed rotation have been
  studied before  \cite{CanadellH17a,CanadellH17b}. 

\end{remark}

Equation \eqref{invariance} indicates that $W(\theta, 0)$ is the
parameterization of
the invariant circle, and $\lambda(\theta)$ describes the dynamics on the
isochrons. More precisely, the isochrons are $I_{\theta} = \{W(\theta, s)\ | \ s
\in \mathbb{R}\}$ for given $\theta \in \mathbb{T}^1$, and we have the orbit of
$I_{\theta}$ under $f$ converges
exponentially fast to the orbit of $W(\theta, 0)$ at the same phase because of the
dissipative property.

\begin{remark} There are many studies of invariant circles
in the case that the dynamics in the circle is conjugate to
irrational rotations.  In this paper, we also allow $a(\theta)$ to be
phase-locked (i.e. it has an attractive periodic orbit). In such a
case, it can happen (indeed, one expects that this is the most common
case in applications) that the invariant circle is only finitely
differentiable even if the map is analytic. See
 Section~\ref{sec:phaselocked} and \cite{Llave97}.
\end{remark}

\subsection{The quasi-Newton Method} \label{subsec_quasiNewton} In this
subsection, we explain how to perform one step of our quasi-Newton method to solve
\eqref{invariance} for $W(\theta, s), a(\theta)$ and $\lambda(\theta)$. More
detailed discussions can be found in \cite{YaoL21a}. 

\subsubsection{Derivation of the Iterative Step} \label{derivation}

Assume that we have an
approximate parameterization of the neighborhood of the invariant circle
$W(\theta, s)$, an approximate expression of the internal dynamics $a(\theta)$
and an approximate dynamics on the isochrons $\lambda(\theta)$ such that
\begin{equation} \label{error} e(\theta, s) = f \circ W(\theta, s) -
W(a(\theta), \lambda(\theta)s),
\end{equation}
where $e(\theta, s)$ is the error for the invariance equation
\eqref{invariance}.
The goal of the quasi-Newton method is to calculate the
corrections $\Delta_W(\theta, s), \Delta_a(\theta) \text{ and }
\Delta_{\lambda}(\theta)$ such that
\begin{equation} \label{newton0} f(W + \Delta_W)(\theta, s) - (W + \Delta_W)((a
+ \Delta_a)(\theta), (\lambda + \Delta_{\lambda})(\theta)s) = 0
\end{equation} up to quadratically small error through expansions up to first
order.

We consider the correction of the torus, $\Delta_W$ expressed in the frame $DW$.
That is, rather than seeking $\Delta_W$, we seek $\Gamma$ related to $\Delta_W$
by 
\begin{equation} \label{adapted_frame}
  \Delta_W(\theta, s) = DW(\theta, s) \Gamma(\theta,s)
\end{equation}

By substituting the new approximation $W(\theta, s) + \Delta_W(\theta)$, $a(\theta)
+ \Delta_a(\theta)$ and $\lambda(\theta) + \Delta_{\lambda}(\theta)$, using first-order
Taylor expansion,
and taking into account of the derivative of equation~\eqref{error}, we have
\begin{equation} \label{error2} De(\theta, s) = Df(W(\theta, s))D(W(\theta, s))
- DW(a(\theta), \lambda(\theta) s) \begin{pmatrix} Da(\theta) & 0 \\ D
\lambda(\theta) s & \lambda(\theta)
  \end{pmatrix},
\end{equation} 
which leads us to the following cohomological equation: 
\begin{equation} \label{cohom}
  \begin{pmatrix} Da(\theta) & 0 \\ D \lambda(\theta) s & \lambda(\theta)
  \end{pmatrix} \Gamma(\theta, s) - \begin{pmatrix} \Delta_a(\theta) \\
    \Delta_{\lambda}(\theta) s
  \end{pmatrix} - \Gamma(a(\theta), \lambda(\theta)s) = \widetilde{e}(\theta,
  s),
\end{equation}
where $\widetilde{e}(\theta, s) \triangleq - (DW(a(\theta),
\lambda(\theta) s)) ^ {-1} e(\theta, s)$.


\begin{remark}
  In the derivation of Equation \eqref{cohom}, the term $D
[\Delta_W(a(\theta), \lambda(\theta)s)] \begin{pmatrix} \Delta_a(\theta) \\
\Delta_{\lambda}(\theta)s
\end{pmatrix}$ and $De(\theta, s)\Gamma(\theta, s)$ are omitted because it is
``heuristically'' quadratically small (rigorous proof can be found in
\cite{YaoL21a}).
\end{remark}

\subsubsection{Corrections of $W(\theta, s), a(\theta), \lambda(\theta)$} \label{solve_coho}

To compute the corrections of $W(\theta, s), a(\theta)$ and $\lambda(\theta)$,
we solve the cohomological equation \eqref{cohom} for $\Gamma(\theta,
s), \Delta_a(\theta)$ and $\Delta_{\lambda}(\theta)$, which is solving for
$\Gamma_1(\theta, s)$, $\Gamma_2(\theta, s)$, $\Delta_a(\theta)$ and
$\Delta_{\lambda}(\theta)$ from
\begin{equation} \label{eq1} Da(\theta) \Gamma_1(\theta, s) - \Delta_a(\theta) -
\Gamma_1(a(\theta), \lambda(\theta) s) = \widetilde{e}_1(\theta, s),
\end{equation}
\begin{equation} \label{eq2} \lambda(\theta) \Gamma_2(\theta, s) -
\Delta_{\lambda}(\theta) s - \Gamma_2(a(\theta), \lambda(\theta) s) = M(\theta,
s),
\end{equation} where $M(\theta, s) = \widetilde{e}_2(\theta, s) - D
\lambda(\theta) s \Gamma_1(\theta, s).$

In this paper, we will not carry out estimates, but  the conditions for 
the results in \cite{YaoL21a} involve measurements expressed in the norms
below. 
Given $r \in \mathbb{R}, r > 0$, recall
the definition of the Banach space
$\mathcal{X}^{r, \delta}$ in \cite{YaoL21a} for some given $\delta > 0$: 

\begin{definition}
  \label{Xr_space} For a function $u(\theta, s)$ with domain $\mathbb{T} \times
[- \delta, \delta]$, we say $u \in \mathcal{X}^{r, \delta}$ if $u(\theta, s) =
\sum_{j = 0}^{\infty}u^{(j)}(\theta)s^j$ with $u^{(j)}(\theta) \in C^r$ and
$\sum_{j = 0}^{\infty}\left\|u^{(j)}\right\|_{C^r}\delta^{j} < \infty$. In other
words,
  \begin{align*} \mathcal{X}^{r, \delta} = \Big\{ u(\theta, s) = &\sum_{j =
0}^{\infty}u^{(j)}(\theta)s^j \mid u^{(j)}(\theta) \in C^r, \text{ and } \sum_{j
= 0}^{\infty}\left\|u^{(j)}\right\|_{C^r}\delta^{j} < \infty \Big\}
  \end{align*} with norm
  \begin{equation*} \left\|u\right\|_{\mathcal{X}^{r, \delta}} = \sum_{j =
0}^{\infty}\left\|u^{(j)}\right\|_{C^r}\delta^{j}.
  \end{equation*}
\end{definition}

\begin{remark}
  The Definition~\ref{Xr_space} is designed to adapt to the numerical
  discretization discussed in Section~\ref{function_representation}, and it is
  our numerical norm of choice in this paper.

We call attention that $\mathcal{X}^r$ has very anisotropic regularity. It is analytic 
in one variable but finite differentiable in another. When the map is 
\eqref{DST}, which is an entire function, it is shown in \cite{YaoL21a}
that, for the $W$, the $\delta$ can be taken arbitrarily large. On the 
other hand, we will see that the regularity $r$ changes and is important for 
the study of breakdown. 
\end{remark}

For every solution $(\Gamma_1(\theta, s), \Gamma_2(\theta, s), \Delta_a(\theta),
\Delta_{\lambda}(\theta))$ in $\mathcal{X}^{r, \delta} \times \mathcal{X}^{r, \delta}
\times C^r \times C^r$ of \eqref{eq1}, \eqref{eq2}, we can discretize
$\Gamma_1(\theta, s)$, $\Gamma_2(\theta, s)$, $\widetilde{e}_1(\theta, s)$ and
$M(\theta, s)$ in Taylor series w.r.t. $s$ with coefficients are $C^r$ functions
in $\theta$ (as in Definition~\ref{Xr_space}). By equating same order powers, we
have that Equation \eqref{eq1},
\eqref{eq2} are equivalent to the following sequence of equations, where the
lower order equations require a bit extra attention as they also contain the
information for $\Delta_a$ and $\Delta_{\lambda}$.

\begin{itemize}
\item For equation \eqref{eq1}:
  \begin{itemize}
  \item[$\circ$] For the coefficients of $s^0$:
    \begin{equation} \label{eq1_order0} Da(\theta) \Gamma_1^{(0)}(\theta) -
\Gamma_1^{(0)}(a(\theta)) - \Delta_a(\theta) = \widetilde{e}_1^{(0)}(\theta);
    \end{equation}
  \item[$\circ$] For the coefficients of $s^j$, $j \geq 1, j \in \mathbb{N}$:
    \begin{equation} \label{eq1_higher} \Gamma_1^{(j)}(\theta) =
\frac{\lambda^j(\theta)}{Da(\theta)} \Gamma_1^{(j)}(a(\theta)) +
\frac{\widetilde{e}_1^{(j)}(\theta)}{Da(\theta)}.
    \end{equation}
  \end{itemize}
\item For equation \eqref{eq2}:
  \begin{itemize}
  \item[$\circ$] For the coefficients of $s^0$:
    \begin{equation*} \lambda(\theta) \Gamma_2^{(0)}(\theta) -
\Gamma_2^{(0)}(a(\theta)) = M^{(0)}(\theta),
    \end{equation*} which can be rewrite as
    \begin{equation} \label{eq2_order0} \Gamma_2^{(0)}(\theta) =
\lambda(a^{-1}(\theta)) \Gamma_2^{(0)}(a^{-1}(\theta)) -
M^{(0)}(a^{-1}(\theta)),
    \end{equation}
  \item[$\circ$] For the coefficients of $s^1$:
    \begin{equation} \label{eq2_order1} \lambda(\theta) \Gamma_2^{(1)}(\theta) -
\Gamma_2^{(1)}(a(\theta)) \lambda(\theta) - \Delta_{\lambda}(\theta) =
M^{(1)}(\theta),
    \end{equation}
  \item[$\circ$] For the coefficients of $s^j$, $j \geq 2, j \in \mathbb{N}$:
    \begin{equation} \label{eq2_higher} \Gamma_2^{(j)}(\theta) = \lambda ^
{j-1}(\theta) \Gamma_2^{(j)}(a(\theta)) +
\frac{M^{(j)}(\theta)}{\lambda(\theta)}.
    \end{equation}
  \end{itemize}
\end{itemize}

\paragraph{2.2.1.1. Solving for $\Delta_a(\theta)$ and
  $\Delta_{\lambda}(\theta)$:} 

Because of the underdeterminancy of \eqref{eq1_order0} and \eqref{eq2_order1}, the
choices for $\Delta_a(\theta)$ and $\Delta_{\lambda}(\theta)$ are not unique. In this paper, we set
$\Gamma_1^{(0)}(\theta) = 0$ and $\Gamma_2^{(1)}(\theta) = 0$, it follows that
$\Delta_a(\theta) = - \widetilde{e}_1^{(0)}(\theta)$ and $\Delta_{\lambda}(\theta)
= - M^{(1)}(\theta)$. This is a reasonable choice of
solution as the norm of the
correction is boundbed by the error.

An interesting question that requires further exploration on
whether there is a way to choose the
a solution of the equation above
that leads to a more numerically stable algorithm. 

\paragraph{2.2.1.2. Solve for $\Gamma_1^{(j)}(\theta), j \neq 0$ and $
\Gamma_2^{(k)}(\theta), j \neq 1$:}

Notice that equation \eqref{eq1_higher}, \eqref{eq2_order0} and
\eqref{eq2_higher} have been reorganized as above such that they share the
general form of a cohomological equation as follows
\begin{equation} \label{coho} \phi(\theta) = l(\theta)\phi(a(\theta)) +
\eta(\theta),
\end{equation} where $\phi(\theta)$ is the unknown and $a(\theta),
\lambda(\theta)$ and $\eta(\theta)$ are given.

As indicated in Algorithm~\ref{algorithm}, the
quasi-Newton method amounts to solving equations of
the form~\eqref{coho} and much more standard operations such
as algebraic operations and derivatives. 

In this paragraph, we  discuss the solution of  Equation
\eqref{coho}. We
postpone the  presentation of a fast 
numerical algorithm 
Algorithm~\ref{algorithm_coho} to  Section~\ref{subsec_algorithm_coho}.

By inductively replacing $\phi(\theta)$ on the right
hand side of \eqref{coho} by the equation itself, we have
\begin{align} \label{coho_steps} \phi(\theta) =\ & \eta(\theta) +
l(\theta)\eta(a(\theta)) + l(\theta)l(a(\theta))\eta(a^{\circ 2}(\theta))\nonumber \\ & + \ldots +
l(\theta)l(a(\theta))l(a^{\circ 2}(\theta))\cdots l(a^{\circ (n - 1)
}(\theta))\eta(a^{\circ n}(\theta))\nonumber \\ & + l(\theta)l(a(\theta))l(a^{\circ 2}(\theta))\cdots
l(a^{\circ n}(\theta))\phi(a^{\circ (n + 1)}(\theta)) \nonumber \\ = & \sum_{j =
0}^{n}l^{[j]}(\theta)\eta(a^{\circ j}(\theta)) + l^{[n + 1]}(\theta)\phi(a^{\circ (n + 1)}(\theta)),
\end{align} where $l^{[j]}(\theta) = l(\theta)l(a(\theta))l(a^{\circ 2}(\theta)) \cdots
l(a^{\circ (j - 1)}(\theta))$, and $l^{[0]}(\theta) = 1$.

As shown in \cite{YaoL21a}, for small enough $r \in \mathbb{R}, r \geq 0$ such that 
$$\left\|l \right\|_{C^0} \left\|Da \right\|_{C^0}^r < 1,$$
there is a $C^r$ solution for Equation \eqref{coho}. More specifically, 
we have that
$\sum_{j = 0}^{n}l^{[j]}(\theta)\eta(a^{\circ j}(\theta))$ converges uniformly and
absolutely in the $C^r$ space as $n \rightarrow \infty$, and $$\left\| l^{[n +
1]}(\theta)\phi(a^{\circ (n + 1)}(\theta)) \right\|_{C^r} \rightarrow 0,$$ hence
equation \eqref{coho} has an unique $C^0$ solution which is in the $C^r$ space:
\begin{equation} \label{coho_solution} \phi(\theta) = \sum_{j =
0}^{\infty}l^{[j]}(\theta)\eta(a^{\circ j}(\theta)).
\end{equation}

\begin{remark} \label{remark_dynamical_average}
  A sufficient condition that implies
  that $\phi(\theta)$  in~\eqref{coho_solution}
  is a $C^0$ solution, one obvious
  criterior is that $\left\| l\right\|_{C^0} < 1$.
  A more general sufficient condition is: 
\begin{equation}\label{dynamical_average} 
\l^* = \lim_{n \rightarrow \infty} (\left\| \l^{[n]} \right\|_{C^0})^{\frac{1}{n}} < 1,
\end{equation}
We will refer $l^*$ as the \textbf{dynamical average} of $l(\theta)$ w.r.t.
$a(\theta)$.

Since $l^{[n + k]} = l^{[n]}\circ a^{\circ k} l^{[k]}$
we have $\ |l^{[n + k]}\|_{C^0} \le \| l^{[n]}\|_{C^0}\| l^{[k]}\|_{C^0} $
and the standard subadditive argument shows the limit  in
\eqref{dynamical_average}  always exists. 
\end{remark}

\subsection{The Algorithm for One Iteration of the quasi-Newtons
  Method} \label{subsec_algorithm} Given
the approximate $W(\theta, s), a(\theta)$ and $\lambda(\theta)$, we now
summarize the pseudo-code of the algorithm in Section~\ref{subsec_quasiNewton}
to compute $\Delta_W(\theta, s), \Delta_a(\theta)$ and
$\Delta_{\lambda}(\theta)$, see Algorithm~\ref {algorithm}. We truncate
functions in $\mathcal{X}^{r,
  \delta}$ up to the $L$-th order in the power series expansion
w.r.t. $s$ (from the analysists's point of view, $L = \infty$).

\begin{algorithm}[!htb]
  \caption{One iteration of the algorithm}
  \label{algorithm}
  \begin{algorithmic}[1]
    \Require{Initial $W(\theta, s), a(\theta)$ and $\lambda(\theta)$}
    \Ensure{Solution $W(\theta, s), a(\theta)$ and $\lambda(\theta)$ to the
    invariance equation \eqref{invariance}}
  \Statex
  \Let{$\sum_{j = 0}^Le^{(j)}(\theta)s^j = e(\theta, s)$}{$f \circ
    W(\theta, s) - W(a(\theta), \lambda(\theta)s)$}, 
  \State Compute $DW(\theta, s)$ and $DW \circ (a(\theta), \lambda(\theta)s)$,
  \Let{$\sum_{j = 0}^L\widetilde{e}^{(j)}(\theta)s^j = \widetilde{e}(\theta, s)$}{$
    (DW(a(\theta), \lambda(\theta)s))^{-1} e(\theta, s)$  (Section~\ref{etilde})},
  \Let{$\Delta_a(\theta)$}{$- \widetilde{e}_1^{(0)}(\theta)$},
  \Let{$\Gamma_1^{(0)}(\theta)$}{0},
  \State Solve $\Gamma_1^{(j)}(\theta)$ from equation \eqref{eq1_higher} for $1 \leq
  j \leq L$ using Algorithm~\ref{algorithm_coho},
  \Let{$\sum_{j = 0}^L M^{(j)}(\theta)s^j = M(\theta, s)$}{$\widetilde{e}_2(\theta,
    s) - D \lambda(\theta) s \Gamma_1(\theta, s)$},
  \Let{$\Delta_{\lambda}(\theta)$}{$- M^{(1)}(\theta)$},
  \Let{$\Gamma_2^{(1)}(\theta)$}{0},
  \State Solve $\Gamma_2^{(0)}(\theta)$ from equation \eqref{eq2_order0} using
  Algorithm~\ref{algorithm_coho},
  \State Solve $\Gamma_2^{(j)}(\theta)$ from equation \eqref{eq2_higher} for $2
  \leq j \leq L$ using Algorithm~\ref{algorithm_coho},
  \Let{$\sum_{j = 0}^L\Delta_W^{(j)}(\theta)s^j = \Delta_W(\theta,
    s)$}{$DW(\theta, s) \Gamma(\theta, s)$},
  \Let{$W(\theta, s)$}{$W(\theta, s) + \Delta_W(\theta, s)$},
  \Let{$a(\theta)$}{$a(\theta) + \Delta_a(\theta)$},
  \Let{$\lambda(\theta)$}{$\lambda(\theta) + \Delta_{\lambda}(\theta)$},
  \State Return updated $W(\theta, s), a(\theta)$ and $\lambda(\theta)$.

\end{algorithmic}
\end{algorithm}

The Algorithm~\ref{algorithm_coho} used in Step (6) (10) and (11) refers to an
algorithm of solving the
cohomological equation~\eqref{coho} base on Equation~\eqref{coho_solution}. As
discussed in Section~\ref{subsec_algorithm_coho}, this can be implemented in a
much faster way than plain summation.

To solve the invariance equation \eqref{invariance}, one simply iterates the
steps in Algorithm~\ref{algorithm} until either
the error $\left\| e \right\|$ is small enough (the algorithm converges) or
$\max\{\left\| W \right\|, \left\| a \right\|, \left\|\lambda \right\|\}$ exceed some
certain value (the algorithm fails to converge).

One advantage of Algorithm~\ref{algorithm} is that both the time and
the memory requirement are $\mathcal{O}(N \times L)$ for a single step of the
iteration, where, again, $N$ is the size of the
grid for $\mathbb{T}^1$, and $L$ is the order of truncation of $W(\theta, s)$.
This is because all the steps in the algorithm are about summation,
multiplication, division, and spline interpolation (as well
as solving cohomology equations \eqref{coho}). Moreover, this algorithm can
then be implemented very easily in high-level languages.

\begin{remark}
  In Algorithm~\ref{algorithm}, cohomological equations are solved using
  Algorithm~\ref{algorithm_coho}. We remark that the while loop in
  Algorithm~\ref{algorithm_coho} will only repeat finite times (bounded above).
  In our implementation (Section~\ref{subsub_run_time_experiment}), such while
  loop only repeat at most 10 times, which is equivalent to applying the
  contraction on the
  cohomological equation $2^{10}$ times, which is sufficient in most of the
  cases even when
  the tolerance is close to round-off error or when the contraction is slow.
\end{remark}

\begin{remark}
  In the case when $a(\theta)$ admits Diophantine rotation number, Fourier
  Transform is commonly used \cite{HdlL13, ZdlL18, GYLlave20, H16}. In these
  scenarios,
  despite the operations required becomes $\mathcal{O}(N\log N)$ for each
  iteration, the constant is
  smaller and the implementation is indeed faster than the spline interpolation.
\end{remark}

\section{Some Implementation Details} \label{sec_implementation}

The implementation of Algorithm~\ref{algorithm} requires some practical
considerations in terms of function representation and functional operations.
In this section, we provide some implementation details
for the algorithm. We start with
the representation of the functions, followed by the discussion of some basic
functional operations: composition, inverse, etc. We then discuss the
algorithm for solving
the cohomological equation. At the end of this section, We also propose a draft
method for parallel
implementation.

\subsection{Function Representation} \label{function_representation} The first
thing we shall do is to choose a way of discretizing the functions. In order
to perform the algorithm, there are two
types of functions we need to deal with: Type-1: $f
\in C^r(\mathbb{T}, \mathbb{R})$; Type-2: $g \in \mathcal{X}^{r, \delta}\big(\mathbb{T} \times [-\delta, \delta]
, \mathbb{R}\big)$.

\subsubsection{Type-1 Functions}
There are two major methods for discretizing functions of Type-1:
\begin{itemize}
  \item Method 1:
Discretizing $\mathbb{T}$ to a grid of points and storing $f(\theta)$ by the
values on the grid.
In this case, one can obtain function evaluations and
derivatives through interpolation techniques (more specifically, periodic
splines);

\item Method 2:
Representing the function under an orthonormal
basis and storing the coefficients
(for example, Fourier coefficients). To store functions
under the spectral representation, one only needs to truncate the series to a
suitable order. The evaluation and the derivatives of functions can then be computed
accordingly.
\end{itemize}
 
In this paper, we will use method 1 to store $f(\theta)$ for the following
reasons:
\begin{itemize}
\item Since the functions we are dealing with may have spikes or may lose
  regularities in the neighborhood of certain points when the parameters of the
  map come close to the breakdown value, using splines can allow us to partition
  $\mathbb{T}$  in a non-even manner to cope with these situations, which is
  not easily achieved by Fourier Transform.
\item Solving \eqref{eq1} and \eqref{eq2} requires computing the
  composition of two functions, say $f \circ a$, where the internal dynamics,
  $a(\theta)$, in
  general does not conjugate to a rotation, thus it is complicated to use Fourier
  Transform (it is still doable, one can refer to \cite{GYLlave20} for further
  discussions). On the other hand, storing $f$ and $a$ via grid points produces
  a simpler, faster and  more reliable way to compute the composition $f \circ a$.
\end{itemize}

\begin{remark}
  One delicate point of using grid points representation for functions is when the
  regularity of the function drops below the order of splines used in
  interpolation. This happens when the perturbation is close to the breakdown.
\end{remark}

\begin{remark}
  \label{rmk_spline_error}
  Classic results (\cite{H76}) for the error of cubic spline approximation shows
  that
  \begin{equation*}
    \left\| (g - \widehat{g})^{(r)}\right\|_{C^{\infty}} \leq C_r \left\| g^{(4)}\right\|_{C^{\infty}} \Big(\frac{1}{N}\Big)^{4 - r},
  \end{equation*}
  where $g \in C^4$, $\widehat{g}$ is the cubic spline approximation, $C_0 =
  \frac{5}{384}$, $C_1 = \frac{1}{24}$ and $C_2 = \frac{3}{8}$. Thus the
  accuracy drops as the regularity increases. The norm in Definition~\ref{Xr_space} is hence
  affected by this round-off error (as indicated in Subsection~\ref{subsec_convergence_rate}).
\end{remark}

\subsubsection{Type-2 Functions}
Following Definition~\ref{Xr_space}, function $g \in
\mathcal{X}^{r, \delta}$ of Type-2, $g(\theta,s)$ can be written as the Taylor's series
w.r.t. $s$, i.e. $g(\theta, s) = \sum_{j = 0}^{\infty}g^{(j)}(\theta)s^j$, we can
truncate $g(\theta, s)$ up to the $L$-th order in order to store it, provided
that $L$ is big enough so that $\sum_{j =
  0}^{L}g^{(j)}(\theta)s^j$ is a good approximation of $g(\theta, s)$. In this
case, storing $g(\theta, s)$ is equivalent as storing $L + 1$ functions of
Type-1 as $g^{(j)}(\theta) \in C^r(\mathbb{T}, \mathbb{R})$ for $j = 0, 1,
\dots, L$, which is essentially storing a 2-d array with size $N \times L$.

\subsection{Composition Between Functions}
As indicated in the algorithm, we need to cope with the composition between
functions both in the $C^r$ space and the $\mathcal{X}^{r, \delta}$ space .

\subsubsection{Coping with Function Composition in $C^r$}
\label{one-d-composition} 
In Algorithm~\ref{algorithm}, the composition between two $C^r$ functions is
required when computing the error of the invariance equation and when deriving
and solving the cohomological equations.

In these cases,  such operation
can be abstracted as computing
$f \circ g$, where $f:\mathbb{T} \rightarrow
\mathbb{R}$ can be functions either of index 0 or 1, and $g: \mathbb{T}
\rightarrow \mathbb{T}$ is always of index 1. Recall index 0 function $f$ satisfies
$f(\theta  + 1) = f(\theta)$, which is equivalent as periodic function, and
index 1 function $f$ satisfies $f(\theta + 1) = f(\theta) + 1$.

If $f(\theta)$ has index 0, $f \circ g$ can be
calculated by splines with periodic boundary conditions. If $f(\theta)$ has index 1, $f
= id + \hat{f}$, where $\hat{f}(\theta)$ is periodic, and thus $f \circ g = g +
\hat{f} \circ g$. In other words, index 1 functions appears as diffeomorphisms
of the circle to itself, while index 0 functions maps from circle to real numbers.

\begin{remark}
  In the unknowns of equation \eqref{invariance}, only $W_1^{(0)}(\theta)$ and
  $a(\theta)$ has index 1, while the other functions: $W_1^{(j)}(\theta)$ for $j = 1,
  \dots L$, $W_1^{(j)}(\theta)$ for $j = 0, \dots, L$ and $\lambda(\theta)$ have
  index 0.
\end{remark}

\begin{remark}
  \label{remark_cubic_spline}
  As  indicated in Section~\ref{subsec_continuation_k},
  the invariant circle losses regularity when the parameter value is close to
    the breakdown. Often, the loss of regularity is 
very localized. This requires using different 
types of splines and we need to have estimates for composition 
for  cubic, quadratic, 
linear  or Akima splines. 
\end{remark}

\begin{remark}
  During the composition between functions, for example, $W^{(0)} \circ
  a(\theta)$, the monotonicity of $a(\theta)$ is required. We remark that although
  the cubic spline does not guarantee such monotonicity, the results in
  \cite{YaoL21a} indicate that the $a(\theta)$ we are using in the
  implementation is accurate enough to assure such monotonicity. Nonetheless,
  although not used in this paper, one
  may use Steffen's interpolation method \cite{S90} if necessary.
\end{remark}

\subsubsection{Coping with Function Composition  in $\mathcal{X}^{r, \delta}$}
We note that the evaluation of  the functional in 
\eqref{invariance} for the example \eqref{DST}  (and many  of 
the  manipulations described in Algorithm~\ref{algorithm}) 
are:
\begin{itemize} 
\item 
Algebraic operations (addition, substraction, multiplication) 
\item Composite $W(\theta, s)$ with functions on the left. 
\item Composite $W(\theta, s)$ on the right by $a(\cdot)$ and $\lambda(\theta) \cdot$. 
\end{itemize} 

Of course, algebraic operations are straightforward but we now 
detail the others operations.
 The paper \cite{GYLlave20} deals with the same
issues of compositions above in similar manners. 

\leftline{\emph{Composite $W(\theta, s)$ with functions on the left}} 
In Section~\ref{sec_example_2d}, we consider the dissipative standard map
\eqref{DST}, so that the only composition in the left needed
is composition with $\sin(\cdot)$. We will detail the composition with 
$\sin(\cdot)$, but it is clear that the method applies to many functions. 
A very throrough discussion of this and other algorithms 
is \cite[Sec. 4.7]{Knuth98} (which traces it back to Euler) 
and \cite[2.3]{H16}.

Given 
$W_1(\theta, s) = \sum_{j =
  0}^{L}W_1^{(j)}(\theta)s^j$, we want to compute an approximation of 
$\sin(W_1(\theta, s)) $.

Let $L$ be the maximum order of $s$ we want to calculate, i.e., we start with
$W_1(\theta, s) = \sum_{j =
  0}^{L}W_i^{(j)}(\theta)s^j$, and the goal is to find $\sin(W_1(\theta, s))$ up to the
$L$-th order. Denote $S(\theta, s) \triangleq \sin(W_1(\theta, s)) = \sum_{j =
  0}^{L}S^{(j)}(\theta)s^j$ and $C(\theta, s) \triangleq \cos(W_1(\theta, s)) = \sum_{j =
  0}^{L}C^{(j)}(\theta)s^j$. By differentiating $S(\theta, s)$ and $C(\theta,
s)$ with respect to $s$, and by noticing that
\begin{align*}
  \partial_sS(\theta, s) &= C(\theta, s) \partial_sW_1(\theta, s), \\
  \partial_sC(\theta, s) &= - S(\theta, s) \partial_sW_1(\theta, s),
\end{align*}
we can have the following iterating formulae for $S^{(j)}(\theta),
C^{(j)}(\theta)$:
\begin{align*}
  S^{(j)}(\theta) &= \frac{1}{j}\sum_{k = 0}^{j - 1}(j - k)W_1^{(j - k)}(\theta) C^{(k)}(\theta), \\
  C^{(j)}(\theta) &= - \frac{1}{j}\sum_{k = 0}^{j - 1}(j - k)W_1^{(j - k)}(\theta) S^{(k)}(\theta).
\end{align*}

Starting with
\begin{align*}
  S^{(0)}(\theta) &= \sin(W_1(\theta, 0)) = \sin(W_1^{(0)}(\theta)),\\
  C^{(0)}(\theta) &= \cos(W_1(\theta, 0)) = \cos(W_1^{(0)}(\theta)).
\end{align*}

we are now able to calculate the coefficient of $\sin(W_1)(\theta, s)$ up to the $L$-th order.

\leftline{\emph{Composite $W(\theta, s)$ on the right with the internal dynamics and the scaled variables}}
Given 
$W_1(\theta, s) = \sum_{j =
  0}^{L}W_1^{(j)}(\theta)s^j$,
we have 
\[
W_1(a(\theta), s) = \sum_{j =
  0}^{L}W_1^{(j)}(a(\theta))s^j
\]
so that the composition refers to composition of one dimensional 
functions discussed in Section \ref{one-d-composition}. 

We have 
\[
W_1(\theta, \lambda(\theta) s) = \sum_{j =
  0}^{L}\left( W_1^{(j)}(\theta) \lambda(\theta)^j \right) s^j
\]
so that the desired operation becomes just an arithmetic operation.

\subsubsection{The Computation of $\widetilde{e}(\theta, s)$} \label{etilde}

For step (4) in Algorithm~\ref{algorithm}, instead of calculating the
inverse of
\begin{equation*}
  DW \circ (a(\theta), \lambda(\theta)s) \triangleq
  \begin{pmatrix} \beta_{11}^{(0)}(\theta) & \beta_{12}^{0}(\theta) \\
    \beta_{21}^{0}(\theta) & \beta_{22}^{0}(\theta) 
  \end{pmatrix},
\end{equation*}
we will solve the following
linear system:
$$(DW(a(\theta), \lambda(\theta)s)) \widetilde{e}(\theta, s) = e(\theta, s),$$
which, after some routine calculation, is just to solve the following linear
system inductively from $k = 0$ to $L$ to get $\widetilde{e}^{(k)}(\theta)$.
\begin{equation*}
  \begin{pmatrix} \beta_{11}^{(0)}(\theta) & \beta_{12}^{0}(\theta) \\
    \beta_{21}^{0}(\theta) & \beta_{22}^{0}(\theta) 
  \end{pmatrix} \begin{pmatrix} \widetilde{e}_1^{(k)}(\theta) \\
    \widetilde{e}_2^{(k)}(\theta) 
  \end{pmatrix} = \begin{pmatrix} \hat{e}_1^{(k)}(\theta) \\
    \hat{e}_2^{(k)}(\theta)
  \end{pmatrix}.
\end{equation*}
where
\begin{equation*}
  \begin{pmatrix} \hat{e}_1^{(k)}(\theta) \\ \hat{e}_2^{(k)}(\theta)
  \end{pmatrix} = \begin{pmatrix} e_1^{(k)}(\theta) - \sum_{j = 0}^{k} (\beta_{11}^{(k-j)}(\theta)\widetilde{e}_1^{(j)}(\theta) + \beta_{12}^{(k-j)}(\theta)\widetilde{e}_2^{(j)}(\theta)) \\ e_2^{(k)}(\theta) - \sum_{j = 0}^{k} (\beta_{21}^{(k-j)}(\theta)\widetilde{e}_1^{(j)}(\theta) + \beta_{22}^{(k-j)}(\theta)\widetilde{e}_2^{(j)}(\theta))
  \end{pmatrix},
\end{equation*}
Thus,
\begin{equation*}
  \begin{pmatrix} \widetilde{e}_1^{(k)}(\theta) \\
    \widetilde{e}_2^{(k)}(\theta) 
  \end{pmatrix} = \frac{1}{\delta(\theta)} \begin{pmatrix} \beta_{22}^{(0)}(\theta) & - \beta_{12}^{(0)}(\theta) \\ - \beta_{21}^{(0)}(\theta) & \beta_{11}^{(0)}(\theta) 
  \end{pmatrix} \begin{pmatrix} \hat{e}_1^{(k)}(\theta) \\
    \hat{e}_2^{(k)}(\theta)
  \end{pmatrix},
\end{equation*}
where $\delta(\theta) = \beta_{11}^{(0)}(\theta)\beta_{22}^{(0)}(\theta) -
\beta_{21}^{(0)}(\theta)\beta_{12}^{(0)}(\theta).$

\subsection{Computation of the Approximate Inverse of the Internal Dynamics}
\label{subsec_inverse_method}

Deriving equation \eqref{eq2_order0} involves the computation of the inverse of
the internal dynamics $a(\theta)$. This is a crucial  piece
of our algorithm. 

We present four methods for computing $a^{-1}$.
Method 1 is to compute $a^{-1}$
directly via reflecting the graph acroess the diagonal 
 using interpolation. Methods 2, 3, and 4 
are perturbative methods, suitable for continuation algorithms.
In method 2, 3, and 4, we assume that $a^+ = a + \Delta_a$ and that we know
$a^{-1}$. We want to compute $a^- \equiv (a^+)^{-1}$. 
Closely related to this is to develop an iterative step
that produces a better approximation given some approximation to $a^-$. 
We will develop several such iterative methods. Note that  these iterative
methods can be applied one after the other, or after the non-perturbative 
one to polish-off the errors. 

\begin{itemize}
\item Method 1: Reflecting the graph.

This method takes advantage of the
fact that the homeomorphism $a: \mathbb{T} \rightarrow \mathbb{T}$ is strictly
increasing and is of index 1, and thus $a(\theta) = T_{a(0)} \circ \hat{a}(\theta)$, where
$T_{\alpha}(x) = x + \alpha$ and $\hat{a}(\theta)$ is a
strictly increasing
function with $\hat{a}(0) = 0$ and $\hat{a}(1) = 1$, thus $\hat{a}$ is invertible. It follows that
$a^{-1}(\theta) = \hat{a}^{-1} \circ T_{a(0)}^{-1} (\theta)$, where
$T_{\alpha}(x) = x - \alpha$. To compute $\hat{a}^{-1}(\theta)$, where 
$\theta$ is a grid of points in $\mathbb{T}$, we can reflect $\hat{a}(\theta)$ over
the line $\hat{a}(\theta) = \theta$ by treating $\hat{a}(\theta)$ as
the new grid of $\mathbb{T}$, treating the original grid of $\theta$ as the
corresponding values for $\hat{a}^{-1}(\theta)$, and then using splines to
evaluate $\hat{a}^{-1}(\theta)$ at the original grid of $\theta$.

This method
performs well provided the original grid is
sufficiently large such that the values for $\hat{a}^{-1}(\theta)$ are spread
over $\mathbb{T}^1$. If $a(\theta)$
has large slope in
certain parts of $\mathbb{T}$, one can use adaptive grid and put more points on
those parts to ensure a good accuracy.  The results can be polished 
off if needed using the iterative methods described below. 

Note that the method of reflecting the graph is non-perturbative and can 
be started with only the knowledge of $a$ and no approximate inverse is needed. 

\def\Id{\text{Id}}

\item Method 2: Compute the  ``Left'' inverse of $a^+$.
  
This method is mentioned in \cite{H16}. 

From
\begin{equation} \label{method2}
  (a^{-1} + \Delta_{a^{-1}}) \circ a^+(\theta) - \theta = 0,
\end{equation}
 
and by omitting the quadratically small terms, we have $ \Delta_{a^{-1}}(\theta) =
- e \circ a^{-1}(\theta)$, where $e(\theta) = a^{-1} \circ a^+(\theta) -
\theta$. It follows that the updated inverse of $a(\theta)$ is $a^{-1}(\theta) +
\Delta_{a^{-1}}(\theta)$.

\item Method 3: Compute ``Right'' inverse of $a^+$.

  Similar to Method 2, we now optimize the right side inverse by solving the
  following objective function:
  \begin{equation} \label{method3}
    a^+ \circ (a^{-1} + \Delta_{a^{-1}})(\theta) - \theta = 0,
  \end{equation}
  which, after omitting quadratically small term, gives $\Delta_{a^{-1}} =
  \frac{\theta - a^+\circ (a^{-1})}{Da^+ \circ (a^{-1})}$.

\item Method 4:
  By combining method 2 and 3, we aim to optimize
  \begin{equation}
    \label{ainv_method4}
    (a^{-1} + \Delta_{a^{-1}}) \circ a^+(\theta) - \beta a^+ \circ (a^{-1} +
  \Delta_{a^{-1}})(\theta) + (\beta - 1) \theta = 0,
  \end{equation}
  where $\beta = \text{ceil}(\left\| \frac{1}{Da^+(a^{-1})} \right\|_{C^0})$, which can 
  be treated as a cohomological equation w.r.t. $\Delta_{a^{-1}}(\theta)$, and can
  then be solved via Algorithm~\ref{algorithm_coho}.  
\end{itemize}

Generally speaking, method 1 tends to have a better performance in the early
iterations,
while methods 2, 3, and 4 work better in the later iterations, as these methods
requires $\Delta_a$ to be small, where $\Delta_a$ is of the same order as the
error for the cohomological equation.
In practice, we try all the methods and pick the one with minimal error.
One can refer to Section~\ref{subsubsection_inverse_method_plot} for some numerical
examples for the performance.

\begin{remark}
  For homeomorphism $a(\theta)$, the left inverse and the right inverse are the
  same. However, the mathematically equivalent equations
  \eqref{method2} and equation \eqref{method3} are very different from
  the numerical point of view. 
  Notably,  the unknown $\Delta_{a^{-1}}$ appears linearly in equation
  \eqref{method2}, but it does not depend differentiably on the independent variable
  $a^+$; on the other hand, $\Delta_{a^{-1}}$ appears non-linearly in equation \eqref{method3},
  and it depend differentiablly on $a^+$.
\end{remark}

\subsection{A Fast Algorithm for Solving Cohomological Equations}
\label{subsec_algorithm_coho}

Calculating $\Gamma^{(j)}_{1, 2}(\theta)$ from step (7) and (10) requires solving the
cohomological equation of the form in equation~\eqref{coho}:
$$\phi(\theta) = l(\theta)\phi(a(\theta)) + \eta(\theta).$$

\begin{remark}
  The cohomological equation \eqref{coho} is a linearization of 
the conjugacy and, therefore it appears in many problems in dynamics 
and in singularity theory. It has been used as the basis of Newton methods, 
and deformation theory. 

When $a(\theta)$ conjugates to a Diophantine rotation,
it can be treated using Fourier series and it admits a solution
even if the $l$ is not contraction.  This particular case, 
is the basis of KAM theory. 
\end{remark}

Recall the discussion in \eqref{coho_steps}, $\sum_{j =
  0}^Ml^{[j]}(\theta)\eta(a^{\circ j}(\theta))$ is a good approximation of the unknown
$\phi(\theta)$ provided the series converges (i.e., when the dynamical average
$\lambda^* < 1$) and $M$ is big enough. 
We can use
the above algorithm to get to $\sum_{j =
  0}^Ml^{[j]}(\theta)\eta(a^{\circ j}(\theta))$ with $\log M$ iterations:

\begin{algorithm}[!htb]
  \caption{Solving the cohomological equation \eqref{coho}}
  \label{algorithm_coho}
  \begin{algorithmic}[1]
    \Require{$l(\theta)$, $a(\theta)$, $\eta(\theta)$ and $tolerance$}
    \Ensure{The solution of equation \eqref{coho}: $\phi(\theta)$}
    \Statex
    \Let{$\phi(\theta)$}{$\eta(\theta)$},
    \Let{$L(\theta)$}{$l(\theta)$},
    \Let{$A(\theta)$}{$a(\theta)$},
    \While{$\left\| \phi(\theta) - l(\theta)\phi(a(\theta)) - \eta(\theta)
      \right\| > tolerance$}
    \Let{$\phi(\theta)$}{$\phi(\theta) + L(\theta) \phi \circ A(\theta)$}
    \Let{$L(\theta)$}{$L(\theta) L \circ A(\theta)$}
    \Let{$A(\theta)$}{$A \circ A(\theta)$}
    \EndWhile
    \State Return $\phi(\theta)$\;
  \end{algorithmic}
\end{algorithm}


\begin{remark}
  We emphasis that Algorithm~\ref{algorithm_coho} allows us to make the
  summation of $M$ terms in $\log M$ steps. We hope such idea can be used  in
  more general applications.
\end{remark}

\begin{remark}
  The main idea of Algorithm~\ref{algorithm_coho} is similar to
  the binary expansion. More specifically, in the begining of the $(k + 1)$-th
  iteration, we have
  \begin{equation*}
    \phi(\theta) = \sum_{j=0}^{2^k}l^{[j]}(\theta)\eta(a^{\circ j}), 
    L(\theta) = l^{[2^k]}(\theta), \text{ and }
    A(\theta) = a^{\circ 2^k}(\theta),
  \end{equation*}
  then inside the while loop, $\phi(\theta), L(\theta)$ and $A(\theta)$ got
  updated to
  \begin{equation*}
    \phi(\theta) = \sum_{j=0}^{2^{k+1}}l^{[j]}(\theta)\eta(a^{\circ j}), 
    L(\theta) = l^{[2^{k + 1}]}(\theta), \text{ and }
    A(\theta) = a^{\circ 2^{k+1}}(\theta).
  \end{equation*}
\end{remark}

\begin{remark}
  An alternative derivation of Algorithm~\ref{algorithm_coho}, which
  yields some extra flexibility is as follows:

  Because of the uniqueness of equation~\eqref{coho}, solving equation~\eqref{coho} is
  the same as solving
  \begin{equation} \label{coho_better}
    \phi(\theta) = l^{[n + 1]}(\theta)\phi(a^{\circ (n + 1)}(\theta)) + \sum_{j = 0}^{n}l^{[j]}(\theta)\eta(a^{\circ j}(\theta)),
  \end{equation}
  for any $n > 0$. Solving equation \eqref{coho_better} for a suitable $n$ gives
  better contraction and the sums
in \eqref{coho_solution} converge faster. 

In other words, for any 
$n$, if  we trade $l$, $a$, $\eta$ 
for  $l^{[n]}$, $a^{\circ n}$, $\sum_{j = 0}^{n}l^{[j]}(\theta)\eta(a^{\circ j}(\theta))$
respectively, we get an equivalent problem with a stronger contraction. Therefore, 
given a given accuracy, the number of terms needed in \eqref{coho_solution} for 
\eqref{coho_better} is a fraction of the number of terms neeed in the  original problem. 

Algorithm~\ref{algorithm_coho} can be described as repeating the trade offs (whith $n=2$) 
till the traded equation has such a strong contraction that one term in the sum
in \eqref{coho_solution} is enough. Of course, performing a trade is
computationally expensive
because it involves composition and multiplication between functions, but each trades
halfs the number of terms needed.

Notice that, each of the trades could be based on a different $n$
(not just in a binary expansion). This extra flexibility is useful
in the case that the dynamics is a rotation. In such a case it is useful
to use as $n$ the numbers in the continued fraction expansion. 

The people familiar with renormalization group will notice the similarity of 
the procedure with the \emph{decimation procedures}. One issue to consider
is the stability of the procedure.
Related algorithms are studied in \cite{HuguetLS13}

\end{remark}

\begin{remark}
  Due to the accumulation of the truncation error and round-off errors,
  Algorithm~\ref{algorithm_coho} may not be able to reach the
required accuracy. To resolve this, one can repeatly apply this algorithm for
$$ \Delta_{\phi}(\theta) = l(\theta)\Delta_{\phi}(a(\theta)) + \widehat{e}(\theta),$$
where $\widehat{e}(\theta) = \phi(\theta) - l(\theta)\phi(a(\theta)) - \eta(\theta)$ is
the error for the cohomological equation from the previous step,
and the solution for such cohomological
equation now becomes $\phi(\theta) + \Delta_{\phi}(\theta)$.
\end{remark}

\begin{remark}
  From the discussion in \cite{YaoL21a}, we have 
  \begin{equation*}
    \left\| \phi - \sum_{j =0}^Ml^{[j]}\eta(a^j)\right\|_{C^r} \leq (r + 1)! \Bigg( \left\| l \right\|_{C^r} + \left\| a \right\|_{C^r})^r\Big(\sum_{j = M + 1}^{\infty}j^{r - 1}(\left\| Da \right\|_{C^0}^r\left\| l \right\|_{C^0}\Big)^n \Bigg)\left\| \eta \right\|_{C^r},
  \end{equation*}
  where the regularity $r$ is bounded above:
  \begin{equation*}
    r < - \frac{\ln \left\| l \right\|_{C^0}}{\ln \left\| Da \right\|_{C^0}},
  \end{equation*}
  indicating that such cohomological equation can only be solved for a finite range
  of regularities.

For the range of regularities that the equation can be solved, the solution is 
bounded from a $C^r$ space to itself. 
\end{remark}

\subsection{Truncation Versus Smoothing Operator}
In the proof of the hard implicit function theorem in \cite{YaoL21a}, the
corrections  to  $W(\theta, s), a(\theta)$ and $\lambda(\theta)$
predicted by the algorithm are modified 
by a smoothing operator. The reason for the need of  smoothing in this problem is different 
than in KAM theory. In KAM theory, the cohomology equations lose derivatives. 
The solutions \eqref{coho_solution}, however, have the same regularity of the data,
but the need of the smoothing comes from the lack of differentiability  the functional involved. 

In our numerical implementation, we did not use any smoothing operator. 
One can argue that the effect of trucation to a finite representation has the same 
effect of smoothing. 
On the other hand, we need to take care that the number of points in the grid is 
rather large so that the effects of the truncation do not affect the result.

\subsection{Validation of the Correctness of the Solution} \label{subsec_validation}
The most naive way of validating the correctness of the solution is to check the
norm of the error from the invariance equation and see if it is of the same
order as the round-off error. However, this is not a thorough
approach.

As mentioned in Section~\ref{function_representation}, one of the drawbacks of using splines is that it
cannot capture all the information of the
function, but only for a selected grid of points. It is possible, especially
when the perturbation is close to the breakdown, that localized singularity
appears in the solution. The nature of the singularity is a highly oscillatory
solution with moderate regularity (Good models of the expected behavior  are
the Weierstrass functions $\sum_n  \alpha^n \sin( \beta^n x) $ with
$0< \alpha < 1$, $\beta \in \mathcal{N}, \beta > 1$ ) 
oscillation. Depending on the choice of the grid, the splines may be smooth (when
the localized singularity is in the gap of the grid) or not (when the
singularity is close to one of the grid points), and the latter scenario leads to
a disaster, especially when the regularity drops below the order of choice of
our spline.

More specifically, in the case when the rotation number is rational, the
localized singularity happens at two hyperbolic points when their
corresponding stable and unstable manifolds meet. This admits a $C^{\alpha}$
manifold, where $0 < \alpha < 1$ is regulated by the eigenvalues at the
attractive periodic orbit. Further discussions can be found in Section~\ref{sec:phaselocked}
and \cite{Llave97}.

By the above discussions, even when the solution we computed induces small
error, it is still not guaranteed that the solution is the true
solution. We considered several methods of validating the correctness of the
solution as follows:

\begin{itemize}
\item Method 1: Grid shifting:
  To guarantee that the grid of choice captures the information of the invariant
  circle and isochrons, i.e., no singularity points in the gap of the grid, we
  can either shift the grid or increase/double the size of the grid (and maybe
  repeatedly doing this) and check
  if the absolute error for the invariance equation remains to be small. In our
  actual implementation, we choose the grid adaptively and use the grid size
  doubling method.

\item Method 2: Monitor the error of the invariance equation:

  Another aspect comes from the proof of our theoritical result in \cite{YaoL21a}
  is that the error for the invariance equation: $e(\theta, s)$ has to satisfy
  $\left\| e\right\|_{C^{m - 2}} \leq v e^{-2 \mu \beta \kappa^n}$ for some
    prescribed positive $v, \mu, \beta > 0,$ and $\kappa > 1$. Thus, if the
    solution is
    valid, the convergence rate of
  our iteration process has to be superlinear.

We emphasize that monitoring the residual of the invariance equation is
not enough. It is very common (especially when we are close to the breakdown 
to obtain many \emph{spurious} solutions of the truncated invariance equation. 
The results of \cite{YaoL21a} show  that a well behaved function that satisfies
very accurately the invariance equation will correspond to a true solution.

\end{itemize}

We have included a brief example regarding the above discussion in
Section~\ref{example_validation}.

\section{Continuation Method} \label{sec_continuation}

By the previous discussions, for a given map $f: \mathbb{T} \times \mathbb{R}
\rightarrow \mathbb{T} \times \mathbb{R}$, starting with some initial
approximation, we are now
ready to apply the iterative Algorithm~\ref{algorithm} several times until the
error is close to the round off error, with the solution for the invariance
equation \eqref{invariance}: $W(\theta, s)$, $a(\theta)$ and $\lambda(\theta)$.

In this section, we do some further discussions of applying our algorithm to the
parameter-dependent problems. The idea is based on the standard continuation method.

\subsection{The Continuation Method} \label{continuation}

In the parameter-dependent problems, it is a standard procedure  to start with a
simple (unperturbed) scenario and perform the continuation method to moving
the parameters gradually towards the desired value.  

The following procedure is quite standard. The only subtlety in our case 
is that the equations for parameterizations are underdetermined 
(the geometric objects considered are unique), so that the result of
running two different continuations can be different (even if they 
parameterize the same geometric object). 

\subsubsection{Basic Idea} \label{continuation_basic}
Let $f_{\epsilon}$ be a parameter family of diffeomorphim, the goal is to find
$W_{\epsilon}$, $a_{\epsilon}$ and $\lambda_{\epsilon}$ such that the invariance
equation:
\begin{equation} \label{invariance_para}
f_{\epsilon} \circ W_{\epsilon}(\theta, s) -
W_{\epsilon}(a_{\epsilon}(\theta), \lambda_{\epsilon}(\theta)s) = 0
\end{equation}
holds, starting from $\epsilon = 0$.

Given the solution for equation \eqref{invariance_para} for some
$\epsilon$: $W_{\epsilon}(\theta, s)$, $a_{\epsilon}(\theta)$ and
$\lambda_{\epsilon}(\theta)$, the goal is to find the starting approximation:
$W_{\epsilon + h}(\theta, s)$, $a_{\epsilon + h}(\theta)$ and $\lambda_{\epsilon
+ h}(\theta)$ for
\begin{equation} \label{invariance_para_h}
f_{\epsilon + h} \circ W_{\epsilon + h}(\theta, s) -
W_{\epsilon + h}(a_{\epsilon + h}(\theta), \lambda_{\epsilon + h}(\theta)s) = 0.
\end{equation}

The naive choice would be the 0-th order approximation, i.e. $W_{\epsilon +
  h}(\theta, s) = W_{\epsilon}(\theta, s)$, $a_{\epsilon + h}(\theta) =
a_{\epsilon}(\theta)$ and $\lambda_{\epsilon + h}(\theta) =
\lambda_{\epsilon}(\theta)$. We remark that this choice is already good enough
if the increment of the perturbation $h$ is small.

One can also try to look for 1-st order approximations, which is looking for
$\frac{\partial W_{\epsilon}}{\partial \epsilon}(\theta, s)$, $\frac{\partial
  a_{\epsilon}}{\partial \epsilon}(\theta)$ and $\frac{\partial
  \lambda_{\epsilon}}{\partial \epsilon}(\theta)$ such that
\begin{align} \label{first_order_approximation}
  W_{\epsilon + h}(\theta, s) &= W_{\epsilon}(\theta, s) + \frac{\partial W_{\epsilon}}{\partial \epsilon}(\theta, s)h, \nonumber \\
  a_{\epsilon + h}(\theta) &= a_{\epsilon}(\theta) + \frac{\partial a_{\epsilon}}{\partial \epsilon}(\theta)h, \nonumber \\
  \lambda_{\epsilon + h}(\theta) &= \lambda_{\epsilon}(\theta) + \frac{\partial \lambda_{\epsilon}}{\partial \epsilon}(\theta)h
\end{align}
satisfies equation \eqref{invariance_para_h} up to quadratic error $\mathcal{O}(h^2)$.

The procedure of this computation is similar to the derivation in
Section~\ref{derivation}: by \eqref{invariance_para} and
\eqref{invariance_para_h}, and
omitting quadratically small terms, we have
\begin{align}
  & Df_{\epsilon} \circ W_{\epsilon} (\theta, s) \frac{\partial W_{\epsilon}}{\partial \epsilon}(\theta, s) - DW_{\epsilon}(a_{\epsilon}(\theta), \lambda_{\epsilon}(\theta)s)\begin{pmatrix} \frac{\partial a_{\epsilon}}{\partial \epsilon}(\theta) \\ \frac{\partial \lambda}{\partial \epsilon}(\theta) s \end{pmatrix} - \frac{\partial W_{\epsilon}}{\partial \epsilon}(a_{\epsilon}(\theta), \lambda_{\epsilon}(\theta)s) \nonumber \\
  = & - \frac{\partial f_{\epsilon}}{\partial \epsilon}(W_{\epsilon})(\theta, s) \triangleq E_{\epsilon}(\theta, s)
\end{align}

By $\frac{\partial W_{\epsilon}}{\partial \epsilon}(\theta, s) =
DW_{\epsilon}(\theta, s) \eta_{\epsilon}(\theta, s)$ and by differentiating
\eqref{invariance_para}, we end up with
\begin{equation} \label{para_coho}
  \begin{pmatrix} Da_{\epsilon}(\theta) & 0 \\ D \lambda_{\epsilon}(\theta) s & \lambda_{\epsilon}(\theta)
  \end{pmatrix} \eta_{\epsilon}(\theta, s) - \begin{pmatrix} \frac{\partial a_{\epsilon}}{\partial \epsilon}(\theta) \\
\frac{\partial \lambda_{\epsilon}}{\partial \epsilon}(\theta) s
  \end{pmatrix} - \eta_{\epsilon}(a(\theta), \lambda(\theta)s) = E_{\epsilon}(\theta,
s),
\end{equation}
By solving equation \eqref{para_coho} with the same technique used in
Section~\ref{solve_coho}, we can have valid $\frac{\partial
  W_{\epsilon}}{\partial \epsilon}(\theta, s)$, $\frac{\partial
  a_{\epsilon}}{\partial \epsilon}(\theta)$ such that
\eqref{first_order_approximation} admits the first-order starting approximation
for equation \eqref{invariance_para_h}.

\begin{remark}
 We have found, as in \cite{H16}, that the first-order continuation does not
  produce a significant improvement on the initial approximation than the 0-th
  order continuation. A possible reason is that, when the solutions have small derivatives,
 the quadratic convergence of the quasi-Newton method 
is so fast that a better initial approximation does not improve much; on the other 
hand, when the solutions have large derivatives, the extrapolation does not work well. 
 Therefore,  we mainly used the
  0-th order continuation in our implementation.

\end{remark}

\subsubsection{The Continuation Algorithm}

In our implementation, it is important to choose adaptatively the stepsize $h$ in
Section~\ref{continuation_basic}.

Inspired by \cite{CC10}, Algorithm~\ref{algorithm_continuation} here is the 
continuation algorithm we used for the family of
maps $f_\epsilon$ indexed by the one dimensional 
parameter $\epsilon$.  In the numerical examples, we will consider 
\eqref{DST} which depends on several parameters $k,\gamma, \eta$.
We will choose functions $k(\epsilon), \gamma(\epsilon), \eta(\epsilon)$
and then apply the continuation algorithm in $\epsilon$.

\begin{algorithm}[!htb]
  \caption{Continuation algorithm before the breakdown}
  \label{algorithm_continuation}
  \begin{algorithmic}[1]
    \Require{$W_{\epsilon_0}(\theta, s), a_{\epsilon_0}(\theta),
      \lambda_{\epsilon_0}(\theta)$ for the integrable case}
    \Require{$\Delta \epsilon$: The initial increment of the parameter
      $\epsilon$}
    \While{both $\Delta \epsilon$ and $\left\| (W, a, \lambda) \right\|$ are
      acceptable}
    \If {the Algorithm~\ref{algorithm} does not converge}
    \State Move back to the solution $(W, a, \lambda)$ before the increment of the parameters,
    \State Decrease the increment: $\Delta \epsilon$,
    \Else
    \State Update the newly computed $(W, a, \lambda)$,
    \If{$\left\| f_{\epsilon} \circ W_{\epsilon} -
        W_{\epsilon}(a_{\epsilon}, \lambda_{\epsilon}s) \right\| > tolerance$}
    \State Move back to the solution $(W, a, \lambda)$ before the increment of the parameters,
    \State Decrease the increment: $\Delta \epsilon$,
    \EndIf
    \If{$\left\| (W, a, \lambda) \right\|$ exceeds a certain value}
    \State Double the size of the grid points,
    \EndIf
    \EndIf
    \EndWhile
  \end{algorithmic}
\end{algorithm}

\begin{remark}
  As we will discuss in Section~\ref{sec_breakdown}, the invariant circle loses regularity
  when the perturbation is large. Thus, in
  Algorithm~\ref{algorithm_continuation}, one also needs to decrease the order of the
  spline accordingly.
\end{remark}

\subsubsection{Adaptative Grid for Continuation}
It is essential that one check for the correctness of the solution before
doing continuations. As discussed in Section~\ref{subsec_validation}, we need to find the
appropriate grid under which the function has a better representation. This can
be done by performing the line search for different grid sizes and find the one
with the minimum error while applying our quasi-Newton algorithm.

\begin{remark}
  An aspect that requires extra caution is that the existence theorem in
  \cite{YaoL21a} does not guarantee the local uniqueness of the solution. Indeed,
  the solution is only unique under conjugacy (Remark~\ref{undertermincy}).
  This
  results the drift of the solution as the parameter changes. So comparing
  solutions computed through different continuation paths are
  difficult to compare.
\end{remark}

\section{Numerical Explorations} \label{sec_example_2d}
In this section, we take the dissipative standard map as an example to run
the algorithm and explore some of the properties.

The dissipative standard map is a family of 
maps $f_{\eta, \gamma, k}: \mathbb{T} \times \mathbb{R} \rightarrow \mathbb{T}
\times \mathbb{R}$ such that

\begin{equation}
  \label{DST}
  \begin{pmatrix}\theta_{n + 1} \\ p_{n + 1}\end{pmatrix} \triangleq f_{\eta, \gamma,
    k}(\theta_n, p_n) = \begin{pmatrix} \theta_n + p_{n + 1} + \eta, \\ \gamma p_{n}  + \gamma k V'(\theta_n)\end{pmatrix},
\end{equation}

with $(\theta_n, s_n) \in \mathbb{T} \times \mathbb{R} $, $\gamma \in (0, 1)$ is
the dissipative parameter, $k > 0$ is the perturbation
parameter, $\eta$ is the drift parameter, and $V(\theta)$ is an analytic,
periodic function representing the kick from the kicked rotator. In this
example, we shall consider the case when $V(\theta) = - \frac{1}{(2 \pi)^2} \cos(2
\pi \theta)$, then $V'(\theta) = \frac{1}{2 \pi} \sin(2 \pi \theta)$.

\begin{remark}
  If $\gamma = 1$ and $\eta = 0$, the map \eqref{DST} 
  reduces to the Chirikov standard map.
\end{remark}

\begin{remark}
  \label{kzero}
  One can easily verify that when $k = 0$, the solution to the invariance
  equation \eqref{invariance} for $f_{\eta, \gamma, k}$  is
  $$W_{\eta, \gamma, k}(\theta, s) = \begin{pmatrix} \theta \\
    \frac{\gamma}{\gamma - 1}s \end{pmatrix},
  \text{ } a_{\eta, \gamma, k}(\theta) = \theta + \eta, \text{ } \lambda_{\eta,
    \gamma, k}(\theta) = \gamma.
  $$
\end{remark}

For a given choice of the parameters, we follow our quasi-Newton
Algorithm~\ref{algorithm} to
solve the invariance equation \eqref{invariance} for $W_{\eta, \gamma,k}$,
$a_{\eta, \gamma, k}$ and $\lambda_{\eta, \gamma,k}$. For the change of
parameters, we can start with the unperturbed case ($k = 0$) and follow the
continuation algorithm Algorithm~\ref{algorithm_continuation}.

In this example, we first discuss the behavior of the algorithm regarding the
aspects pointed out in Section~\ref{sec_implementation}. First, we
run the algorithmfffffffffffffffffffffffffffffffffffffffffff for a fixed choice of parameters, then we implement the
continuation algorithm through 
three continuation paths in parameter space. 
\begin{enumerate}
\item Continuation with respect to the perturbation parameter $k$, with fixed
  $\gamma$ and $\eta$;
\item Continuation with respect to the drift parameter $\eta$, with fixed
  $k$ and $\gamma$;
\item Continuation with fixed rotation number for the internal dynamics
  $a_{\eta, \gamma, k}$, with fixed $k,$ and $\gamma$, with $\eta$ tuned
  to ensure the inner dynamics is conjugate to a fixed rotation number.
  
  Notice that in this exploration, the parameter $\eta$ is also one of
  the unknowns that have to be determined in the algorithm. This requires
  a small adaptation of the algorithm presented before. 
\end{enumerate}

\begin{remark}
  In practice, we choose the step size and grid size for the continuation of the
  parameters
  dynamically according to Algorithm~\ref{algorithm_continuation}.
\end{remark}

\subsection{Example Solution}
In this subsection, we will show the numerical performace of
Algorithm~\ref{algorithm} regarding the aspects discussed in
Section~\ref{sec_implementation}. In the following examples, we set $\gamma =
0.5$, $\eta = 0.3$.

\subsubsection{Convergence of the quasi-Newton Iteration}
\label{subsec_convergence_rate}
To demonstrate the convergence rate, we set $k = 0.3$ and use the solution when
$k = 0$ (Remark~\ref{kzero}) as the initial approximation, and iterate the algorithm several times
until the error $e(\theta, s)$ for the invariance equation meets the tolerance
(see Table~\ref{table_convergence_rate}). In this example, $N = 1024$, $L = 10$
and $\delta = 0.001$.

\begin{table}[!htb]
  \centering
  \begin{tabular} {|c |c |c |c|}    
    \hline
    Number of Iteration & $\left\| e \right\|_{\mathcal{X}^{0, \delta}}$ & $\left\| e \right\|_{\mathcal{X}^{1, \delta}}$ & $\left\| e \right\|_{\mathcal{X}^{2, \delta}}$ \\
    \hline
    1 & 9.710402e-03 & 6.101232e-02 & 3.833525e-01\\
    \hline
    2 & 2.860761e-04 & 3.274913e-03 & 5.806850e-02\\
    \hline
    3 & 5.587798e-06 & 9.021029e-05 & 2.031524e-03\\
    \hline
    4 & 4.152389e-10 & 7.825378e-09 & 2.376669e-07\\
    \hline
    5 & 2.645506e-14 & 9.540554e-13 & 2.810590e-09\\
    \hline
    6 & 3.889196e-16 & 3.030427e-13 & 2.737512e-09\\
    \hline
  \end{tabular}
  \caption{Convergence of the quasi-Newton Iteration}
  \label{table_convergence_rate}
\end{table}

\begin{remark}
The choice of $k = 0.3$ and the initial approximation ($k = 0$)
  are only used to demonstrate the convergence of the algorithm. In practice, we
  usually have a much smaller step size (less than $10^{-3}$) and the initial
  error usually is of order $10^{-6}$. The step size gets even much smaller when
  the perturbation is big and when $k$ is close to the breakdown.
\end{remark}

\begin{remark} 
  As discussed in Remark~\ref{remark_cubic_spline}, we mainly used cubic splines
  in our implementation. 
  Remark~\ref{rmk_spline_error} illustrates the reason for the
  increase of the round-off error for $\left\| \cdot \right\|_{\mathcal{X}^1}$ and $\left\|
    \cdot \right\|_{\mathcal{X}^2}$ in Table~\ref{table_convergence_rate}.
\end{remark}

\subsubsection{Validation for Correctness of the Solution} \label{example_validation}
In this part, we
continue with the solution achieved in 
Section~\ref{subsec_convergence_rate} (when $k = 0.3$). Following the
discussions in Section~\ref{subsec_validation}, we validate the solution
through monitoring both the error when the grid size is doubled and the norm of
the solution. For the grid size equals 1024, the corresponding data are
recorded in Table~\ref{table_validation}. Since the error remains to be
relatively small, we are more confident to say that the solution we achieved is
indeed the true solution.

\begin{table}[!htb]
  \begin{tabular} {|c| c|c|c|}
    \hline
    & $\left\| \cdot \right\|_{\mathcal{X}^{0, \delta}}$ & $\left\| \cdot \right\|_{\mathcal{X}^{1, \delta}}$ & $\left\| \cdot \right\|_{\mathcal{X}^{2, \delta}}$ \\
    \hline
    Error for the invariance equation & 1.622146e-13 & 1.494237e-10 & 3.059292e-06\\
    \hline
    Norm for the solution & 1.07970 & 1.225545 & 1.768379\\
    \hline
  \end{tabular}
  \caption{The performance when the grid size is doubled}
  \label{table_validation}
\end{table}

\begin{remark}
  As indicated in Table~\ref{table_validation}, the norm for the solution is
  relatively small. This happens because the choice of our perturbation is small.
  In fact, as discussed in Section~\ref{subsec_continuation_k}, such norm will
  blow-up when $k$ is close to the breakdown value.
\end{remark}

\subsubsection{Methods for Computing the Inverse of $a(\theta)$}
\label{subsubsection_inverse_method_plot}
\begin{figure}[!htb]
  \includegraphics[width=0.95\textwidth]{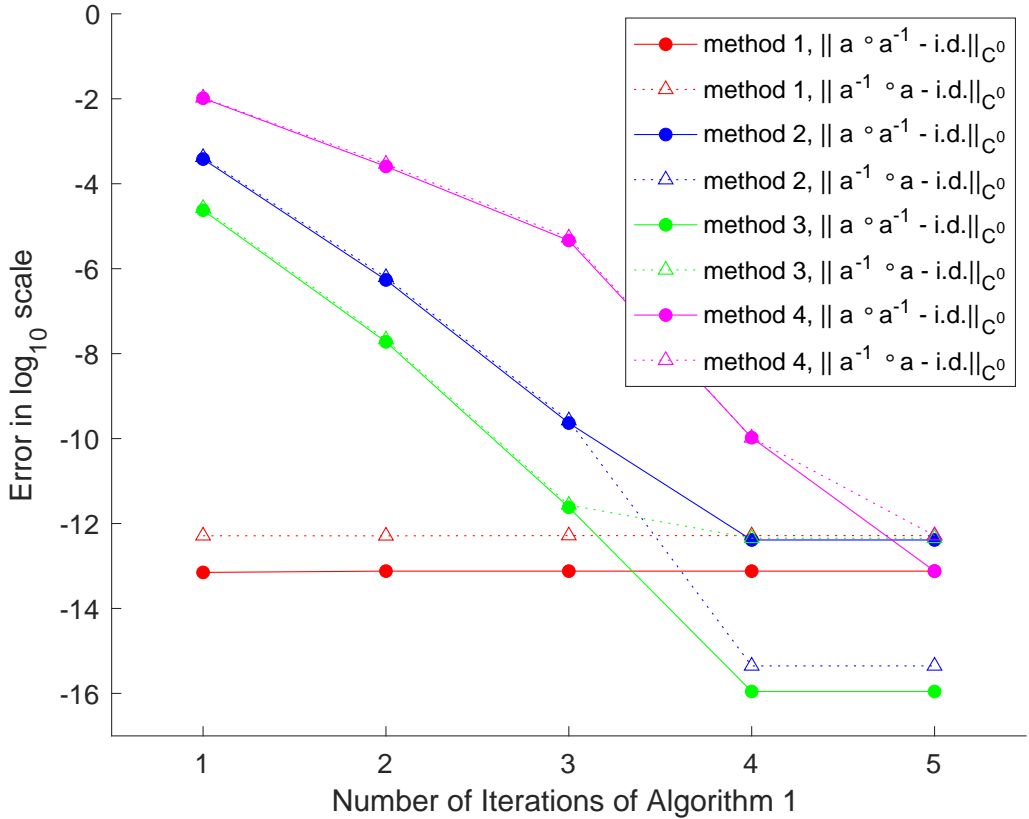}
  \caption{Compariation of the 4 methods in Section~\ref{subsec_inverse_method}}
  \footnotesize
  \emph{Method 1: Graph Reflection; Method 2: ``Left'' Inverse; Method 3:
    ``Right'' Inverse; Method 4: Inverse through Solving Cohomological Equation}
  \label{ainv_a_a_ainv}
\end{figure}

In this paragraph, we present some numerical results comparing 
the methods proposed in 
Section~\ref{subsec_inverse_method} to compute the inverse of 
the dyamics. 
The detailed results 
are presented in Figure~\ref{ainv_a_a_ainv}.

\begin{remark}
  As indicated in Figure~\ref{ainv_a_a_ainv}, Method 1 always has a great
  performance as it is independent of the error of the cohomological equation,
  while Method 2, 3, 4 slowly get better as the convergence of the algorithm.
  Despite that Method 4 is designed for optimizing both $a \circ a^{-1}$ and
  $a^{-1} \circ a$, it turns out Method 4, in general, does not outperform the
  other methods. This is also true when $k$ is bigger or even near the
  breakdown.

  As stated in Section~\ref{subsec_inverse_method}, in practice, we try all of the
  methods and use the one with the best performance. More specifically, through
  out the iterations of the algorithm, we
  always start by method 1, and then replace it with method 2 or 3.
\end{remark}

\subsubsection{Run Time Analysis}
\label{subsub_run_time_experiment}
As discussed in Section~\ref{subsec_algorithm}, both the time and space
complexity for Algorithm~\ref{algorithm} are $\mathcal{O}\big(N \times L\big)$.

For the same choice of parameters as in Section~\ref{subsec_convergence_rate}, the
average time for running one iteration of Algorithm~\ref{algorithm} can be found
in Table~\ref{table_run_time}. The code is written in C using the GNU Scientific
Library (GSL),
and this set of data in Table~\ref{table_run_time} is generated by a Mid
2014 13-inch Macbook Pro with 2.6 GHz Dual-Core Intel Core i5 Processor and 8 GB
1600 MHz DDR3 Memory.

\begin{table}[!htb]
  \centering
  \begin{tabular} {|c|c|c||c|c|c|}
    \hline
    $L$ & $N$ & Avg Time & $L$ & $N$ & Avg Time  \\ 
    \hline \hline
    \multirow{4}{*}{2} & 1024 & 0.062405 & \multirow{4}{*}{5} & 1024 & 0.132773 \\
        & 4096 & 0.284577 & & 4096 & 0.641300 \\
        & 16384 & 1.682731 & & 16384 & 3.674957 \\
        & 65536 & 9.136261 & & 65536 & 18.175884 \\
        & 262144 & 55.717758 & & 262144 & 122.279080 \\
    \hline

    \hline \hline
    \multirow{4}{*}{10} & 1024 & 0.392576 & \multirow{4}{*}{20} & 1024 & 0.779805 \\
        & 4096 & 1.862491 & & 4096 & 4.708662 \\
        & 16384 & 9.303865 & & 16384 & 21.924758 \\
        & 65536 & 44.038226 & & 65536 & 99.157648 \\
        & 262144 & 226.787753 & & 262144 & 635.228799 \\
    \hline
  \end{tabular}
  \caption{Average run time (in seconds) for one iteration of Algorithm~\ref{algorithm}}
  \label{table_run_time}
\end{table}

\subsubsection{Plot of the Invariant Circle and Isochrons}
\label{subsubsec_color_plot}

From Remark~\ref{kzero}, the isochrons are approximately linear in the
neighborhood of the invariant circle when the perturbation is small. To achieve a
relatively nontrivial
solution, we consider $k = 1.1037$ and $\eta = 0.30532$, where $\eta$ is
tuned manually in order to guarantee a rational rotation number for
$a(\theta)$ using the Brent zero-finding algorithm (see Section~\ref{subsec_continuation_k_prev}). As
one will see in Section~\ref{subsubsec_globalization}, such cases require more
considerations in terms of globalizing the isochrons. 
The solution here is computed through
continuation method with step size $10^{-3}$ and error tolerance $10^{-14}$ in
$C^0$.

\begin{figure}[!htb]
  \centering
  \subfloat[$\theta_0 = 0.797$]{
    \includegraphics[width=0.48\textwidth]{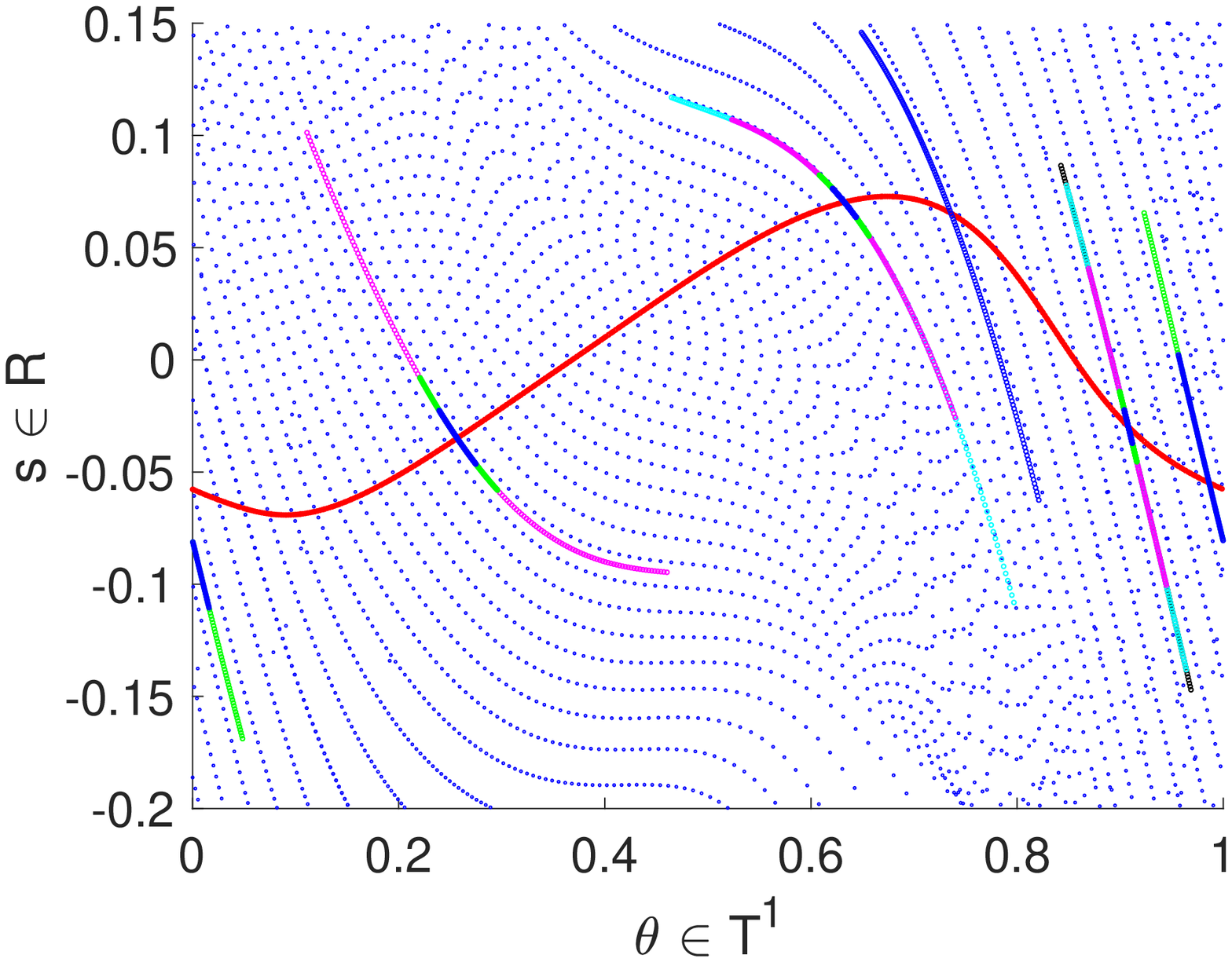}
    \label{fig:cic1}
  }
  \subfloat[$\theta_0 = 0.78$]{
    \includegraphics[width=0.48\textwidth]{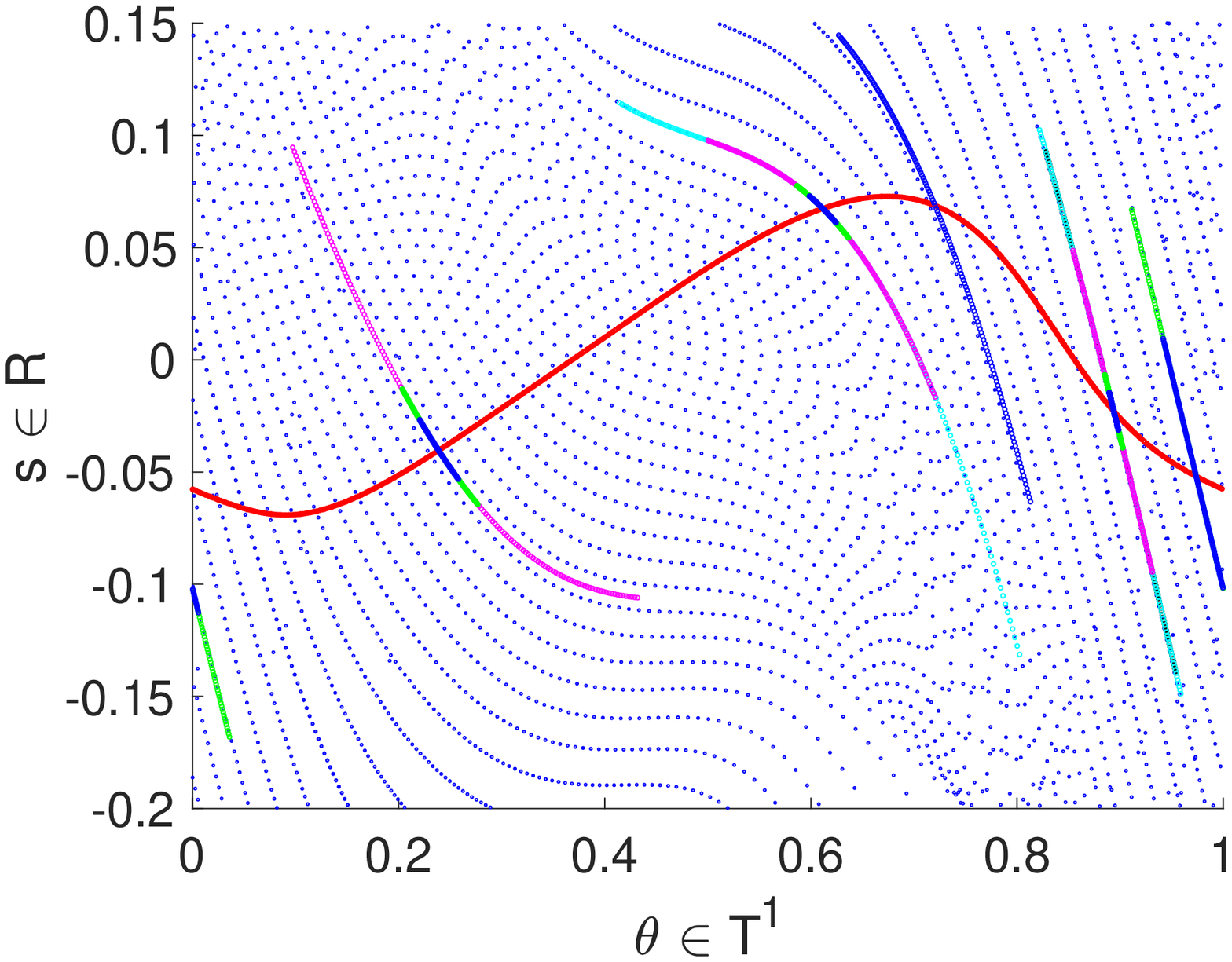}
    \label{fig:cic2}
  }\\
  \footnotesize
  \emph{where $I_{\theta_0}$ corresponds to the blue isochron.}
  \\
  \subfloat[$\lambda$ and $\lambda^*$]{
    \includegraphics[width=0.48\textwidth]{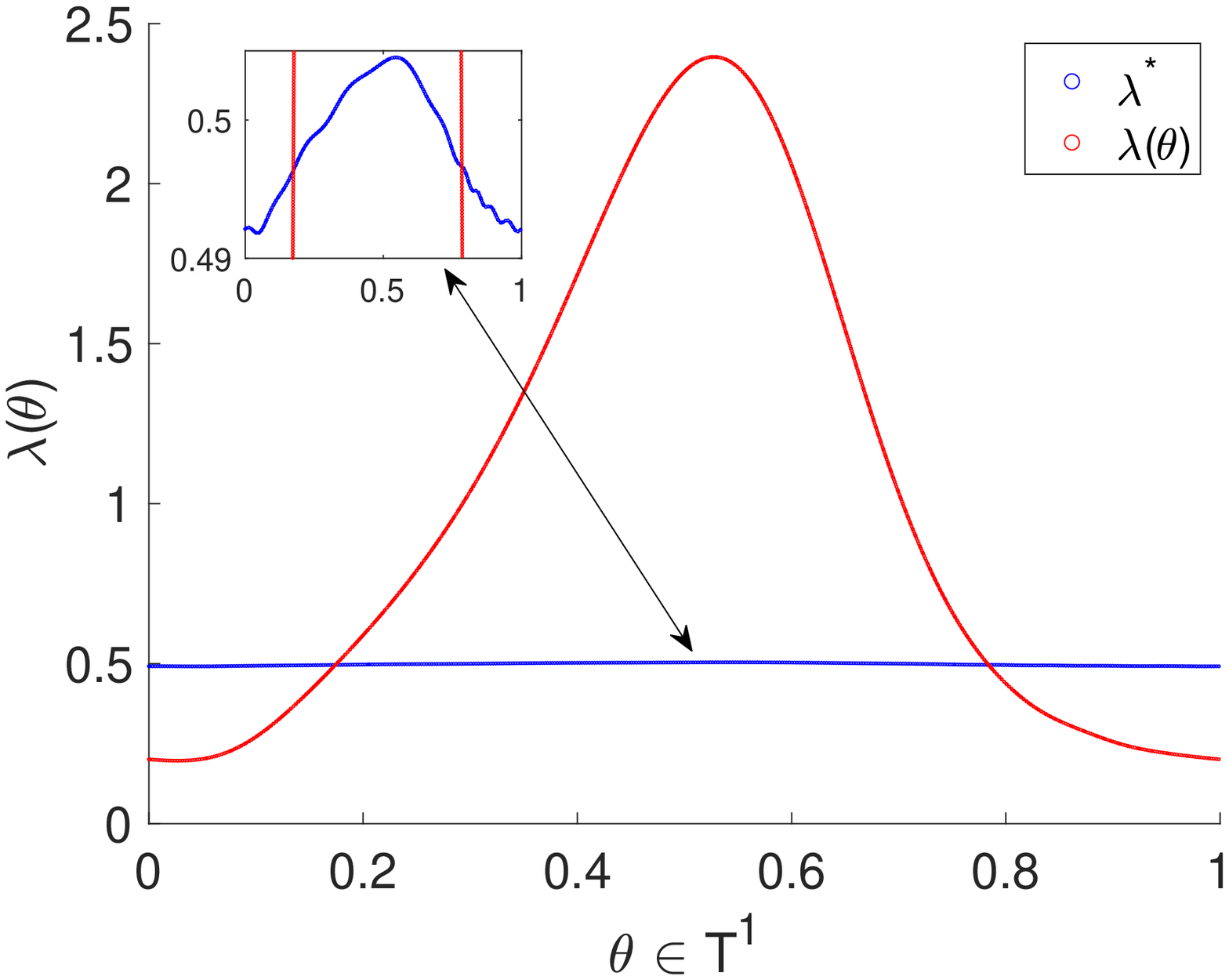}
    \label{fig:cic3}
  }
  \subfloat[rational rotation]{
    \includegraphics[width=0.48\textwidth]{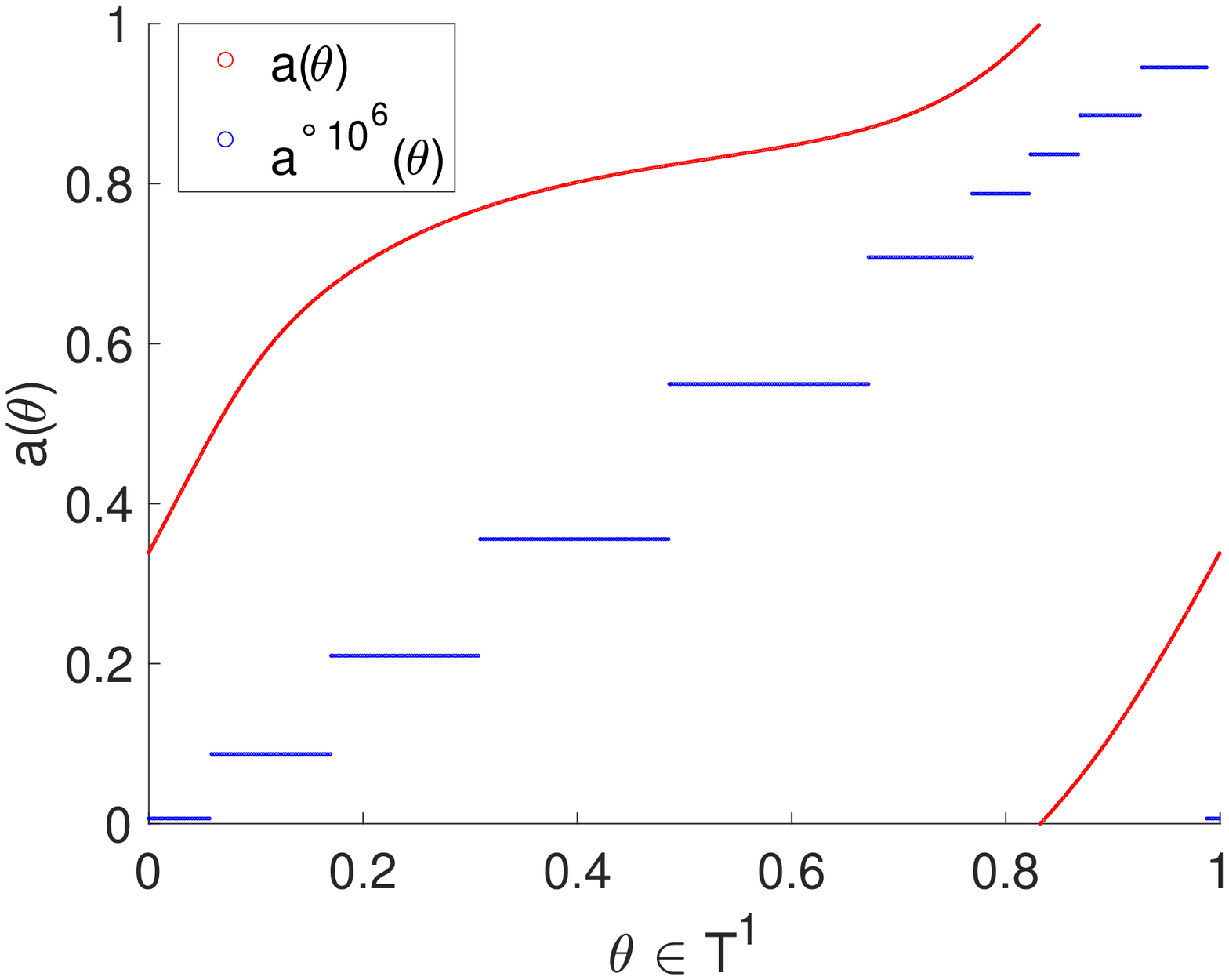}
    \label{fig:cic4}
  }
  \caption{Invariant Manifold and Isochrons for the dissipative standard map
    \eqref{DST}}
  \footnotesize
  \emph{The last figure is only here to indicate that the rotation number for
    $a(\theta)$ is indeed rational. Please refer to
    Section~\ref{subsec_continuation_eta} for further reasons.}
  \label{plot_all}
\end{figure}

In Figure~\ref{fig:cic1}, we present the invariant circle (in red) and some of the
isochrons (in the  order: blue, green, magenta, cyan, black along the
internal dynamics). Notice how the isochrons contract when
the dissipative standard map is applied. More specifically, the blue isochron
($I_{\theta = 0.797}$) got
mapped to the interior of the green isochron, and then the green-blue isochron
got mapped to the interior of the magenta one, and then the cyan isochron
followed by the black one.

With a different starting point on the grid, Figure~\ref{fig:cic2} describes a case
when the isochrons do not contract for every single step (indeed, one may find
that the cyan isochron is actually expanded rather than contracted). This is
perfectly normal: As discussed in Remark~\ref{remark_dynamical_average}, our
algorithm allows $\left\|\lambda \right\|_0$ to be bigger than one, as long as
the dynamical average (discussed in Remark~\ref{remark_dynamical_average})
$\lambda^* < 1$ (as shown Figure~\ref{fig:cic3}). Notice also that $\|\lambda\|_{C^0}$
changes when we perform a change of variable in the map, but $\lambda^*$ does
not. Thus, by the underdetermination
of the solution, one can also find a suitable
$\lambda$ with a  contraction that do not depend on the point. 

\begin{remark}
  In the case when $\eta = 0.3$, while the other parameters remain unchanged,
  one can observe that the dynamical average is in fact not a constant
  function, on the contrary to the irrational rotation case (for the same reason
  as in Remark~\ref{remark_constant_lambda}).
\end{remark}

\subsubsection{Globalization}
\label{subsubsec_globalization}

\begin{figure}[
  !htb]
  \subfloat[Globalization through Backward Propagation]{
    \includegraphics[width=0.49\textwidth]{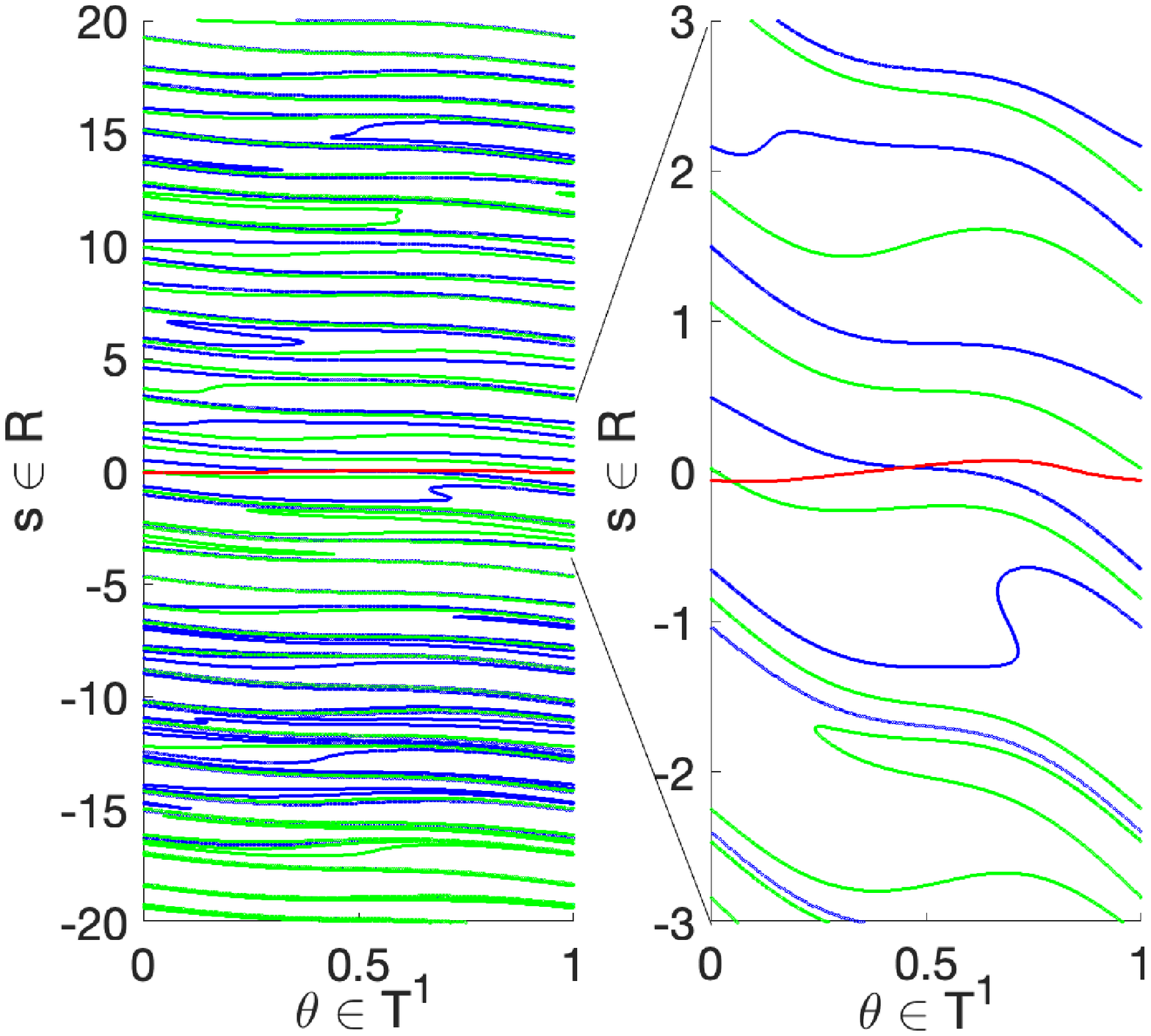}
    \label{fig:glo1}
  }
  \subfloat[Globalization through Increasing Approximation Order]{
    \includegraphics[width=0.49\textwidth]{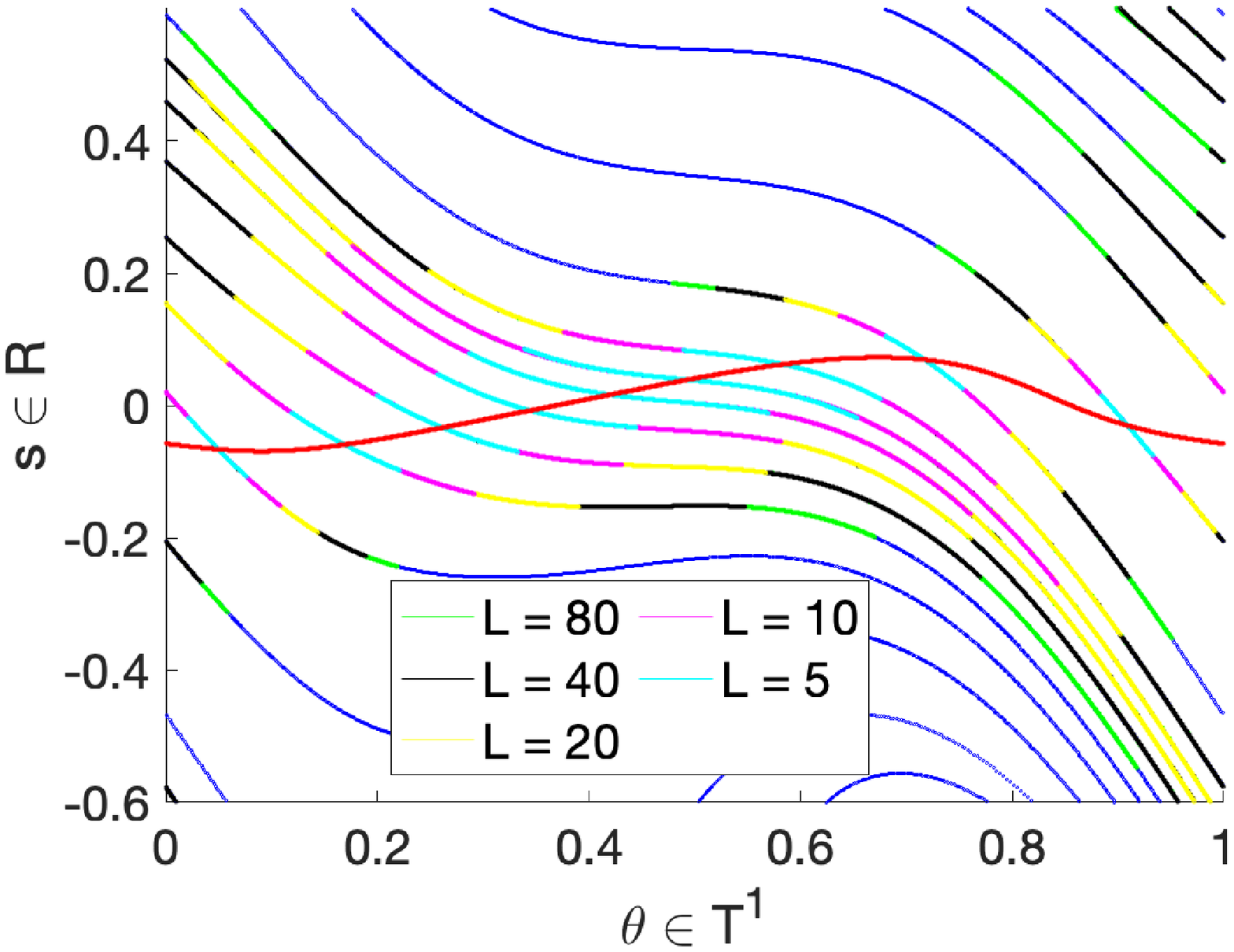}
    \label{fig:glo2}
  }
  \label{fig:globalization_inv}
  \caption{Globalization of Isochrons}
\end{figure}

\begin{figure}[!htb]
  \includegraphics[width=0.99\textwidth]{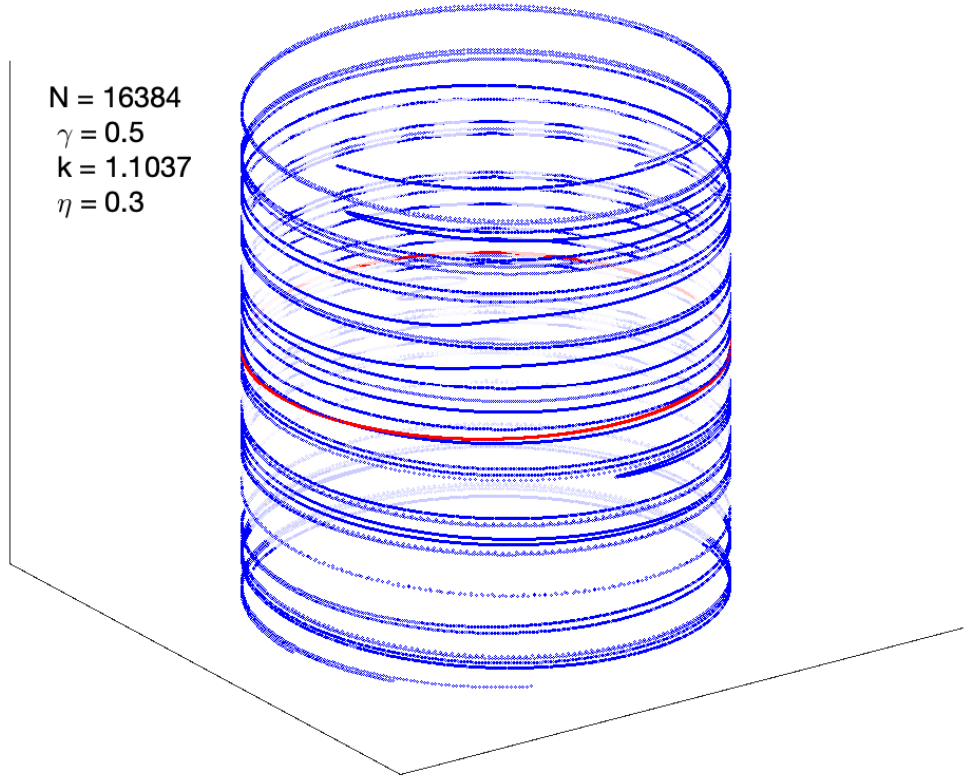}
  \label{fig:best_fig}
  \caption{Invariant Circles and the Corresponding Isochrons in 3D}
\end{figure}

Due to the truncations with respect to $s$,  the approximate solution 
 $W(\theta, s) = \sum_{i = 0}^LW^{(i)}(\theta)s^i$ (a Taylor approximation of 
the parameterization)  only 
approximates the isochrons for a small range of $s$.  To obtain a representation 
of the isochrons in a larger domain, we perform \emph{globalization}.

Basically, we use the functional equation satisfied by the isochrons and obtain 
representations in larger domain.  The extended isochrons are the iterates under 
the inverse map of the isochrons accurately computed (in a small enough domain of $s$). 
Even if the isochrons will be very regular in  neighborhood of the invariant circle, 
when we globalize them, they will have interesting interactions with other invariant objects in 
the map (e.g. if there are several attractors, the isochrons will accumulate on 
the boundary of the basin of attraction of the limit cycle.)

Similar procedures have been used extensively. In this case, there are some peculiarities.

\begin{remark}
  When the internal dynamics $a(\theta)$ has rational rotation number (thus so
  does $a^{-1}(\theta)$), the
  phase-locking phenomenon occurs. In this case, the globalized isochrons will
  accumulate near the periodic orbit instead of distributing in the whole
  $\mathbb{T}$.  

  To resolve this, one needs to initialize the isochrons in 
a grid which has many points near the stable point. 
 
A practical algrorithm is to start with an even grid on $\mathbb{T}$, first apply
  $a(\theta)$ $n$ times, compute the stable manifolds on the resulting points and 
 then do the backward propagation
  of the map $n$ times on the computed stable manifolds. 
\end{remark}

Figure~\ref{fig:glo1} gives the result of  the globalization of two isochrons (in blue and
green) using the above method.

Another way to enlarge the validity  region of $s$ is via increasing
$L$, the order of truncation in $s$.
Again, we can use the solutions for an order as 
approximate solutions for a truncation of higher order and apply the Newton 
method. See Algorithm~\ref{algorithm_nothing}.
Figure~\ref{fig:glo2} presents the isochrons
up to order 80.

\begin{remark}
Since \eqref{DST}
  is entire,  an argument in   \cite{YaoL21a}
shows that $W(\theta, s)$ is entire in $s$ for every $\theta \in \mathbb{T}$. 

If the map $f$ was not entire  (an important case is
when the maps appear as the time-one of 
a perturbation of an ODE), then the radius of convergence of the expansion 
in $s$ of $W(\theta, s)$ could be finite and, in such a
case, increasing the order $L$ would not improve the domain of validity. 

Even if the function is entire, the numerical computation could be 
affected by round-off errors if the coefficients of the expansion have very different sizes. 
We have found it convenient to introduce a scale factor $b$ so that we consider 
the expansions on $b s$ rather than on $s$. By adjusting $b$ we can get calculations
that, even if mathematically equivalent to the original ones, are less affected 
by round-off effects. 
\end{remark}

\begin{algorithm}[!htb]
  \caption{Increasing the Order of Parameterization}
  \label{algorithm_nothing}
  \begin{algorithmic}[1]
    \Require{Solution for Equation \eqref{invariance} $W(\theta, s), a(\theta)$
      and $\lambda(\theta)$ with
      truncation order $L_{start}$}
    \Require{Tolerance for the error of Equation \eqref{invariance} in
      $\mathcal{X}^{\delta_{start}, r}$ norm, with some precribed $\delta_{start}$}
    \Ensure{Solution $W(\theta, s), a(\theta)$ and $\lambda(\theta)$ to the
    invariance equation \eqref{invariance} with truncation order $L_{end}$}
  \Let{$\lambda(\theta)$}{$\lambda(\theta) + \Delta_{\lambda}(\theta)$},
  \For{$l \gets L_{start}, L_{end}$}
  \State Find the biggest $\delta$ (using any root finding method: Brent,
  Bisect, etc.) that essures the $\mathcal{X}^{\delta, r}$
  norm of the error for Equation \eqref{invariance} is within the tolerance.
  \EndFor
  \State Return updated $W(\theta, s), a(\theta)$ and $\lambda(\theta)$.

\end{algorithmic}
\end{algorithm}

We present Figure~\ref{fig:glo2} with the invariant circle (in
red) and 10 isochrons (in blue), and as indicated in Figure~\ref{fig:cic4}, the
internal dynamics have rational rotation number.

\subsection{Continuation w.r.t. $k$}
\label{subsec_continuation_k}

In this subsection, we
perform the continuation scheme as in
Algorithm~\ref{algorithm_continuation} for different $k$, with fixed $\gamma$
and $\eta$.

We start at $k = 0$, with Remark~\ref{kzero} as the initial approximation, and
keep computing until the
quasi-Newton algorithm, Algorithm~\ref{algorithm}, stops. 
The invariant circle and the corresponding stable manifolds are demonstrated in
Figure~\ref{plot_all}.

\begin{figure}[!htb]
  \centering
  \subfloat[$k = 0.0000$]{
    \includegraphics[width=0.313\textwidth]{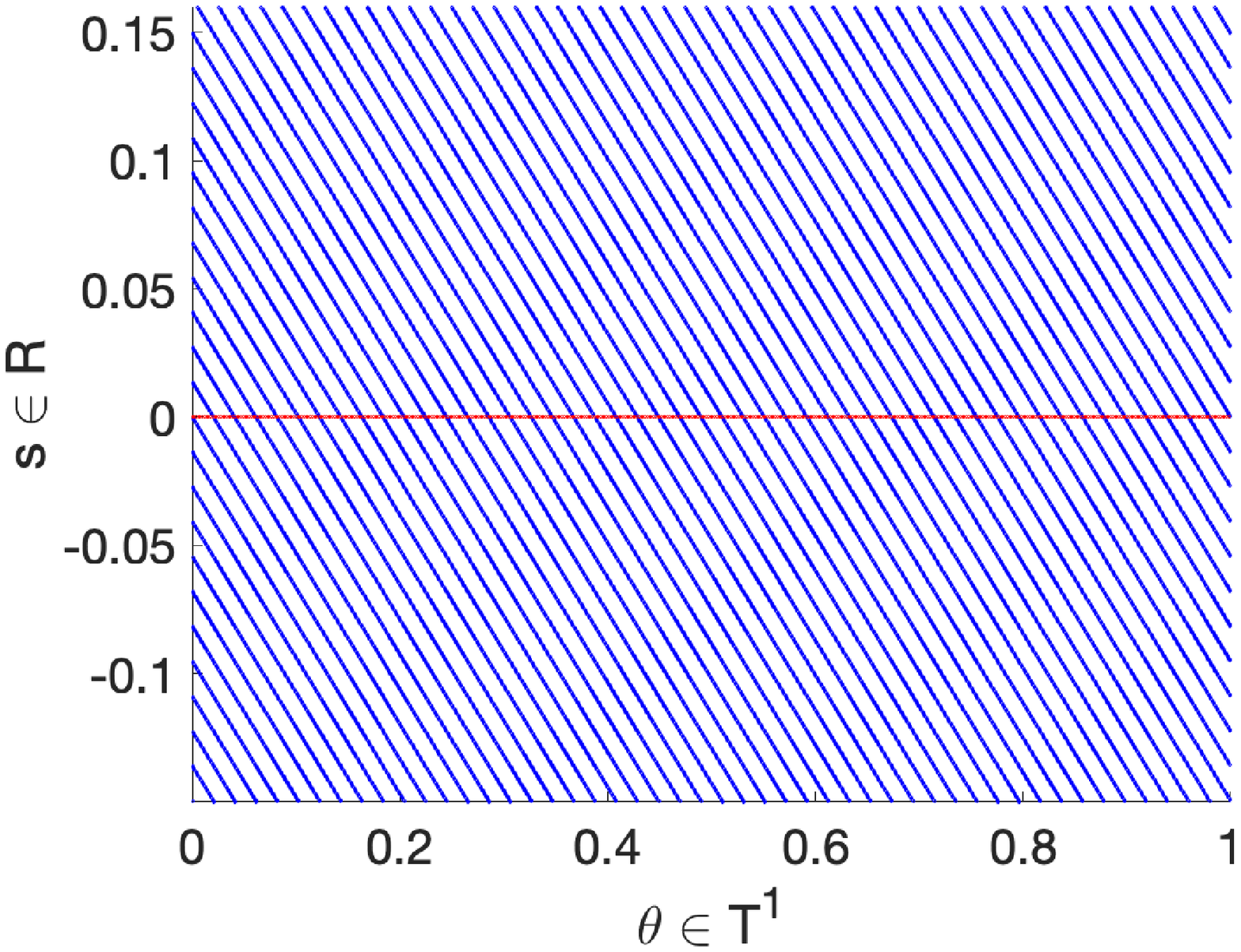}
    \label{fig:f1}
  }
  \subfloat[$k = 0.2000$]{
    \includegraphics[width=0.313\textwidth]{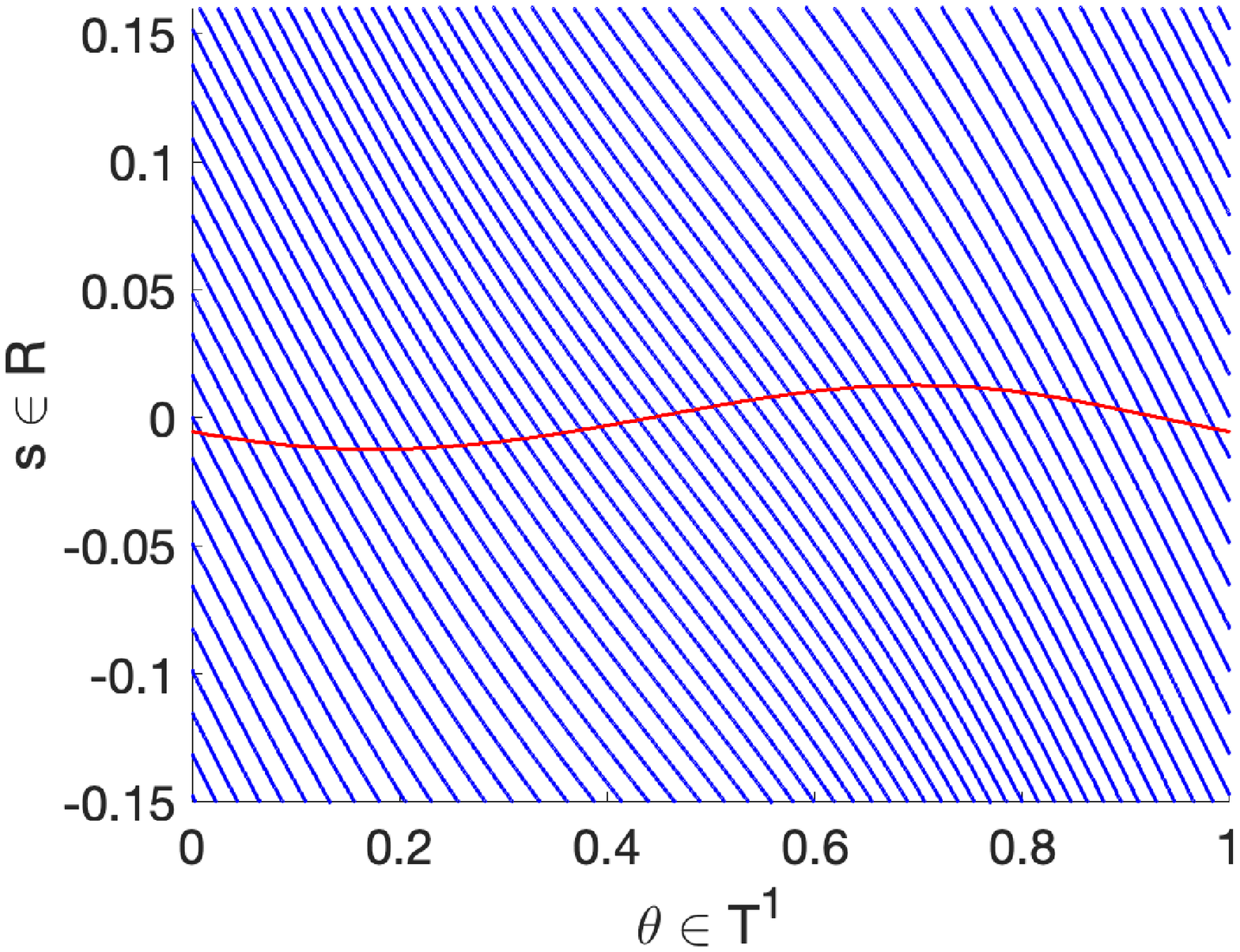}
    \label{fig:f2}
  }
  \hfill
  \subfloat[$k = 0.4000$]{
    \includegraphics[width=0.313\textwidth]{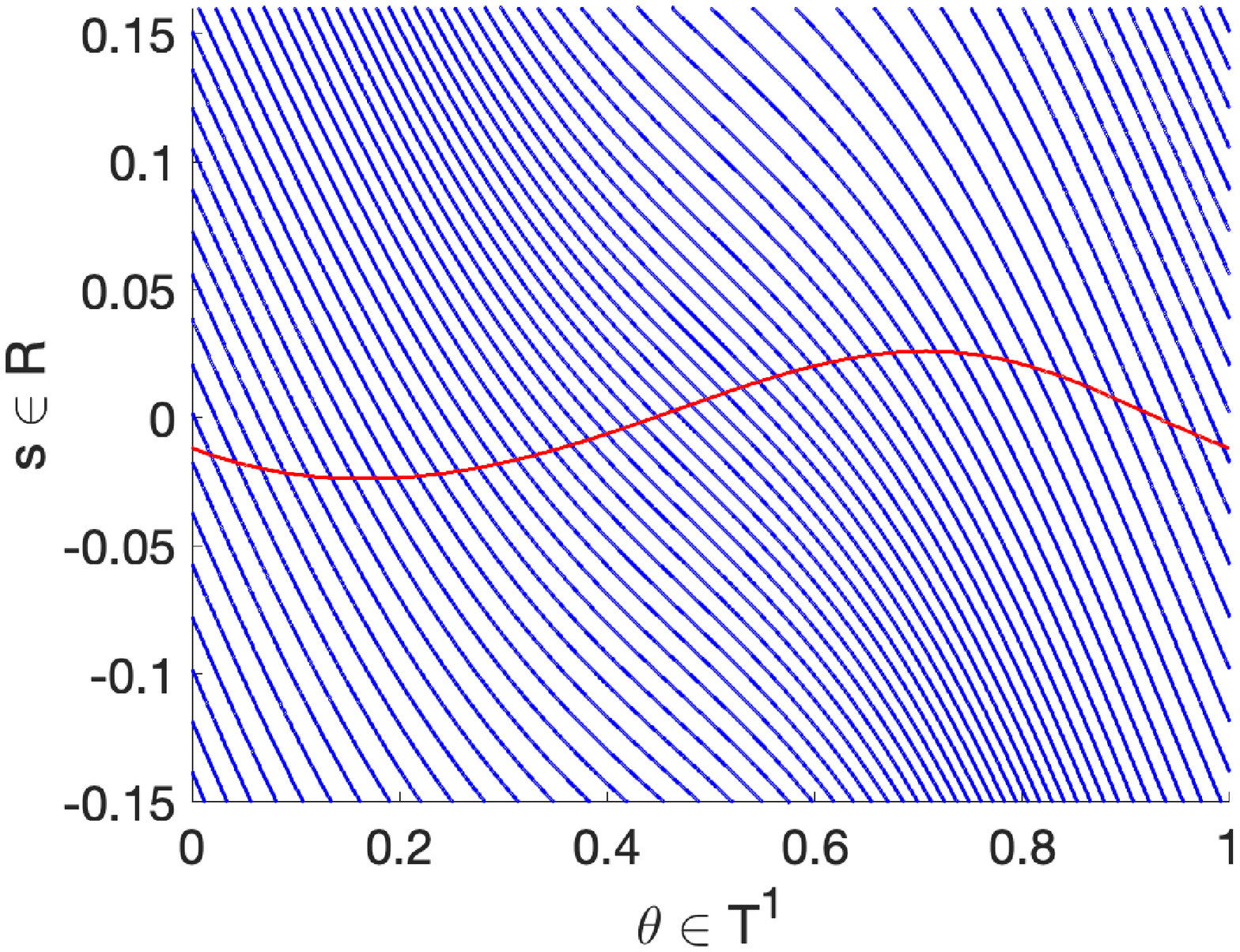}
    \label{fig:f3}
  }\\
  \subfloat[$k = 0.6000$]{
    \includegraphics[width=0.313\textwidth]{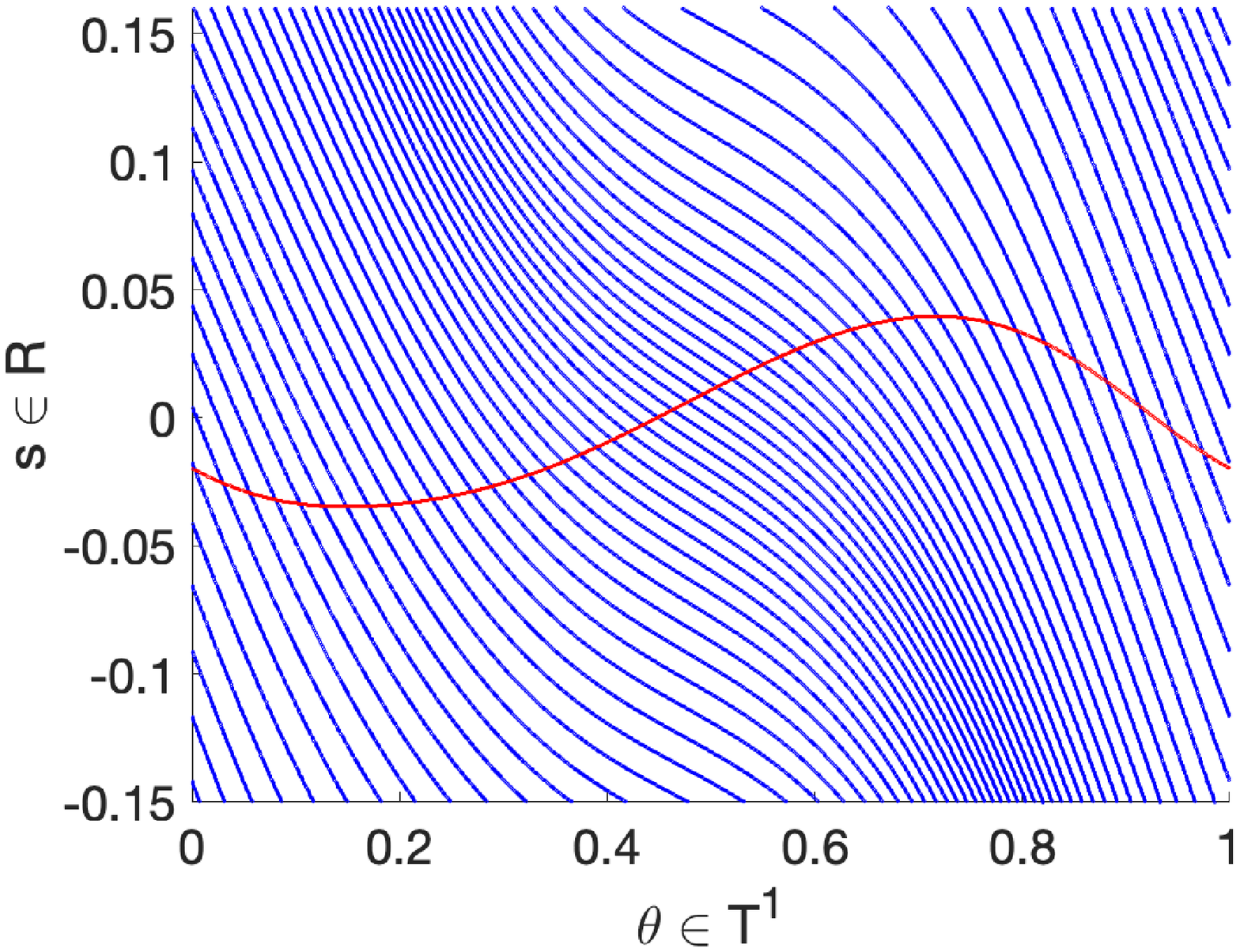}
    \label{fig:f4}
  }
  \subfloat[$k = 0.8000$]{
    \includegraphics[width=0.313\textwidth]{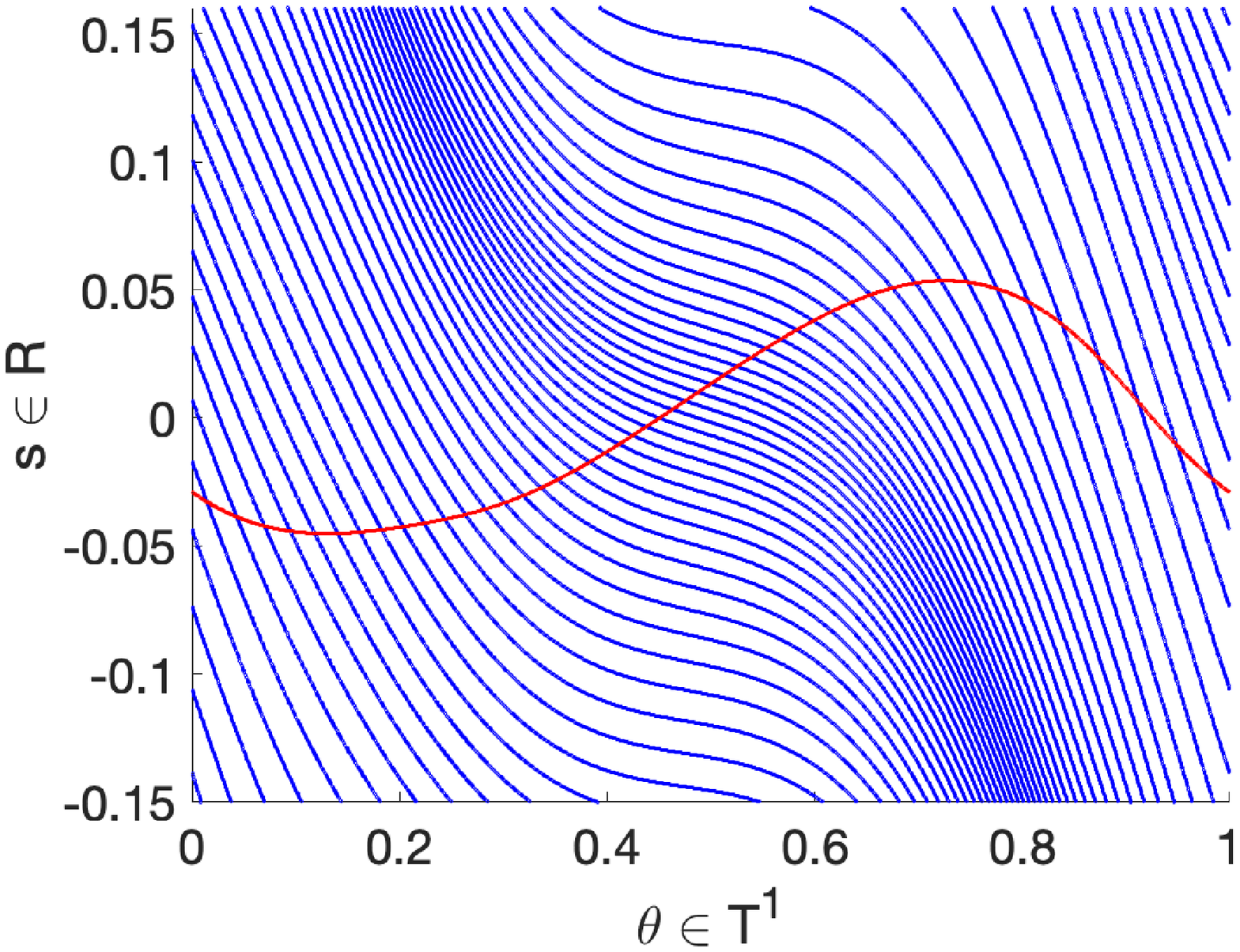}
    \label{fig:f5}
  }
  \hfill
  \subfloat[$k = 1.0000$]{
    \includegraphics[width=0.313\textwidth]{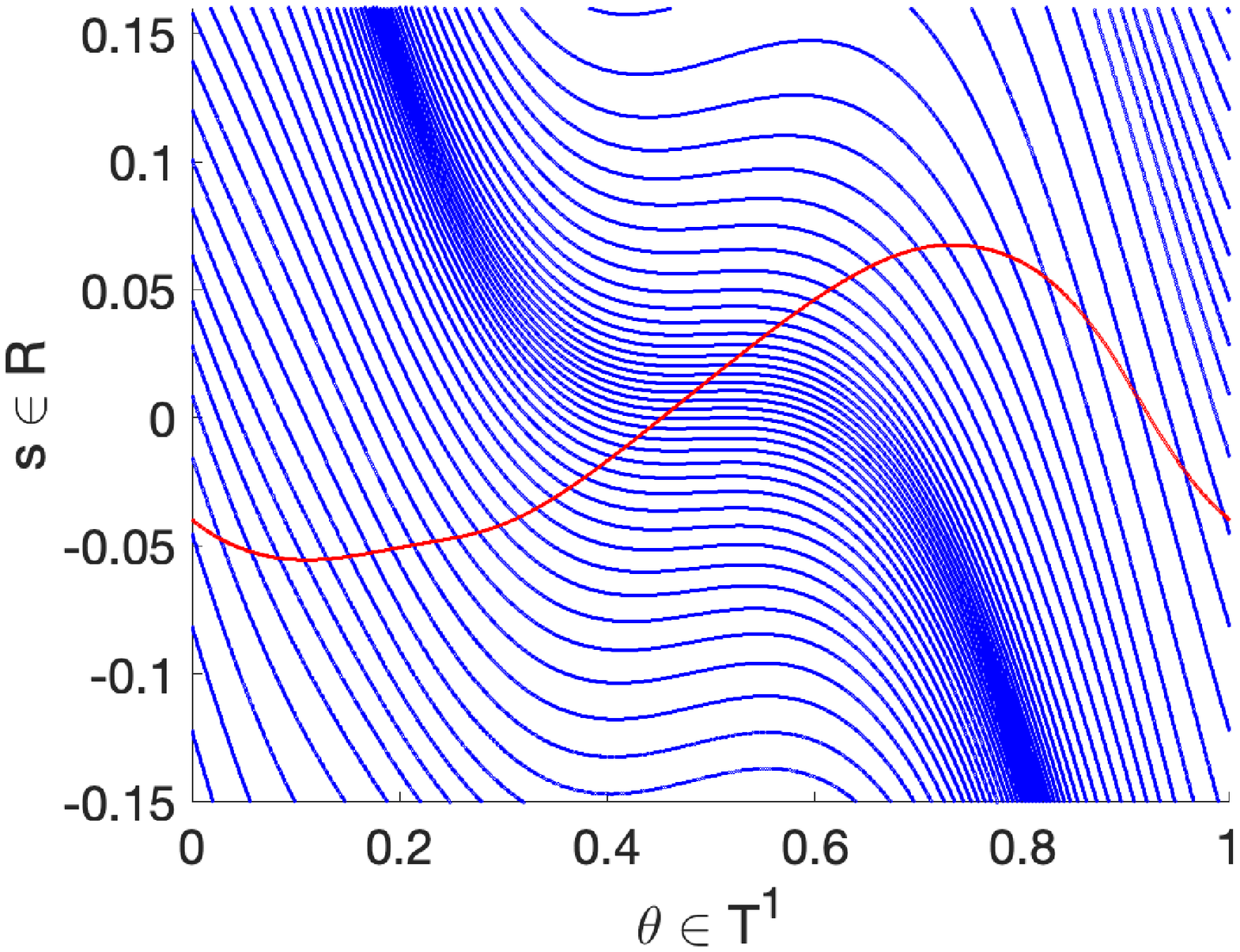}
    \label{fig:f6}
  }\\
  \subfloat[$k = 1.2000$]{
    \includegraphics[width=0.313\textwidth]{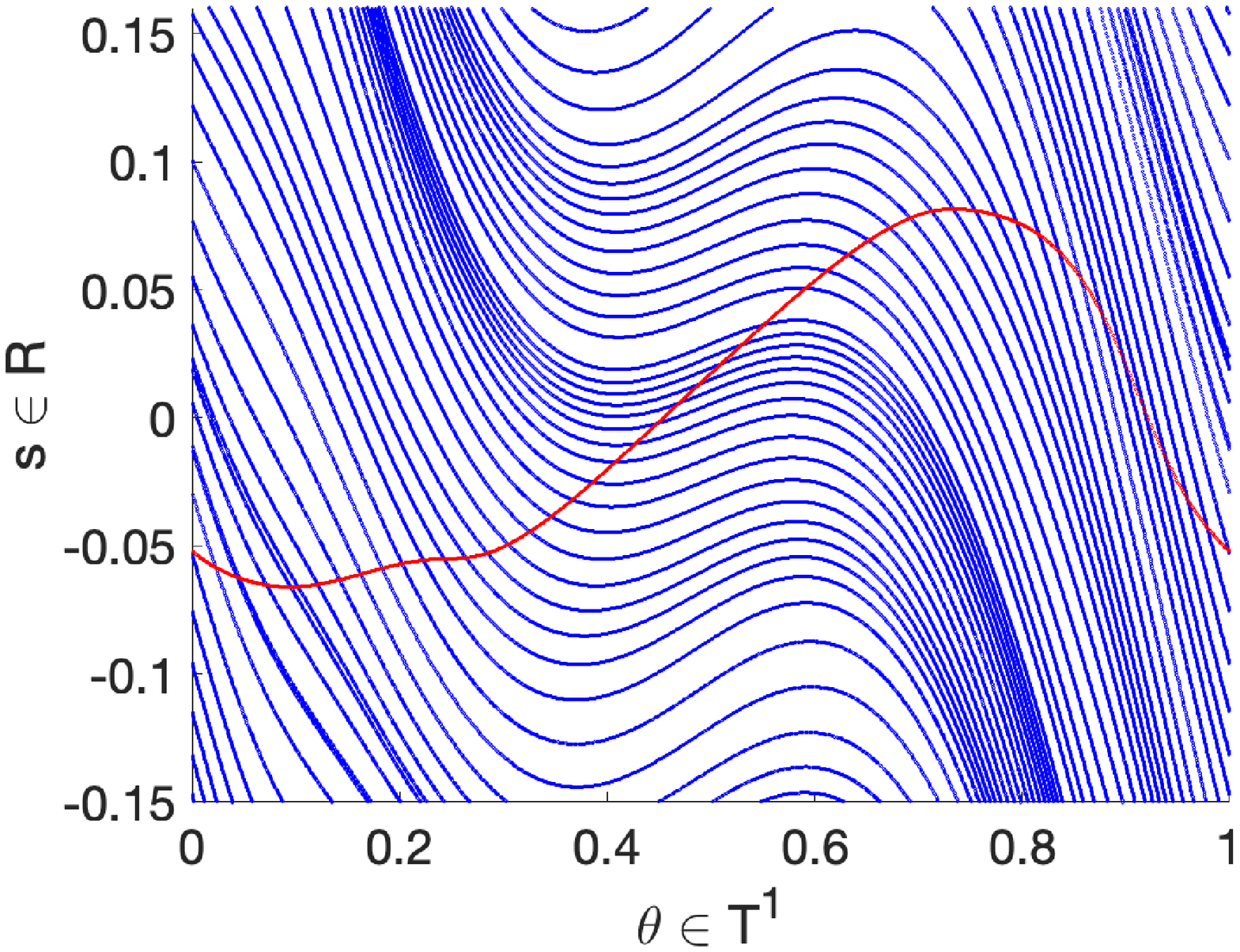}
    \label{fig:f4}
  }
  \subfloat[$k = 1.4000$]{
    \includegraphics[width=0.313\textwidth]{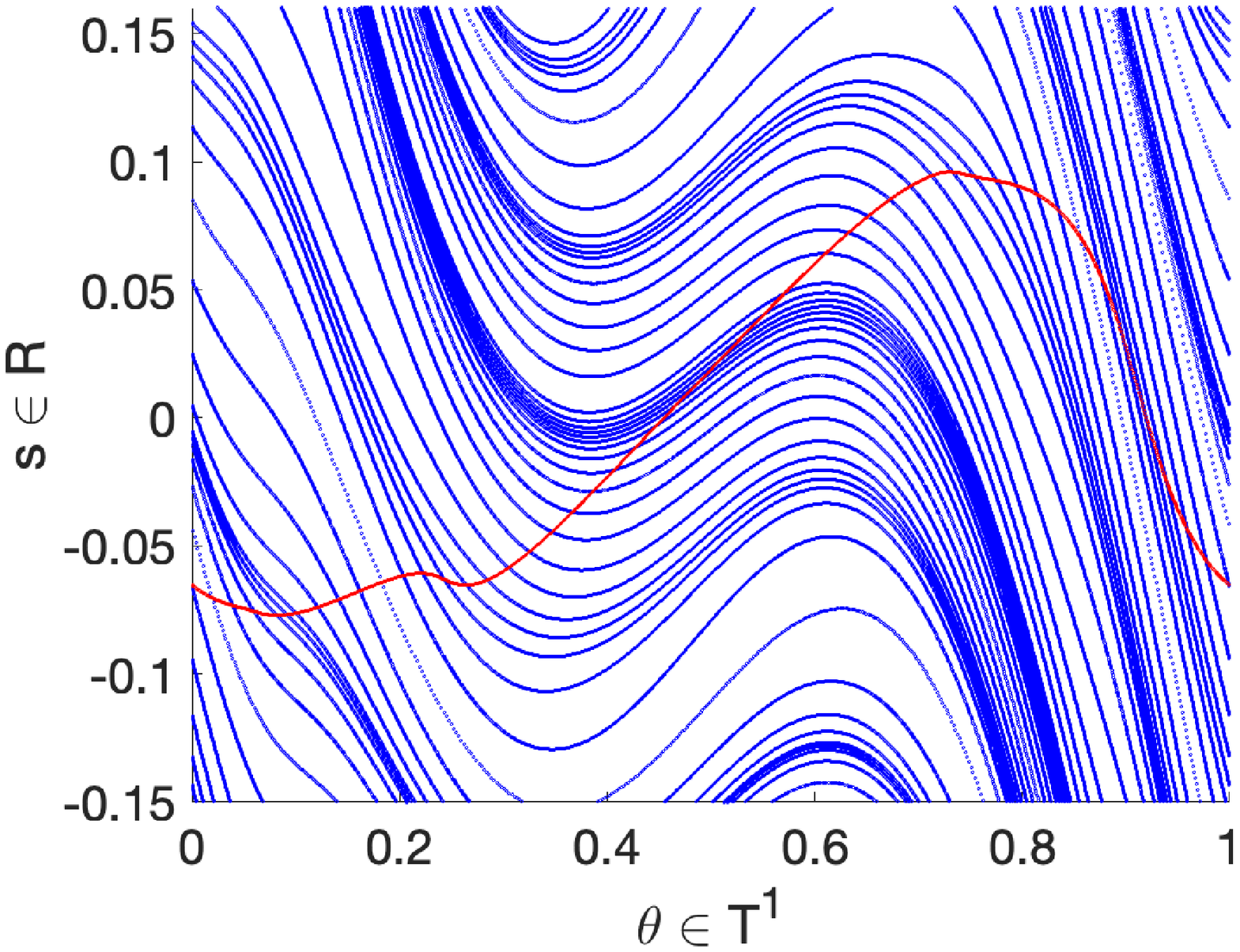}
    \label{fig:f5}
  }
  \hfill
  \subfloat[$k = 1.4927$]{
    \includegraphics[width=0.313\textwidth]{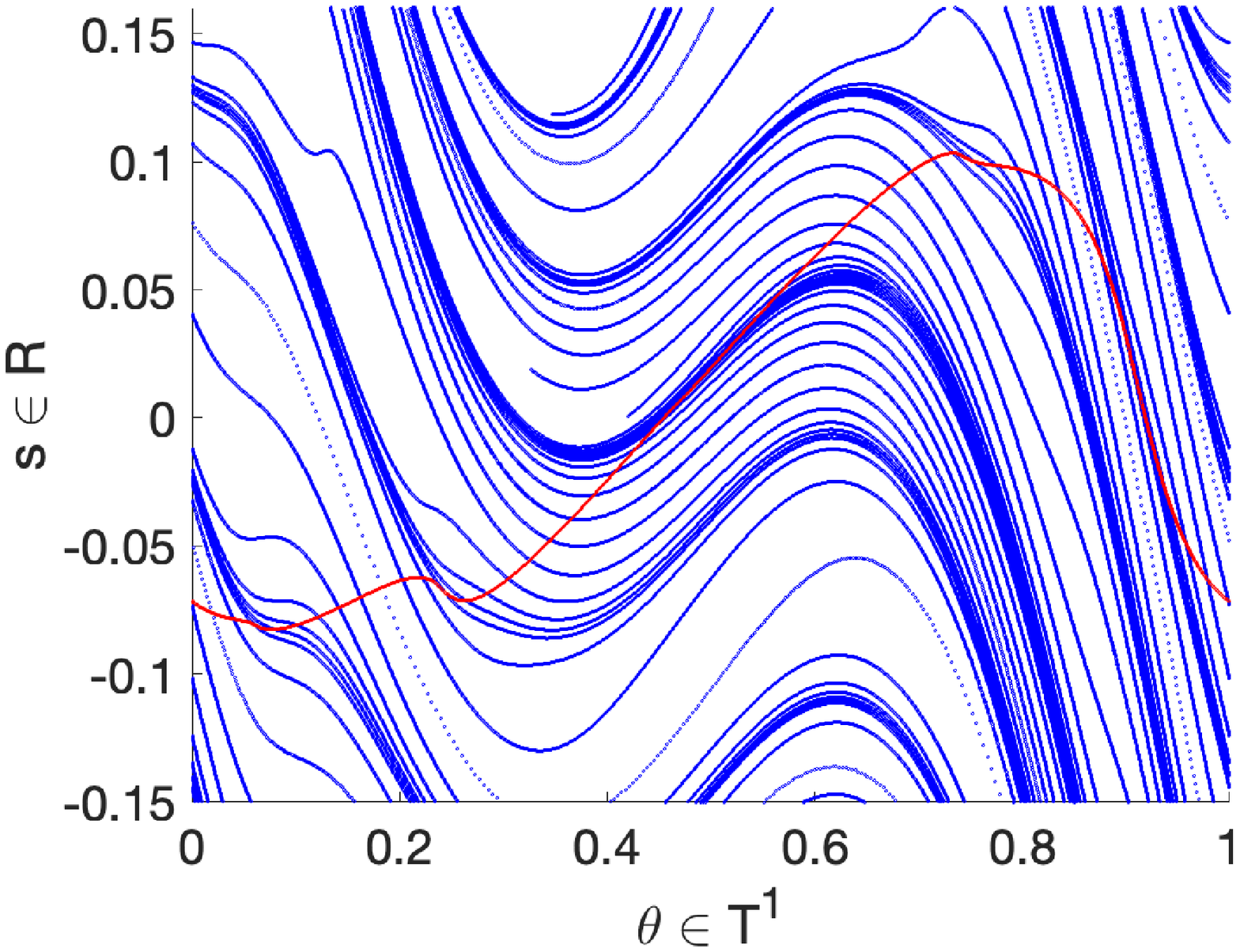}
    \label{fig:f6}
  }
  \caption{Invariant Manifold and Isochrons for the dissipative standard map
    \eqref{DST} where $\gamma = 0.6$, $\eta = 0.4$ and $k$ varies from 0 to
    1.4927}
  \label{plot_all}
\end{figure}

As the perturbation strength $k$ gets larger, the invariant objects become more
irregular, thus smaller continuation steps and
larger grid size are chosen adaptatively according to
Algorithm~\ref{algorithm_continuation}.

As one can observe, while the invariant circle is fluctuating as the
perturbation increases, the isochrons becomes more irregular, and the
minimum angle between the invariant circle and the isochrons is getting closer
to 1. We postpone the detailed discussions regarding this ``bundle
collapsing'' scenario to Section~\ref{sec_breakdown}.

\subsection{Continuation w.r.t. $\eta$}
\label{subsec_continuation_eta}
\begin{figure}[!htb]
  \centering
  \includegraphics[width=0.9\textwidth]{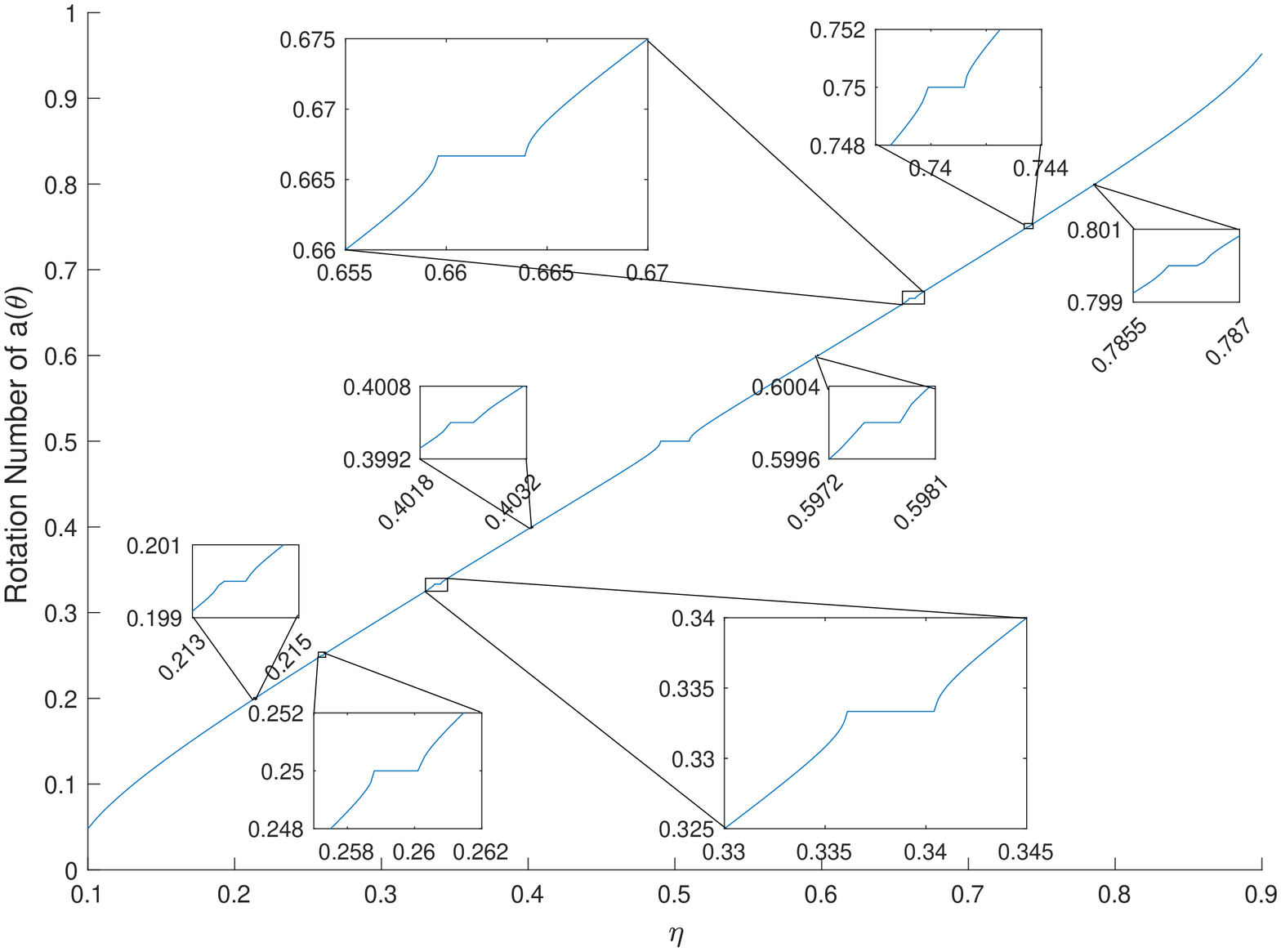}
  \caption{Rotation number of $a(\theta)$ w.r.t. $\eta$}
  \label{rot_num_eta}
  \footnotesize
  \emph{In this example, $k = 5$ and $\gamma = 0.1$.}
\end{figure}

An interesting aspect one can observe by  fixing $k$ and $\gamma$ while
changing $\eta$ is the change of 
rotation number of the internal dynamics $a_{\eta, \gamma, k}(\theta)$.
Figure~\ref{rot_num_eta} presents such $\tau_a$, the rotation number of
$a(\theta)$ as
a function of $\eta$.

The familar Devil's staircase in Figure~\ref{rot_num_eta} is from the classic
results in rotation numbers \cite{K95}, where $\tau_a(\eta)$ is a monotone
function (since our map $f_{\eta, \gamma, k}$ forms a family of
orientation-preserving homeomorphisms), and the stairs (where $\tau_a(\eta)$
have constant values) are corresponding to the rational rotation numbers, in
which case $a(\theta)$ has periodic orbits.

\subsubsection{Computation of the Rotation Number}
The most naive way of computing the rotation number is by the definition:
\begin{equation}
  \tau_a = \lim_{M \rightarrow \infty} \frac{a^{\circ M}(x) - x}{M} \text{ mod } 1
\end{equation}
where $x \in \mathbb{T}$. By computing such rotation number with the help of
the Birkhoff average, one can show that such sequence admits
$\mathcal{O}(\frac{1}{M})$ convergence rate \cite{KU78}.

In this paper, we follow the approach in \cite{DSSY16} and \cite{DSSY18}, which
use the weighted Birkhoff average instead of the regular Birkhoff sum. 

It is shown that, when the inner dynamics is analytically conjugated to a Diophantine rotation, 
for any positive integer $m$, there exists $C_m> 0$ such that 
\begin{equation*}
  \left| \frac{1}{A_M} \sum_{n = 0}^{M - 1} w\left(\frac{n}{M}\right) (a^{\circ (n + 1)}  - a^{\circ n}) \text{ mod 1 }  - \tau_a\right| \leq C_mM^{-m},
\end{equation*}
where $w(t)$ is the exponential weighting function
\begin{equation*}
  w(t) = \begin{cases} \exp\left(\frac{1}{(t(t - 1))^p}\right) & t \in (0, 1), \\ 0 & \text{ elsewhere.}
  \end{cases},
\end{equation*}
$A_N = \sum_{n = 0}^{M - 1} w(\frac{n}{M})$ and our choice of $p$ is 2.

On the other hand, when the system is phase-locked, the algorithm has a much slower convergence. 
This can be taken to our advantage. If the results of two truncated sums are widely different,
we can conclude that the system has a rational rotation number. 

\subsubsection{Distinguishing Irrational Rotations from Rational Rotations}

By classical results \cite{K95}, there are three types of orbits for the
internal dynamics with rational rotation number $\tau_a$:
\begin{enumerate}
\item Periodic orbits with the same period and same order as the rotation
  $\theta + \tau_a$;
\item Homoclinic orbit;
\item Heteroclinic orbit.
\end{enumerate}

\begin{figure}[!htb]
  \centering
  \subfloat[Internal Dynamics]{
    \includegraphics[width=0.45\textwidth]{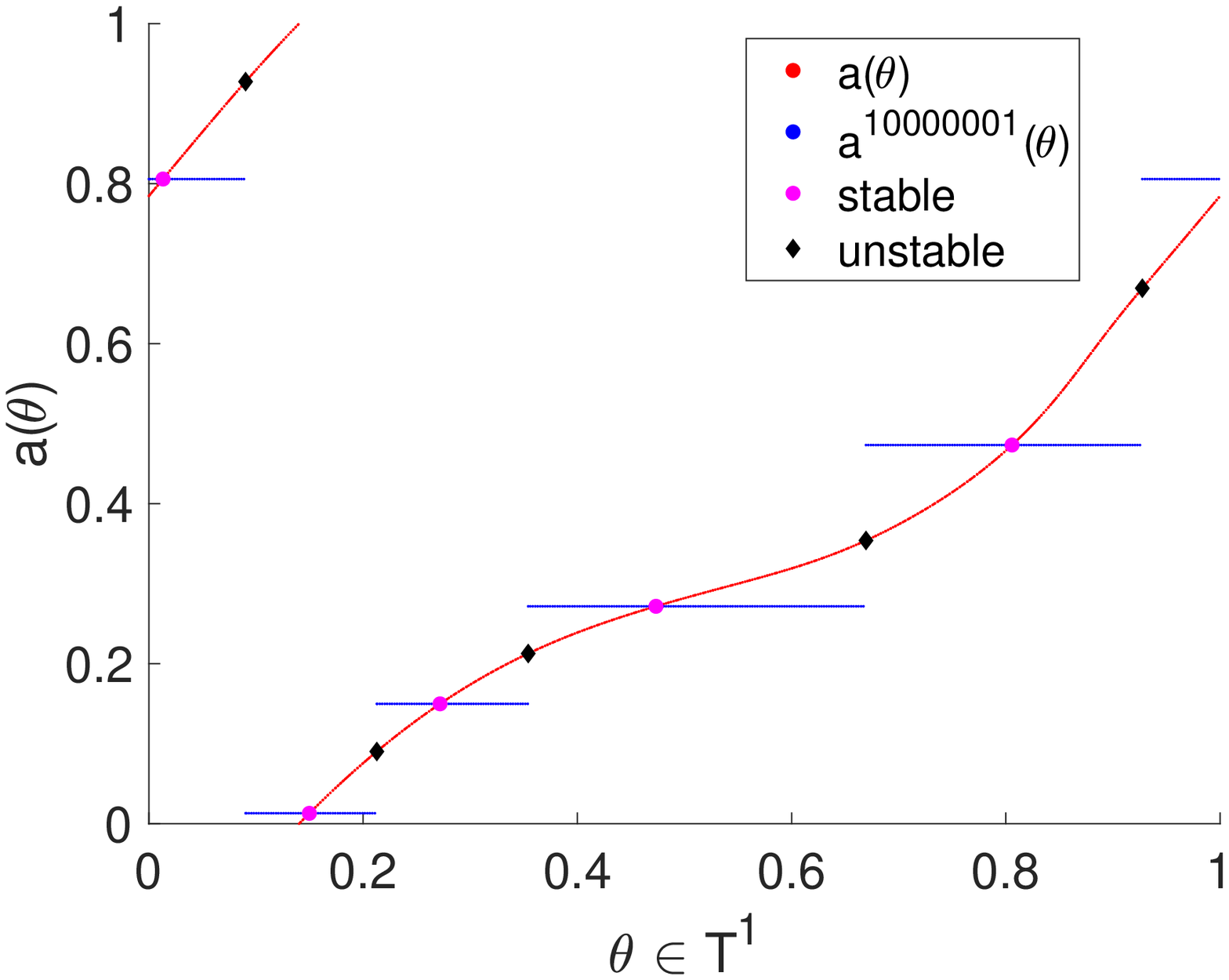}
    \label{fig:f1}
  }
  \hfill
  \subfloat[Invariant Circle and Isochrons]{
    \includegraphics[width=0.45\textwidth]{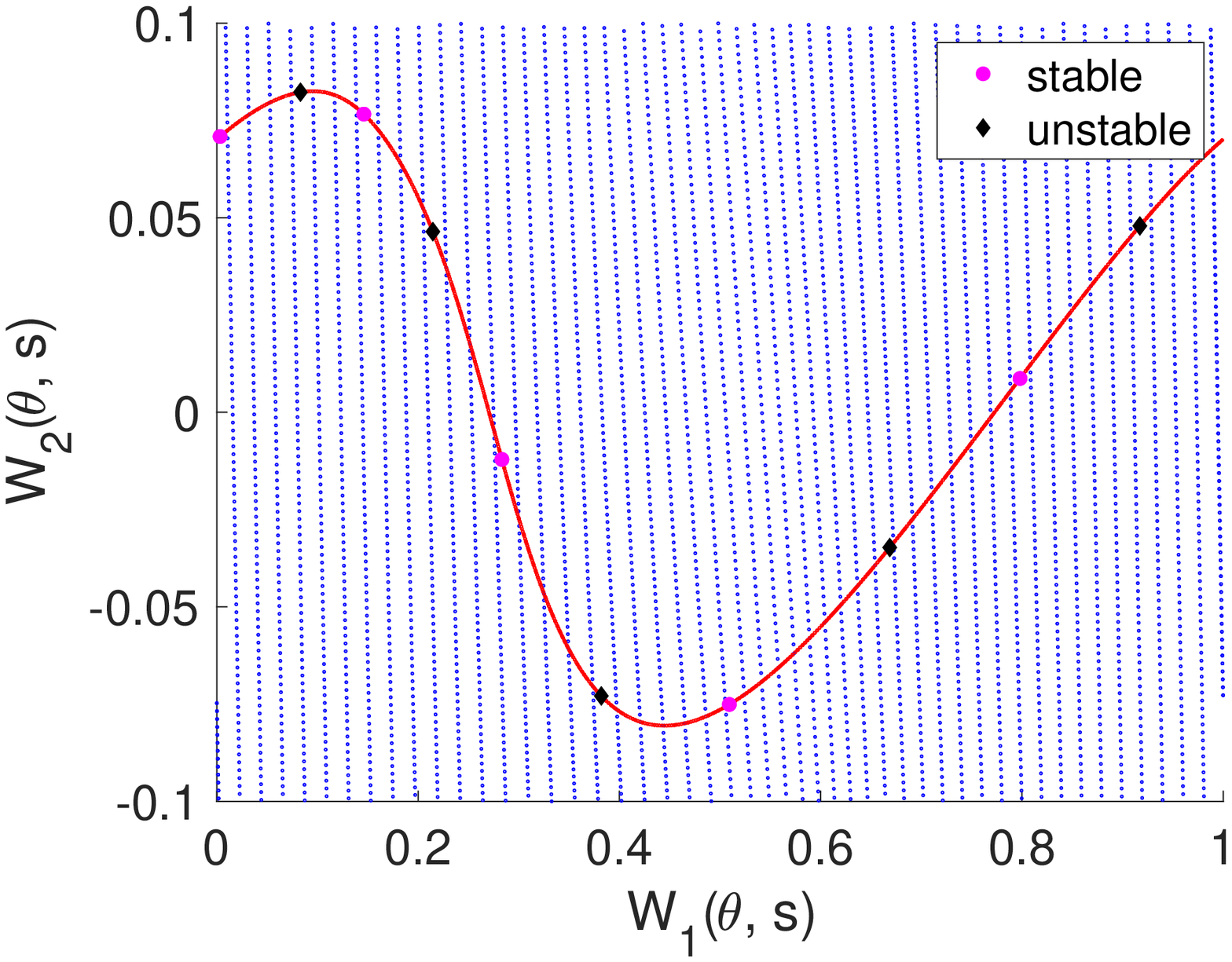}
    \label{fig:f2}
  }
  \caption{Rational Rotation}
  \label{rational_rot_num}
\end{figure}

Thus,  by applying the internal dynamics $a(\theta)$ $M$ times on $\mathbb{T}$
(where $M$ is picked to be a large enough number in case the period has large denominator),
$a^{\circ M}(\theta)$ is either a piecewise constant function (as shown in 
Figure~\ref{rational_rot_num}) or it will be a rational rotation (the latter case is
very unstable and we do not encounter it.) 

On the other hand, if $\tau_a$ is irrational, by the Poincar\'e Classification
theorem, the rotation $R_{\tau_a}$ is at least 
a topological factor of $a^{\circ M}(\theta)$, thus we expect $a^{\circ
  M}(\theta)$ to be at least continuous.

\begin{figure}[!htb]
  \centering
  \subfloat[Internal Dynamics]{
    \includegraphics[width=0.45\textwidth]{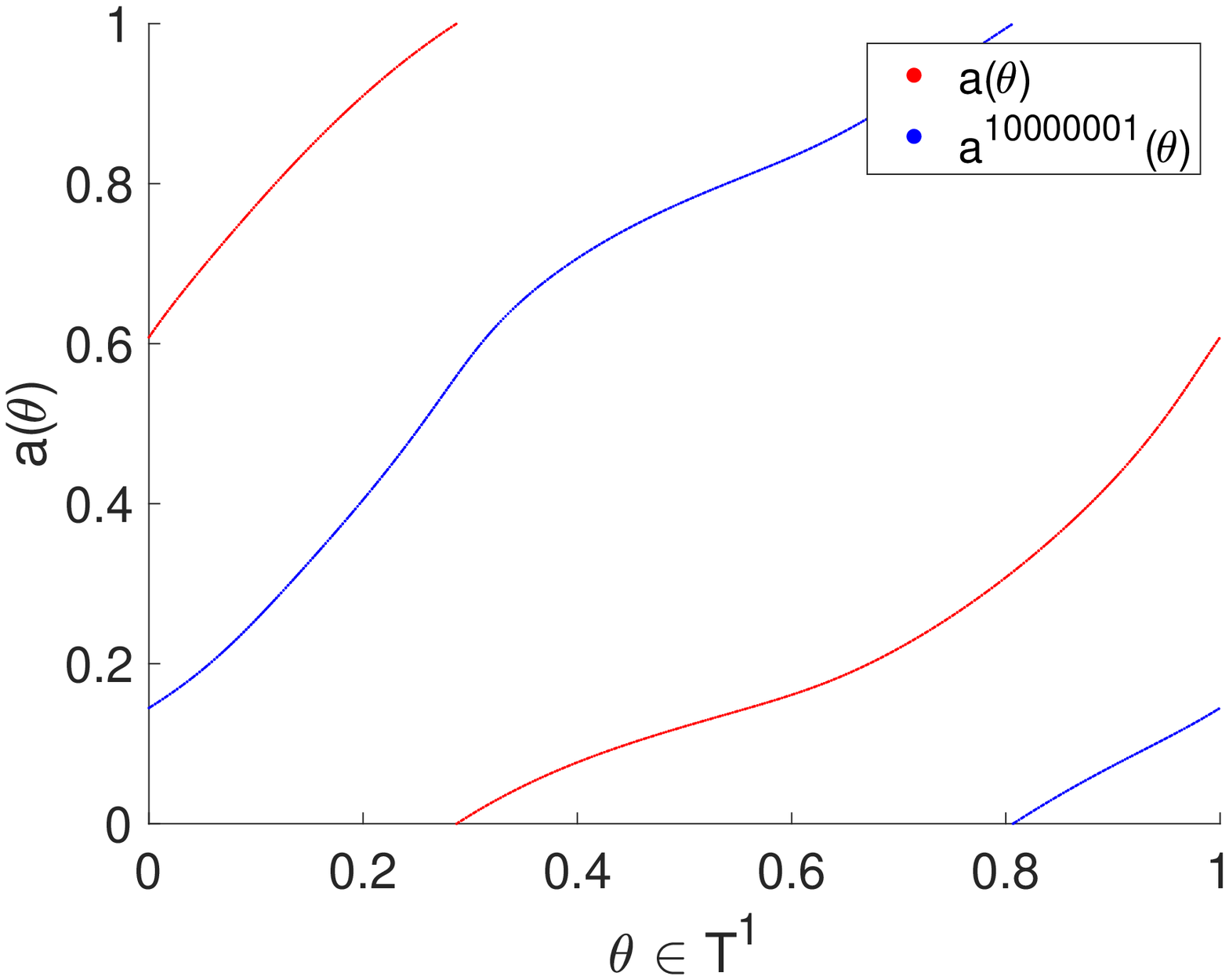}
    \label{fig:f1}
  }
  \hfill
  \subfloat[Invariant Circle and Isochrons]{
    \includegraphics[width=0.45\textwidth]{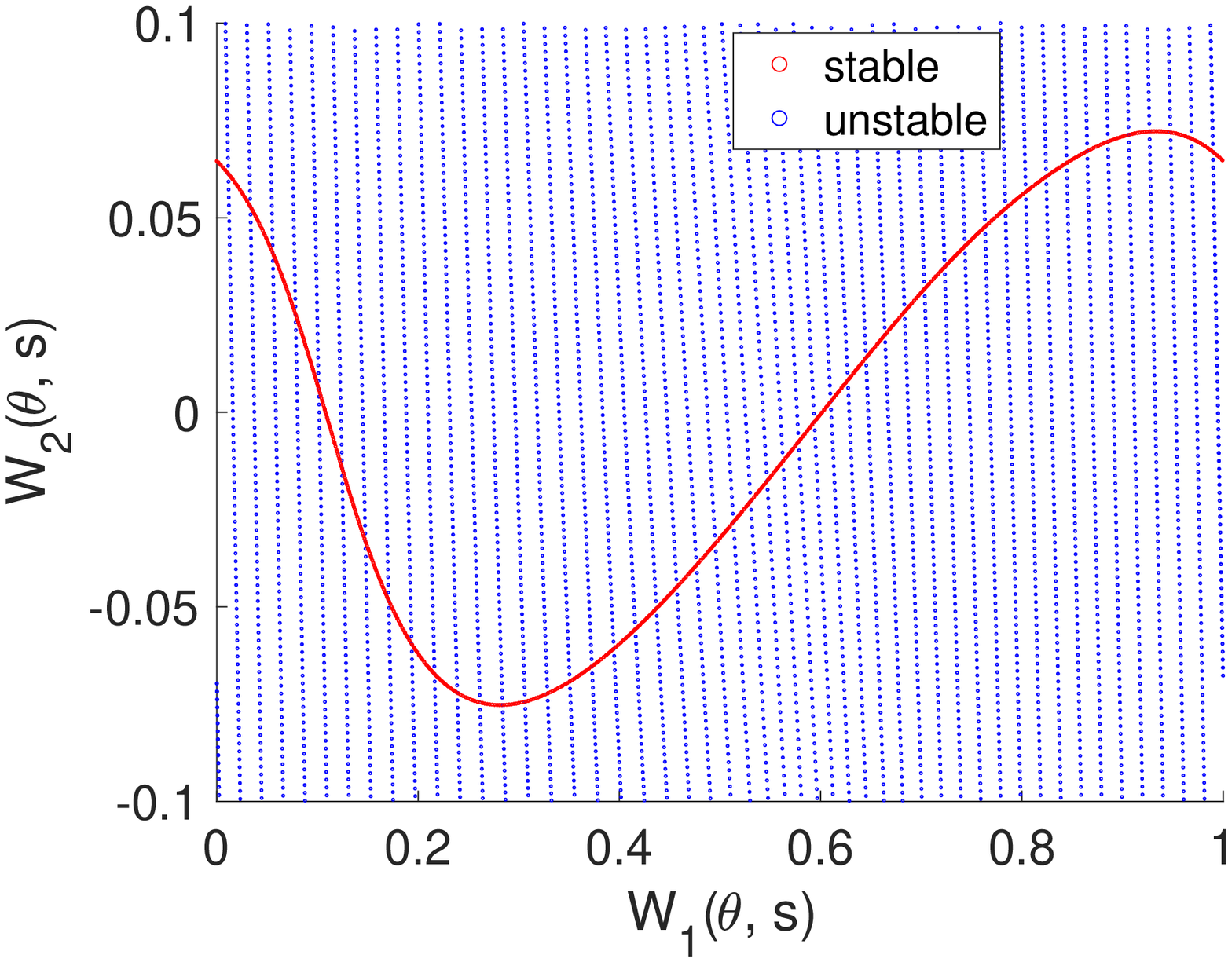}
    \label{fig:f2}
  }
  \caption{Irrational Rotation}
  \label{irrational_rot_num}
\end{figure}

\subsection{Continuation with Prescribed Rotation Number}
\label{subsec_continuation_k_prev}
Many applications (quasiperiodic attractors, oscillators require the
rotation frequency to be a fixed Diophantine number, for example, the golden
mean.

In such cases, as stated in Remark~\ref{remark_constant_lambda}, $a(\theta) =
\theta + \omega$,
$\lambda(\theta) = \lambda$, and the phase-locked phenomenon does not appear.

\begin{figure}[!htb]
  \centering
  \includegraphics[width=0.45\textwidth]{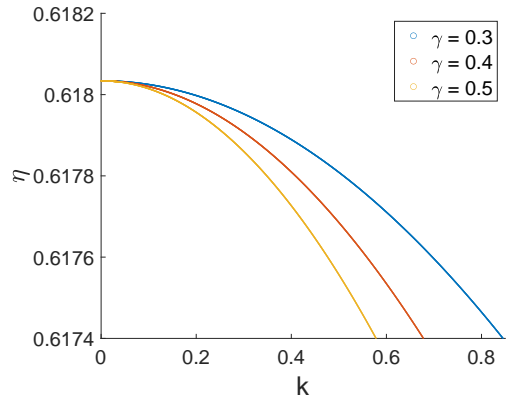}
  \caption{$\eta$ vs $k$ when the rotation number is prescribed by the golden
    mean, with fixed $\gamma$}
  \label{eta_prev_rot}
\end{figure}

In order to cope with the fixed rotation number, we can vary $\eta$ using any
root-finding algorithms (for example, the Brent zero-finding method). The plot for $\eta$ as
$k$ varies is in Figure~\ref{eta_prev_rot}.

In principle, with a prescribed rotation number, one can follow
Algorithm~\ref{algorithm_continuation}, increase the
range for $k$ in Figure~\ref{eta_prev_rot} and make explorations when $k$ is
close to the breakdown.

\def\torus{\mathbb{T}}
\def\real{\mathbb{R}}
\def\C{\mathcal{C}}

\section{Some Explorations on the Breakdown of Circles} \label{sec_breakdown}

Given an algorithm, it is quite natural to study its limits of
applicability. Of course,
the mathematical phenomenon of interest may itself breakdown.
The two breakdowns (of the algorithm and of the real object)
 can, of course, happen at different places.

Furthermore, as we will see, there may be different mathematical 
definitions of breakdown. A natural, widely-studied definition of 
mathematical breakdown is when the invariant set stops from being 
a $C^1$ manifold \cite{Mane78}, but it is also natural 
to consider the case when the circles stop being continuous curves.
The remarkable paper 
\cite{CapinskiK}  uses topological methods to construct 
analogs of the attractors which are not continuous curves, 
but indecomposable continua, so it would also be natural 
to study the breakdown when the  attractor in the system ceases to be 
an indecomposable continua.  Going even further in the perturbations, 
\cite{YoungW1,YoungW2, Levi, JarnikK69} present statistical or topological
descriptions of attractors that are not even continua.

We think it is fair to say that any of these phenomena are not completely
understood (\cite{Palis08} called them \emph{the dark realm}). 
Nevertheless, there are several \emph{``scenarios''} -- we borrow 
the name from \cite{EckmannR85} in a related, but different context -- 
which, with extra hypotheses, are understood (either mathematically or 
numerically).  We summarize them in 
Section~\ref{scenarios}.  Later, in Section~\ref{subsec_breakdown_numerical} 
we will report our numerical  explorations of a route -- chosen more or less at random -- 
which does not seem to fit the previously studied scenarios. We hope
that this will stimulate further research.

\subsection{Several Scenarios for the Breakdown of Invariant Circles} \label{scenarios}

\subsubsection{General Breakdown of Normal Hyperbolicity} 
The theory of normally hyperbolic manifolds \cite{Fenichel71, HirschPS77,Pesin04} 
shows that a $C^1$ manifold persists if 
\begin{itemize}
\item the stable  and the tangent  bundles have different contraction rates;
\item the angle between the stable  and tangent  bundles is  bounded from blow.
\end{itemize}
To measure how well each of these two conditions are
satisfied we can associate a number that quantifies them
(for example the ratio of the rates and the minumum angle between the bundles). 
The size of the
perturbations that are allowed
depends on the sizes of these numbers

The remarkable paper \cite{Mane78} can serve to organize
the questions in the area.  More precisely, 
\cite{Mane78}  showed that if the hyperbolicity conditions 
fail, then a ``generic'' perturbation will not have a $C^1$ invariant manifold nearby 
which is conjugate to the original one. 

Note that \cite{Mane78} leaves several open possibilities
that, indeed happen. 
\begin{itemize} 
\item Persistence of an object which is not a $C^1$ manifold, but 
is still rather regular (e.g a $C^1$ manifold with a few singularities, a $C^{0.99999}$ manifold ). 
\item 
If one restricts the perturbations to be very smooth, it is 
possible to use other arguments, such as KAM theory, to obtain persistence under  in finite   codimension perturbations (i.e. adjusting a finite number of
parameters).  The KAM argument does not need that the invariant manifold
is uniformly hyperbolic. 
\item 
Note also that the loss 
of the hyperbolicity can happen when the stable and tangent bundles become very 
close (even if the contraction rates in the stable direction are separated from 
rates in the tangent direction; we will see numerical observations  of this phenomenon). 

As we will see in Section~\ref{bundlecollape}, in the case while there are
uniform rates of contraction but the stable bundles are close to the tangent, the
breakdown of hyperbolicity can be studied numerically, but have not yet been
understood mathematically.
\end{itemize} 

\subsubsection{Loss of Hyperbolicity because of the Vanishing Rates of Contraction} \label{vanishrate}

In this case, the loss of hyperbolicity happens because the rates of contraction/expansion
in the normal directions get close to $1$, while the bundles remain well 
separated. This scenario has been studied in  \cite{ChencinerI79a, ChencinerI79b,Sell79}. 
In \cite{Los88}, one can find an analog of period-doubling. 

In the above papers -- under higher regularity assumptions -- 
one basically reduces (up to small error) the system to a product
one of whose factors experiences a bifurcation.  Then, one shows
that the errors do not change the conclusion.

\subsubsection{Loss of Hyperbolicity by Vanishing of the Angle: Bundle-collapse
  Scenario} \label{bundlecollape}
One can apply the KAM theory  \cite{CallejaCL13} to continue  tori with 
a Diophantine rotation. This happens in a codimension 1 family in 
the map \eqref{DST} (adjusting the parameter $\eta$). 

Because the map on the circle conjugates to a rotation, the tangential expansion rate is
$1$. Because \eqref{DST} has constant determinant $\lambda$, the normal 
contraction rate has to be $\lambda$.  Hence, these KAM tori are also 
Normally Hyperbolic Invariant Manifolds and their normal exponent is bounded away from $1$
in all their range of existence (we also remark that 
\cite{fenichel73} shows that the circles have to be very smooth -- 
hence no Denjoy counterexamples and with the global conjugacy results of
Herman,  it is shown in \cite{CallejaCL13} that the 
circles  remain analytic and analytically conjugate to a rotation up to the breakdown. It is 
conjectured numerically and predicted by renormalization group  that at the breakdown, 
there is a finitely differentiable circle.)

It follows from the above 
rigorous arguments that the only way that the invariant circles in \eqref{DST} with Diophantine 
rotation number can cease to exist is 
when the hyperbolicity losses due to the angle between 
the stable direction and the tangent to the manifold goes to zero. 

This phenomenon was discovered and studied empirically in \cite{CallejaF12}. 
They not only found numerical evidence that bundle collapse happens at a precise value of 
the parameters.  They also found evidence of universal scaling relations in 
the size of norms of the conjugacy, the angles of the bundles, and the change in 
rates of contractions. There was a limiting regularity of the torus which 
seemed to be universal. 

Similar phenomena happened in other contexts in \cite{HaroL07, FiguerasH13,FiguerasH15}. 
In \cite{HaroL07}, one can find an example where bundle collapse is rigorously shown 
to happen.  The details of the phenomenon and the scaling relations found 
remain a challenge for rigorous mathematics.  Important progress is 
in  \cite{BjerklovS08, Ohlson17}.

\subsubsection{Breakdown of Phase-locked invariant Circles}
\label{sec:phaselocked}

We also refer to \cite[Section 3.1]{BroerST98} for 
a very detailed numerical study of some examples and very illuminating
general discussions.  The paper \cite{BroerST98} presents arguments about 
why the present scenario may be relevant for higher-dimensional phenomena. 

Consider a family of analytic maps of the plane. 
We assume that for all the values of the parameter, 
there is a hyperbolic periodic orbit $p_h$  and an attractive periodic orbit $p_a$ of
the same period $p_h$. Denote 
the basin of attraction of $p_a$ by domain $U$. 

Furthermore, we
assume that $W^u$, the unstable manifold of $p_h$, intersects 
$U$. 
Note that the above situation is stable under small $C^1$ perturbations. 

In such a case, because of the invariance of the unstable manifold of 
$p_h$, we have that $W^u$  together with $p_a$ 
form a circle $\C$ which is normally contractive. Notice that $\C$ is
analytic away from $p_a$. 

We denote the eigenvalues at $p_a$ as $0 <  \lambda_1, \lambda_2 < 1$. 

We recall that if $\lambda_2 < \lambda_1$, there are classical results 
showing that there is an analytic   \emph{strong stable manifold}  tangent to 
the eigenspace of $\lambda_2$. As a matter of fact, in 
\cite{Llave03}, it is shown that the standard method applies when 
$\lambda_2^2 < \lambda_1$. This improved result allows 
that $\lambda_2 > \lambda_1$ (this will prove useful 
when we consider dependence on parameters).  
  Related results when  the eigenvalues are non-resonant 
appear in \cite{CFdlL03a,Llave97}.

We denote $r^* = \log(\lambda_1)/\log(\lambda_2)$. 

We will assume that the $W^u$  intersects transversally the 
stable manifold associated with $\lambda_2$. The above assumption is very robust.

Then, we will now show 
that  the union of $W^u$ and $p_a$ is a circle which is 
$C^{r^* - \delta}$ for every $\delta> 0$ and which is not $C^{r^* +\eta}$ 
for any $\eta > 0$. \footnote{
We can assume, without loss of generality, that the 
manifold associated to $\lambda_2$ is the $x$ axis, and, by  Sternberg 
linearization of contractions, the map is just linear. 
If  $W^s$ is expressed as the graph of a function $\phi$ on an interval
$x \in [\lambda_1 x_0, x_0]$, we see that the full 
manifold will be the graph of $\sum_j \lambda_1^j \phi(\lambda_2^{-j} )$. 
Note that each terms in the summation is defined in  a non-overlapping interval. 
It is easy to see that the derivatives of  
$\lambda_1^j \phi(\lambda_2^{-j} )$ for order more than $r^*$ grow unbounded.
Those derivatives of order less than $r^*$ 
converge exponentially fast to zero. 
}

It is not difficult to construct families of maps $f_\mu$ 
with the two fixed points as above, for which $W^u$ intersecting 
transversally the manifold it associates to and the exponents of
$p_{a,\mu}$ satisfies
\[
\begin{split} 
& \lambda_2 (\mu)  > \lambda_1(\mu) > \lambda_2(\mu)^2 \text{ for } \mu \in [0, \mu_0), \\ 
& \lambda_2 (\mu)  <  \lambda_1(\mu) \text{ for } \mu \in ( \mu_0, 1),  \\ 
& \lim_{\mu \to 1} \lambda_1 (\mu) = 1.
\end{split} 
\] 
Note that such assumptions are also persistent when we change the family slightly.
To fix ideas, we can think of $\lambda_2(\mu)$ independent of $\mu$  and $\lambda_1(\mu)$
increasing from $\lambda_2^2 $ to $1$. 

In the interval $\mu \in [0,\mu_0)$, the manifold attached to $\lambda_2$
is the weak stable manifold, and the results of \cite{Llave03} allow us 
to construct it.  In this interval, $r^*(\mu) > 1$ and the circle $\C$ is 
a normally hyperbolic manifold. 

At $\mu = \mu_0$, we are at the critical value of \cite{Mane78} 
and indeed, there is no persistence as $C^1$ manifold. 

For $\mu \in (\mu_0, 1]$, we have that $r^*(\mu) < 1$ and that the circle
is indeed not a $C^1$ manifold, but it still has a H\"older regularity, 
which can be very close to $1$ near the critical value of loss of hyperbolicity. 

We also have that $\lim r^*(\mu) = 0 $ so that the H\"older regularity
is getting closer to zero.

If the family continues for $\mu > 1 $ and $p_a$ experiences a saddle-node bifurcation
at $\mu = 1$,  
one can see that the limiting invariant object is not a continuous circle 
since the  oscillations accumulate. On the other hand, it is a continuum.

\subsection{Some Numerical Explorations} \label{subsec_breakdown_numerical}

In the previous sections, we have presented some different scenarios
leading to breakdown. Note that both of them assume 
that the rotation number of the invariant circle is maintained constant 
till the breakdown (for the phase-locked maps, this can happen on open sets 
in families, but for Diophantine irrational rotation numbers, it is a codimension $1$ phenomenon. 

In this section, we study numerically  the breakdown of
a family inside \eqref{DST}, chosen arbitrarily. In this family, the 
rotation number changes and goes from rational to irrational.
So, the scenarios Sections~\ref{bundlecollape} and \ref{sec:phaselocked} 
will alternate and the phenomena observed will be a combination of 
the two scenarios.  

This section reports our findings when taking the numerical 
algorithms  to the limit.  We also raise the need for a more detailed mathematical 
theory. 

For invariant attractors, the Lyapunov multiplier is always smaller than 1, thus
in the dissipative standard map \eqref{DST},
we are expecting either the ``bundle merging'' scenario or the ``rate meeting''
scenario (when $\lambda^s = \mu^s$ as in Section~\ref{sec:phaselocked}).

\begin{figure}[!htb]
  \subfloat[Angles between the tangent and normal bundles for Figure~\ref{plot_all}]{
    \includegraphics[width=0.43\textwidth]{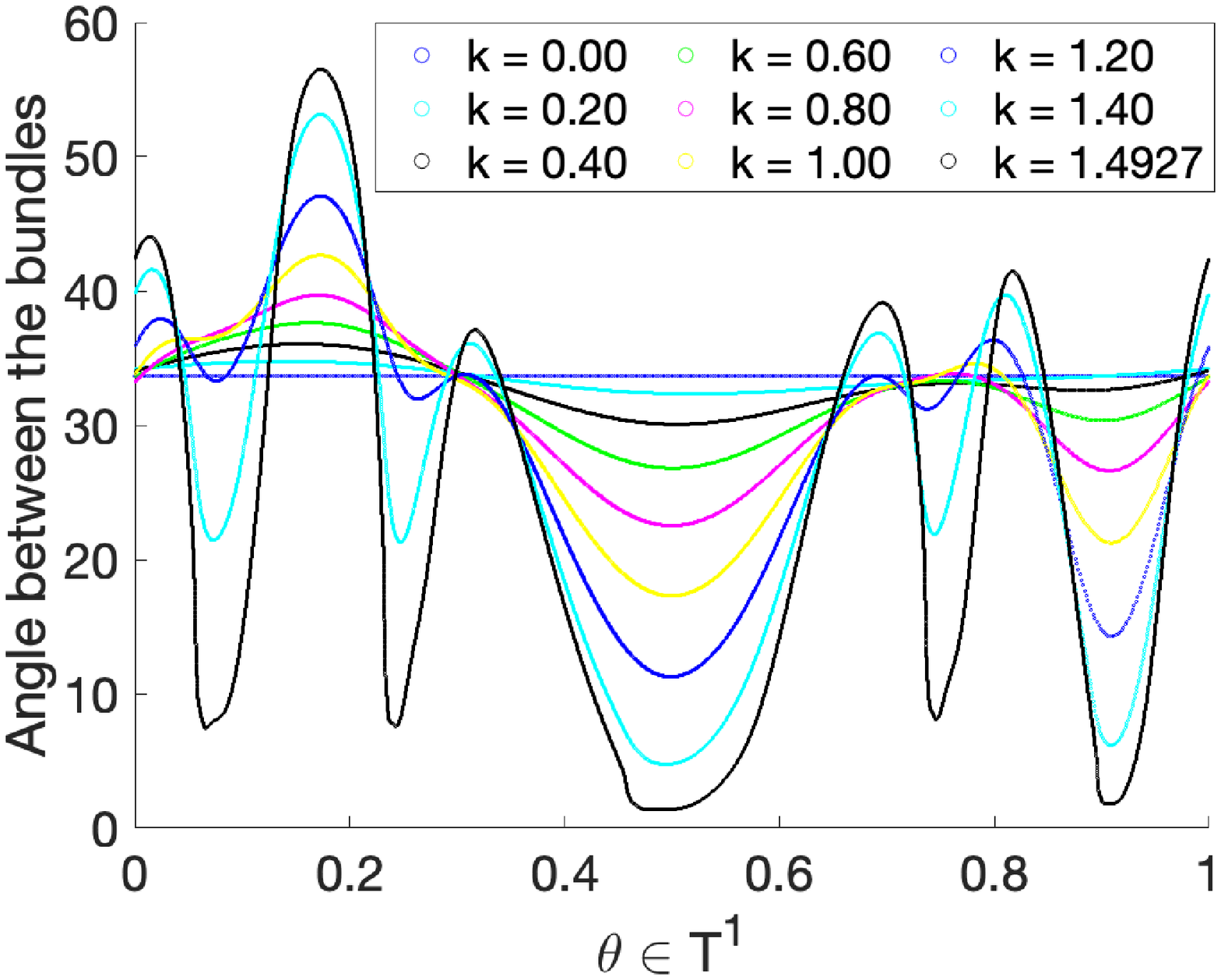}
    \label{fig:all_angle}
  }
  \subfloat[The minimum angle between the tangent and normal bundles]{
    \includegraphics[width=0.43\textwidth]{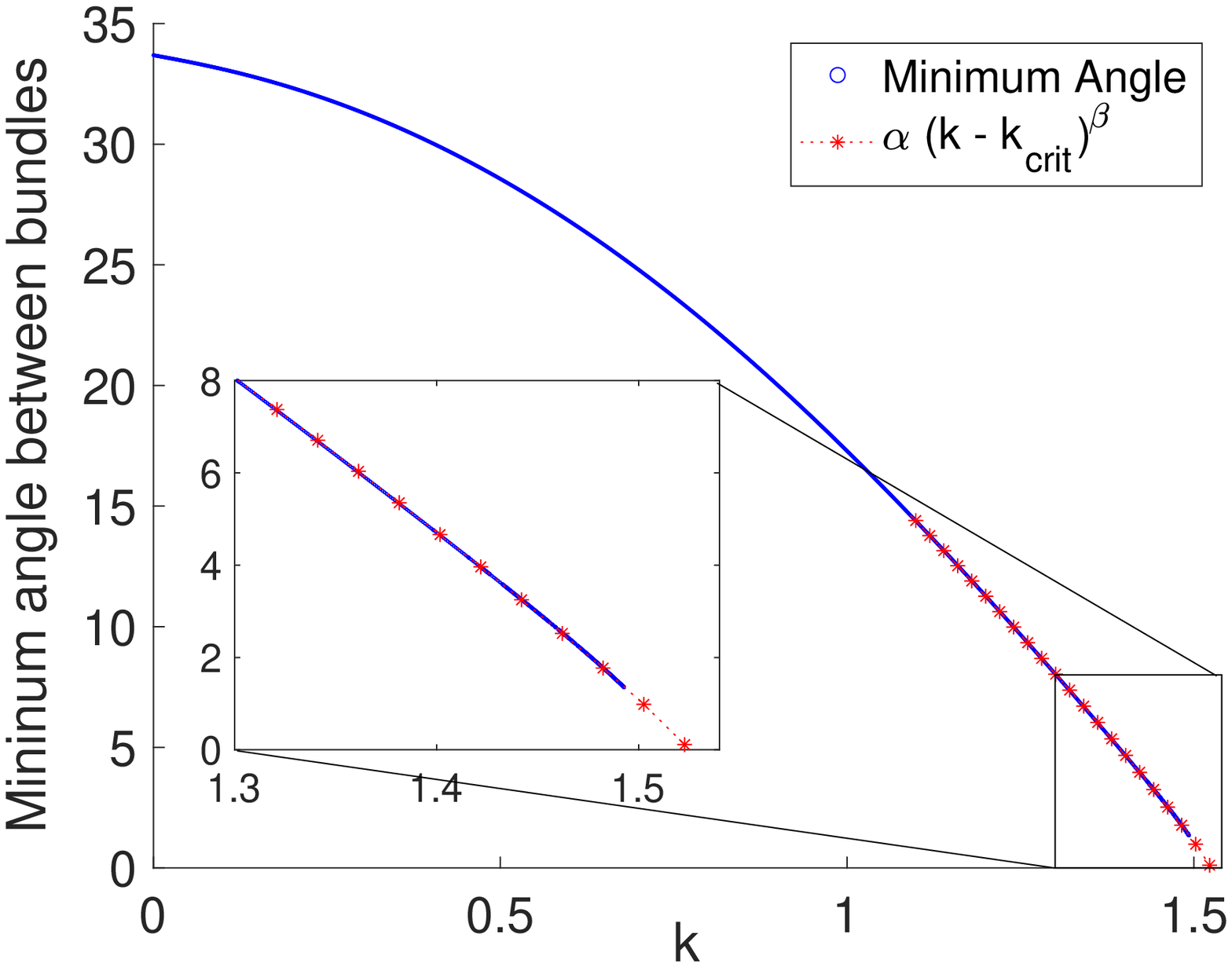}
    \label{fig:minAngleVal}
  }\\
  \subfloat[The position 
  such that the minimum angle in Figure~\ref{fig:minAngleVal} is achieved]{
    \includegraphics[width=0.43\textwidth]{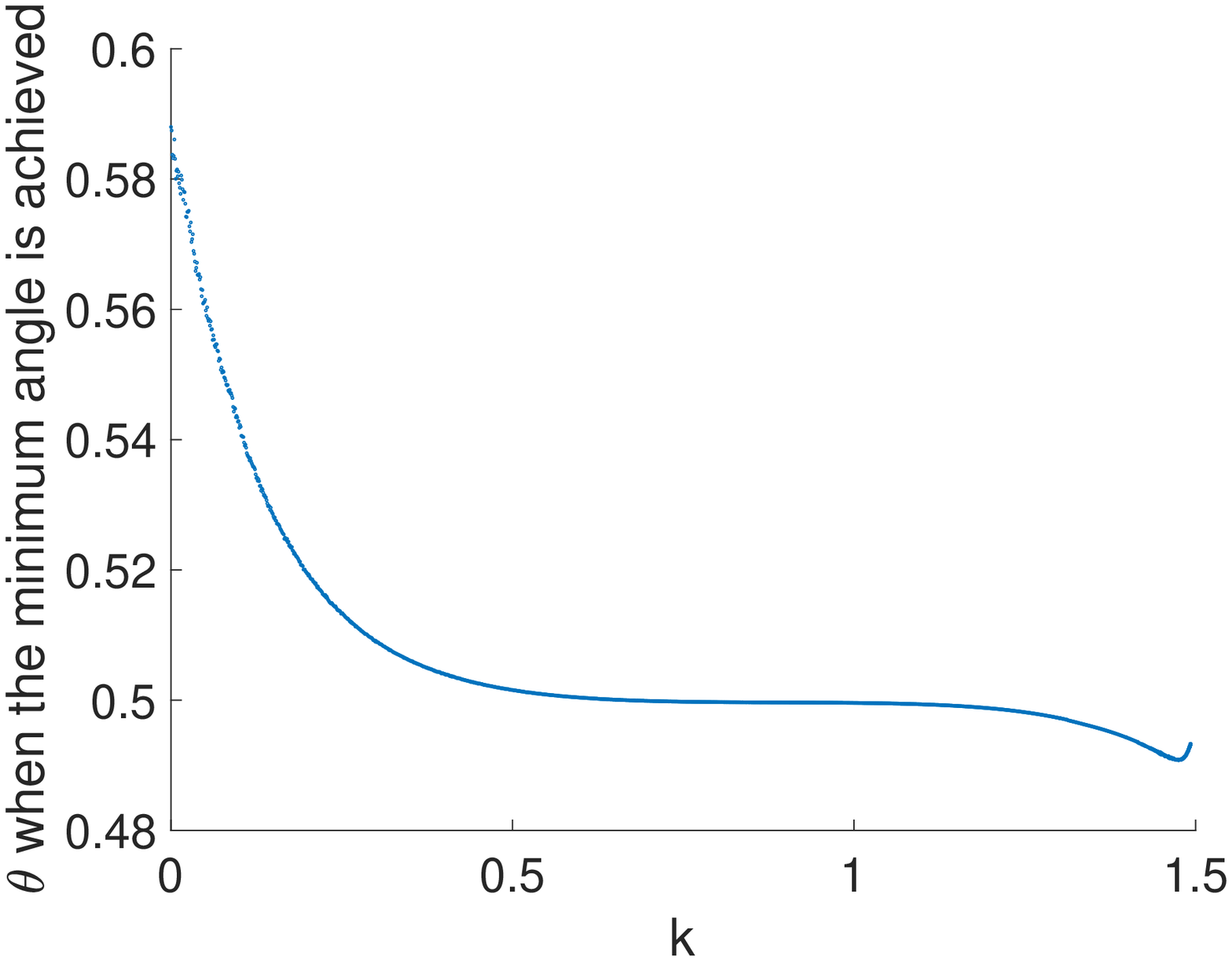}
    \label{fig:minAngleTheta}
  }
  \subfloat[Rotation number w.r.t. the perturbation parameter $k$]{
    \includegraphics[width=0.43\textwidth]{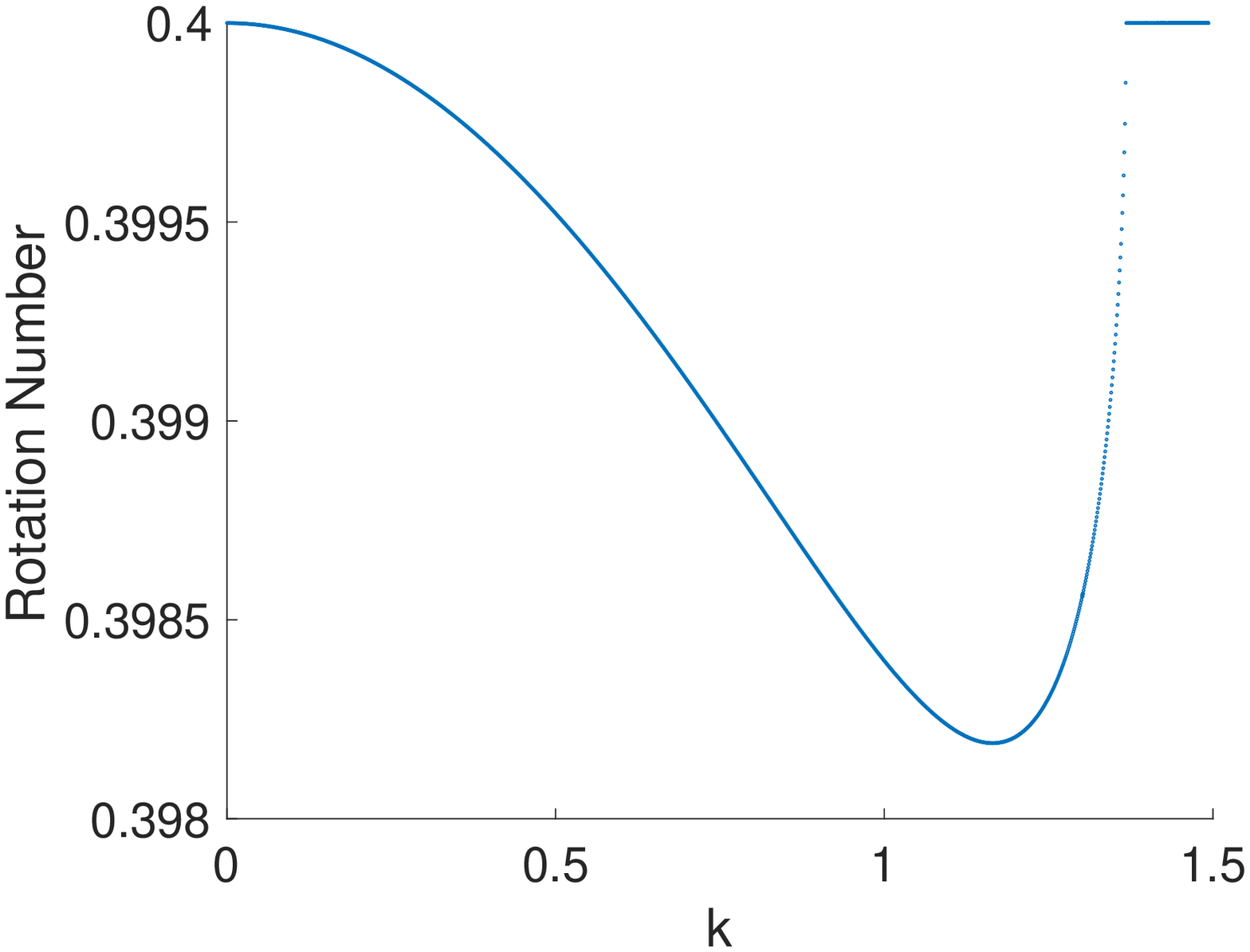}
    \label{fig:rot_wrt_k}
  }\\

  \subfloat[$r^*$ and eigenvalues of $f^{\circ 5}$ for the periodic orbits when $k > 1.3712$.]
  {
    \includegraphics[width=0.43\textwidth]{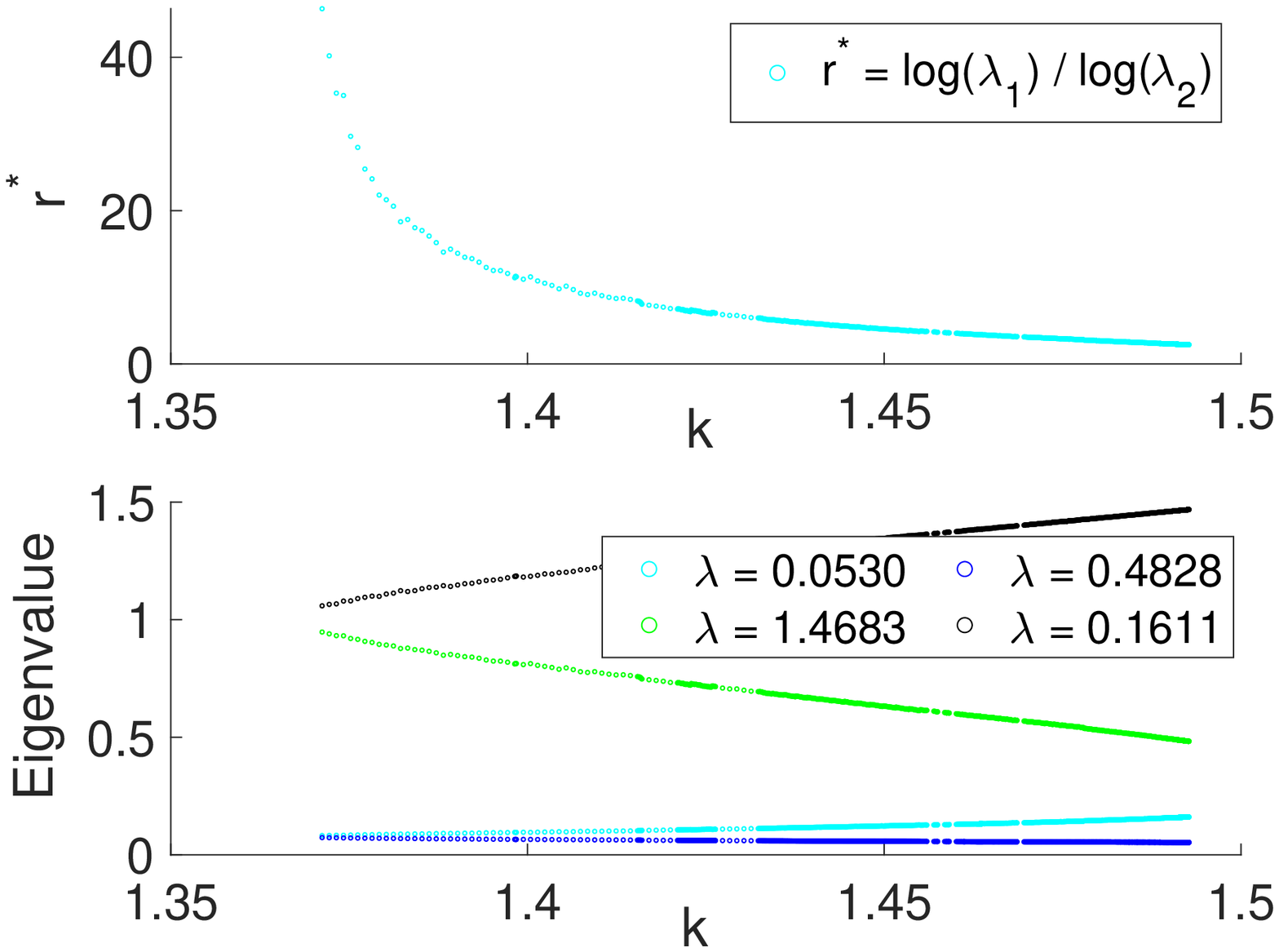}
    \label{fig:eigenPlot}
  }
  \subfloat[Stability of the attractive and hyperbolic periodic orbit when $k = 1.49271$]{
    \includegraphics[width=0.43\textwidth]{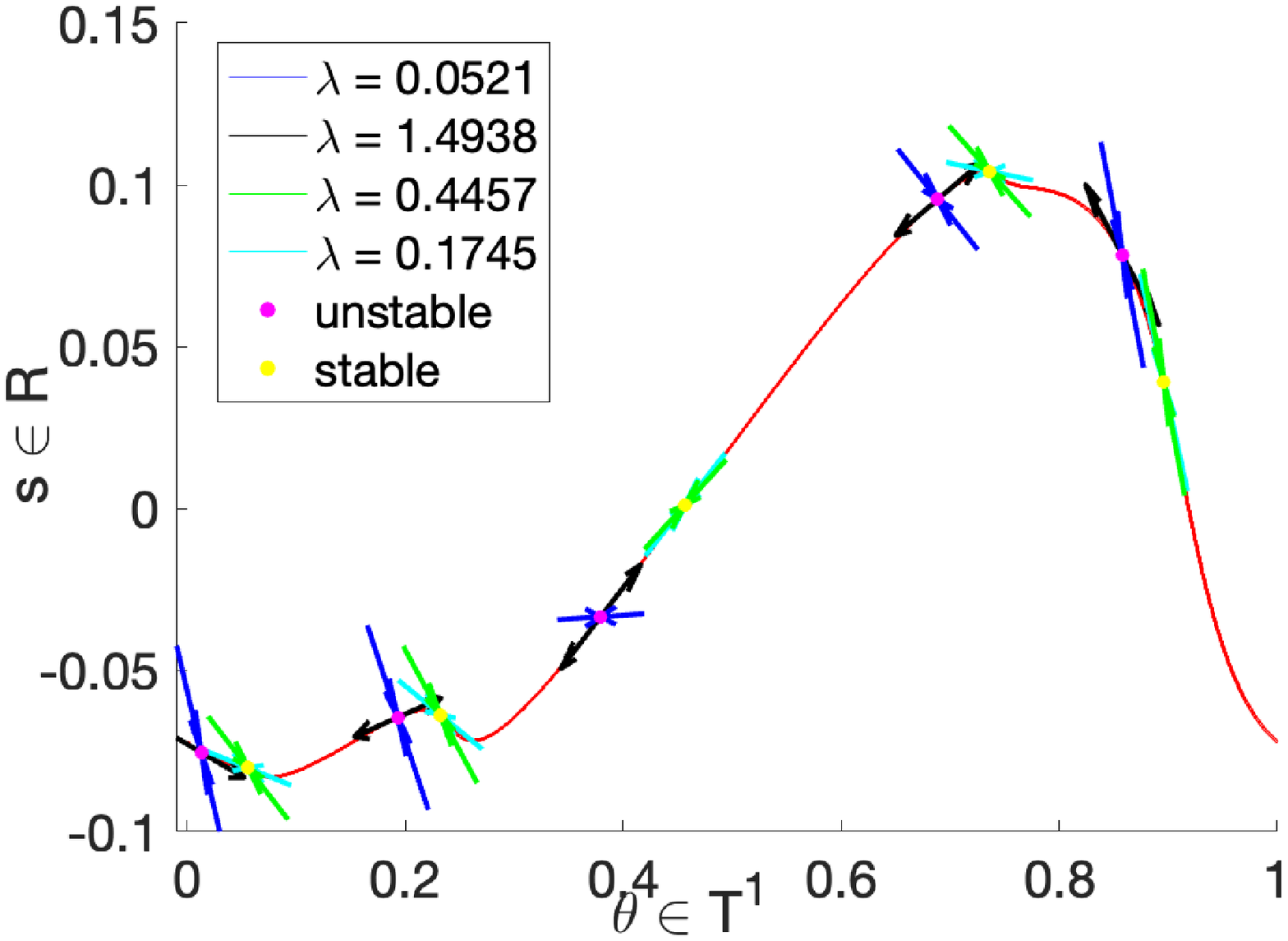}
    \label{fig:eigenPlotEnd}
  }\\
  \caption{``Bundle merging'' Explorations}
  \label{fig:breakdown}
\end{figure}

Following a random continuation path as computed in
Section~\ref{subsec_continuation_k}, where the rotation number is not
prescribed, we
have Figure~\ref{fig:breakdown}.

To quantify the separation between bundles, we measure the minimum angle between
the invariant circle and the corresponding stable manifolds.
Figure~\ref{fig:all_angle} indicates the angles for the $k$ values as in
Figure~\ref{plot_all}. Figure~\ref{fig:minAngleVal} plots the minimum angle with
respect to $k$, and Figure \ref{fig:minAngleTheta} plots the corresponding
$\theta \in \mathbb{T}$ such that the minimum angle is reached.

As indicated in Figure~\ref{fig:minAngleVal}, the slope of the blue curve
(minimum angle) is approaching negative infinity (it is not obvious in the plot
due to the different scales in the horizontal and vertical axis). As a matter of
fact, the later
portion of the minimum angle curve can be fitted by the following asymptotic
expression using nonlinear regression:
\begin{equation}
  \label{breakdown_eq}
  \alpha (k - k_{crit})^{\beta},
\end{equation}
where $\alpha$ is the scaling constant, $\beta$ is the scaling exponent, and
$k_{crit}$ is the breakdown threshold. In fact, $\alpha = 45.8879$, $\beta =
0.9085$ and $k_{crit} = 1.5247$ gives us an estimation for the critical value
such that the tangent and normal bundle collapse.

\begin{remark}
  The asymptotic expression \eqref{breakdown_eq} appears in \cite{CC10} when it
  comes to estimating the blow-up of the Sobolev norm of the parameterization of
  the invariant circle. Empirically, Figure~\ref{fig:minAngleVal} shows that the
  minimum angle also approaches 0 following such power laws with universal exponents.
\end{remark}

Since the rotation number in this example is not prescribed, it may change along
the continuation. In our example, the change of the rotation number as $k$
increases is presented in Figure~\ref{fig:rot_wrt_k}. 

\begin{remark}
  Note that the
  continuous family of
  orientation-preserving circle  homeomorphism $a_{k}$ is not
  monotone, in the sense that the lift of $a_{k}$ can not be
  compared universally on $\mathbb{T}$ \cite{K95}, thus the rotation number in
  Figure~\ref{fig:rot_wrt_k} is not expected to be monotone.
\end{remark}

\begin{remark}
   We are still expecting ``stairs'' (intervals for $k$ in which case the
   rotation number is rational). In Figure~\ref{fig:rot_wrt_k}, the only obvious
   one is when the rotation number is $0.4$, while the rest of the
   ``stairs'' (for example, corresponding to rotation number $0.3999$) is too
   small to be spotted.
 \end{remark}

As indicated in Figure~\ref{fig:rot_wrt_k}, the rotation number is rational when $k
> 1.3712$. In this region, as discussed in Section~\ref{sec:phaselocked},
there are two periodic orbit, the attractive one with eigenvalues $\lambda^s,
\mu^s$, and the hyperbolic one with eigenvalues $\lambda^u, \mu^u$. In
Figure~\ref{fig:eigenPlot}, one can see how the eigenvalue changes as $k$
increases, as well as how the regularity $r^*$ (as discussed in Section~\ref{sec:phaselocked})
changes.
The separation between $\lambda^s$ and $\mu^s$ indicates that the contraction
rate in the tangential and normal direction will never meet in this example.
Figure~\ref{fig:eigenPlotEnd} present the eigenvectors on each periodic point.

\begin{remark}
  \cite{HaroL07} has explored the bundle-merging phenomenon for quasi-periodic
  maps.  Here, we have seen that a similar  phenomenon happens for families in which 
the circles in the family alternate between phase-locked and rotational. Some 
quantitative aspects of the phenomenon can be different. 
\end{remark}

\section{Numerical  algorithm in  3-D maps}
\label{sec_example_3d}
In this section, we further generalize our algorithm to three-dimensional maps
which possesses one-dimensional invariant circles.

Following the same idea introduced previously, we first briefly
discuss the derivation of the algorithm, and then implement such algorithm to
the 3-dimensional Fattened Arnold Family (3D-FAF) maps. 

\subsection{Basic Derivation of the Algorithm} \label{3d_derivation}
Given a map $f: \mathbb{T}^1 \times \mathbb{R}^2 \rightarrow \mathbb{T}^1 \times
\mathbb{R}^2$ that induces a one-dimensional invariant circle, the goal in this
subsection is to derive an algorithm computing both the invariant circle and the
corresponding isochrons.

Because of the increase of dimension on the isochrons, the analogue of  invariance equation
\eqref{invariance} is: 
\begin{equation} 
  \label{invariance_3d}
  f \circ W(\theta, s_1, s_2) = W(a(\theta), \lambda_1(\theta)s, \lambda_2(\theta)s),
\end{equation}
where $W: \mathbb{T}^1 \times \mathbb{R}^2 \rightarrow \mathbb{T}^1 \times \mathbb{R}^2$, $a:
\mathbb{T}^1 \rightarrow \mathbb{T}^1$, $\lambda_1: \mathbb{T}^1 \rightarrow
\mathbb{T}^1$ and $\lambda_2: \mathbb{T}^1 \rightarrow \mathbb{T}^1$ are the
unknowns. Again, $a$ is the internal dynamics on the invariant circle, and for
any given $\theta_0 \in \mathbb{T}$, $W(\theta_0, s_1, s_2)$ parameterizes the isochron,
$\lambda_1(\theta_0)s_1, \lambda_2(\theta_0)s_2$ are the dynamics on the isochron along the
eigen-directions.

Following the same procedure as in Section~\ref{sec_algorithm}, we are looking
for $\Delta_a(\theta)$, $\Delta_{\lambda_1}(\theta)$, $\Delta_{\lambda_2}(\theta)$ and
$\Delta_W(\theta, s_1, s_2) = DW(\theta, s_1, s_2) \Gamma(\theta, s_1, 
s_2)$ satisfing the
following three equations:
\begin{equation}
  \label{eq1_3d}
  Da(\theta) \Gamma_1(\theta, s_1, s_2) - \Delta_a(\theta) - \Gamma_1(a(\theta), \lambda_1(\theta)s_1, \lambda_2(\theta)s_2) = M_1(\theta, s_1, s_2),
\end{equation}
\begin{equation}
  \label{eq2_3d}
  \lambda_1(\theta)\Gamma_2(\theta, s_1, s_2) - \Delta_{\lambda_1}(\theta)s_1 - \Gamma_2(a(\theta), \lambda_1(\theta)s_1, \lambda_2(\theta)s_2) = M_2(\theta, s_1, s_2),
\end{equation}
\begin{equation}
  \label{eq3_3d}
  \lambda_2(\theta)\Gamma_3(\theta, s_1, s_2) - \Delta_{\lambda_2}(\theta)s_2 - \Gamma_3(a(\theta), \lambda_1(\theta)s_1, \lambda_2(\theta)s_2) = M_3(\theta, s_1, s_2),
\end{equation}
where
\begin{equation*}
  M(\theta, s_1, s_2) = \widetilde{e}(\theta, s_1, s_2) - \begin{pmatrix} 0 \\ D\lambda_1(\theta)s_1\Gamma_1(\theta, s_1, s_2) \\ D\lambda_2(\theta)s_1\Gamma_1(\theta, s_1, s_2) 
  \end{pmatrix},
\end{equation*}
and
\begin{equation*}
 \widetilde{e}(\theta, s_1, s_2) = - [DW(a(\theta), \lambda_1(\theta)s_1,
\lambda_2(\theta)s_2)]^{-1}e(\theta, s_1, s_2). 
\end{equation*}

By discretizing function $g$ with three variables $(\theta, s_1, s_2)$ as
\begin{equation}
  \label{discretize_3d}
  g(\theta, s_1, s_2) = \sum_{x = 0}^{\infty}\sum_{y = 0}^{\infty}g^{(x, y)}(\theta)s_1^xs_2^y,
\end{equation}
we again can further discritize Equation \eqref{eq1_3d}, \eqref{eq2_3d} and
\eqref{eq3_3d} order by order as follows:
\begin{itemize}
\item Equation \eqref{eq1_3d}:
  \begin{itemize}
  \item Order $(0, 0)$:
    \begin{equation} \label{eq1_00}
      Da(\theta)\Gamma_1^{(0, 0)}(\theta) - \Delta a(\theta) - \Gamma_1^{(0, 0)}(a(\theta)) = M_1^{(0, 0)}(\theta),
    \end{equation}
  \item Order $(x, y)$:
    \begin{equation} \label{eq1_xy}
      \Gamma_1^{(x, y)}(\theta) = \frac{\lambda_1^x(\theta)\lambda_2^y(\theta)}{Da(\theta)} \Gamma_1^{(x, y)}(a(\theta)) + \frac{1}{Da(\theta)}M_1^{(x, y)}(\theta)
    \end{equation}
  \end{itemize}
\item Equation \eqref{eq2_3d}:
  \begin{itemize}
  \item Order $(0, 0)$:
    \begin{equation} \label{eq2_00}
      \lambda_1(\theta)\Gamma_2^{(0, 0)}(\theta) - \Gamma_2^{(0, 0)}(a(\theta)) = M_2^{(0, 0)}(\theta),
    \end{equation}
  \item Order $(1, 0)$:
    \begin{equation} \label{eq2_10}
      \lambda_1(\theta)\Gamma_2^{(1, 0)}(\theta) - \Delta_{\lambda_1}(\theta) - \Gamma_2^{(1, 0)}(a(\theta))\lambda_1(\theta) = M_2^{(1, 0)}(\theta),
    \end{equation}
  \item Order $(x, y)$:
    \begin{equation} \label{eq2_xy}
      \Gamma_2^{(x, y)}(\theta) = \lambda_1^{x - 1}(\theta)\lambda_2^{y}(\theta) \Gamma_2^{(x, y)}(a(\theta)) + \frac{M_2^{(x, y)}(\theta)}{\lambda_1(\theta)},
    \end{equation}
  \end{itemize}
\item Equation \eqref{eq3_3d}:
  \begin{itemize}
  \item Order $(0, 0)$:
    \begin{equation} \label{eq3_00}
      \lambda_2(\theta)\Gamma_3^{(0, 0)}(\theta) - \Gamma_3^{(0, 0)}(a(\theta)) = M_3^{(0, 0)}(\theta),
    \end{equation}
  \item Order $(0, 1)$:
    \begin{equation} \label{eq3_01}
      \lambda_2(\theta)\Gamma_3^{(1, 0)}(\theta) - \Delta_{\lambda_2}(\theta) - \Gamma_3^{(1, 0)}(a(\theta))\lambda_2(\theta) = M_3^{(1, 0)}(\theta),
    \end{equation}
  \item Order $(x, y)$:
    \begin{equation} \label{eq3_xy}
      \Gamma_3^{(x, y)}(\theta) = \lambda_1^x(\theta)\lambda_2^{y - 1}(\theta) \Gamma_3^{(x, y)}(a(\theta)) + \frac{M_3^{(x, y)}(\theta)}{\lambda_2(\theta)}.
    \end{equation}
  \end{itemize}
\end{itemize}

Equations \eqref{eq1_00}, \eqref{eq2_10} and \eqref{eq3_01} are
underdetermined equations that can be solved by letting
\begin{equation*}
  \Gamma_1^{(0, 0)}(\theta) = \Gamma_2^{(1, 0)}(\theta) = \Gamma_3^{(0, 1)}(\theta) = 0,
\end{equation*}
and thus $\Delta_a(\theta) = - M^{(0, 0)}(\theta)$, $\Delta_{\lambda_1}(\theta)
= - M^{(1, 0)}(\theta)$ and $\Delta_{\lambda_2}(\theta)
= - M^{(0, 1)}(\theta)$.

Equations \eqref{eq1_xy}, \eqref{eq2_00}, \eqref{eq2_xy}, \eqref{eq3_00},
\eqref{eq3_xy} can be written in the format
\begin{equation*}
  \phi(\theta) = l(\theta)\phi(a(\theta)) + \eta(\theta),
\end{equation*}
if the dynamical average $l^* < 1$ (Remark~\ref{remark_dynamical_average}), or
be rewritten in the format
\begin{equation*}
  \phi(\theta) = \frac{1}{l(a^{-1})(\theta)}\phi(a^{-1}(\theta)) - \frac{\eta(a^{-1}(\theta))}{l(a^{-1}(\theta))},
\end{equation*}
if $(\frac{1}{l(a^{-1})})^* < 1$, 
and such equations can then be solved through contraction as in
\eqref{coho_steps} and Algorithm~\ref{algorithm_coho}.

\begin{remark}
  In the 3-dim case, our quasi-Newton method does not work if $\lambda_1(\theta)$ and
  $\lambda_2 (\theta)$ are resonant.

  In fact, 
  in this case, the fibered version of the Sternberg theorem
  fails. We are expecting a fibered version of Chen's Theorem \cite{KC63}. More
  specifically, for $|\lambda_2|^q <
  |\lambda_1|$, we are expecting a polynomial
  $p(\theta, s)$ with degree $\leq q$ such that the invariance equation becomes
  $$
  f(W(\theta, s)) - W(a(\theta), p(\theta, s)) = 0.
  $$
  The same discussions remain valid for higher-dimensional cases. We hope to come
  back to this problem and develop a complete theory.
\end{remark}

\subsection{Numerical Exploration: 3D-Fattened Arnold Family} 
Inspired by \cite{H16}, in this subsection, we implement the algorithm discussed
in Section~\ref{3d_derivation} on a 3-dimensional Fattened Arnold Family (3D-FAF)
$f_{\alpha, \epsilon}: \mathbb{T}^1 \times \mathbb{R}^2 \rightarrow \mathbb{T}^1 \times \mathbb{R}^2$:

\begin{equation}
  \label{3D-FAF}
  f_{\alpha, \epsilon}(x, y, z) =
  \begin{pmatrix}
    x + \alpha + \frac{\epsilon}{2 \pi} (\sin(2 \pi x) + y + \frac{z}{2}) \\
    \beta(\sin(2 \pi x) + y) \\
    \gamma(\sin(2 \pi x) + y + z)
  \end{pmatrix},
\end{equation}
where $\epsilon$ is the perturbation parameter, $\alpha$ is the drift parameter
and $\beta$, $\gamma$ are the eigenvalues. We have implemented the algorithm for
the stable (i.e. $\beta < 1$, $\gamma < 1$) case, unstable ($\beta > 1$, $\gamma >
1$) case and the saddle (i.e. $\beta < 1 < \gamma$) case.

\subsubsection{The Unperturbed Case}
In order to apply the continuation method, we start with the unperturbed case
where $\epsilon = 0$. As discussed in \cite{H16}, in such case, the solution
to the invariance equation \eqref{invariance_3d} is:
\begin{align*}
  W_1(\theta, s_1, s_2) &= \theta, \\
  W_2(\theta, s_1, s_2) &= - S(\alpha, \beta) \cos(2 \pi \theta) + (C(\alpha, \beta) - 1) \sin(2 \pi \theta) + s_1, \\
  W_3(\theta, s_1, s_2) &= \frac{\gamma}{\beta} \big(S(\alpha, \beta) (C(- \alpha, \gamma^{-1}) - 1) + (C(\alpha, \beta) - 1) S(-\alpha, \gamma^{-1})\big) \cos(2 \pi \theta) \\
                        &\ \ \ \frac{\gamma}{\beta} \big( (C(\alpha, \beta) - 1)(C(-\alpha, \gamma^{-1}) - 1) - S(\alpha, \beta) S(-\alpha, \gamma^{-1})\big) \sin(2 \pi \theta)\\
                        &\ \ \ + \frac{\gamma}{\beta - \gamma} s_1 + s_2,\\
  a(\theta) = \theta + \alpha&, \text{ } \lambda_1(\theta) = \beta, \text{ } \lambda_2(\theta) = \gamma, \\
\end{align*}
where 
\begin{equation*}
  S(x, y) = \frac{y \sin(2 \pi x)}{1 - 2 y \cos(2 \pi x) + y^2}, \ 
  C(x, y) = \frac{1 - y \cos(2 \pi x)}{1 - 2 y \cos(2 \pi x) + y^2}. 
\end{equation*}

\subsubsection{The Perturbed Case}
\begin{figure}[!htb]
  \includegraphics[width=0.75\textwidth]{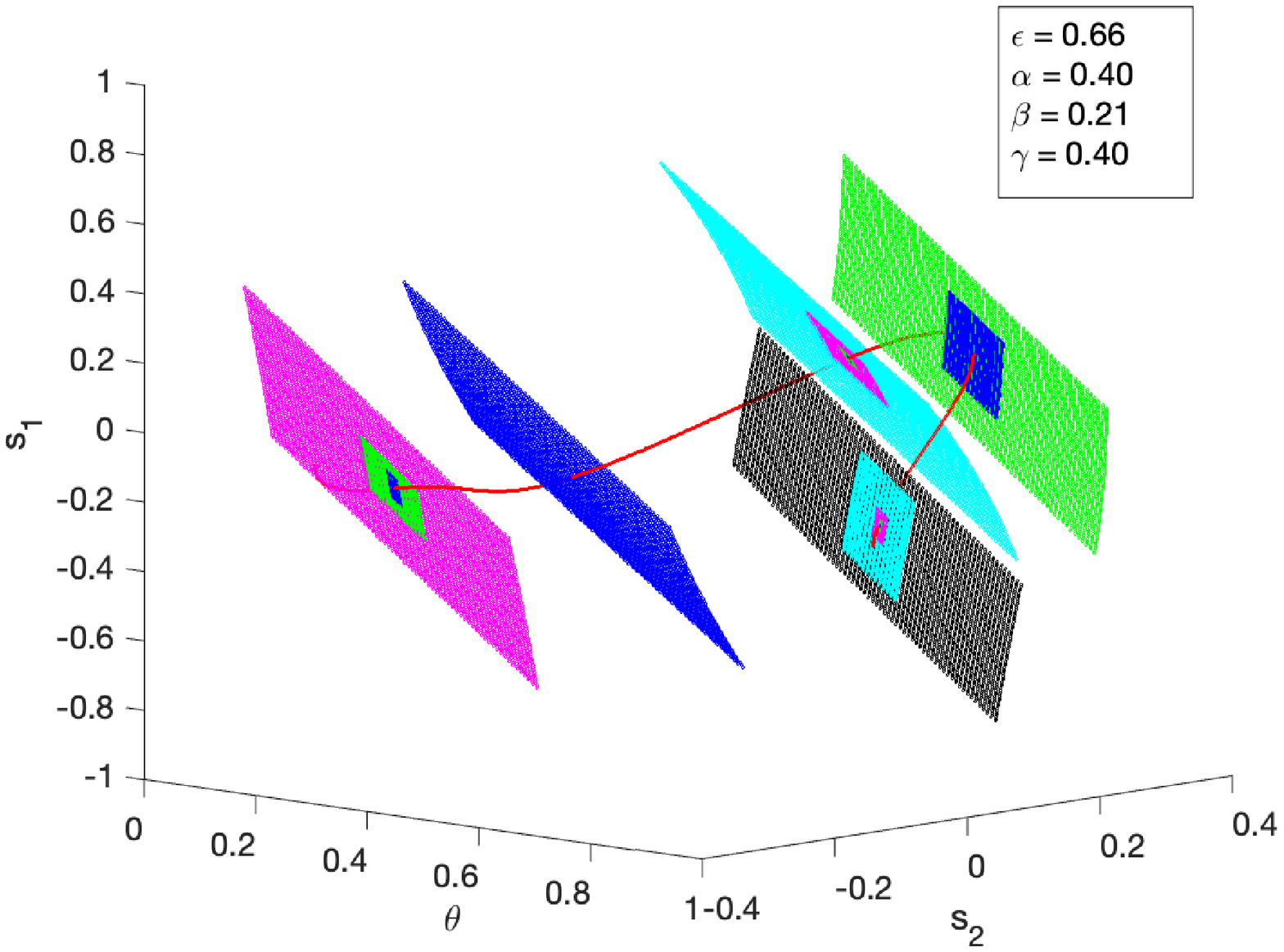}
  \caption{stable}
  \label{fig:stable3d} 
 \footnotesize
  \emph{The isochrons contract along both the eigen-directions.}
\end{figure}

\begin{figure}[!htb]
  \includegraphics[width=0.75\textwidth]{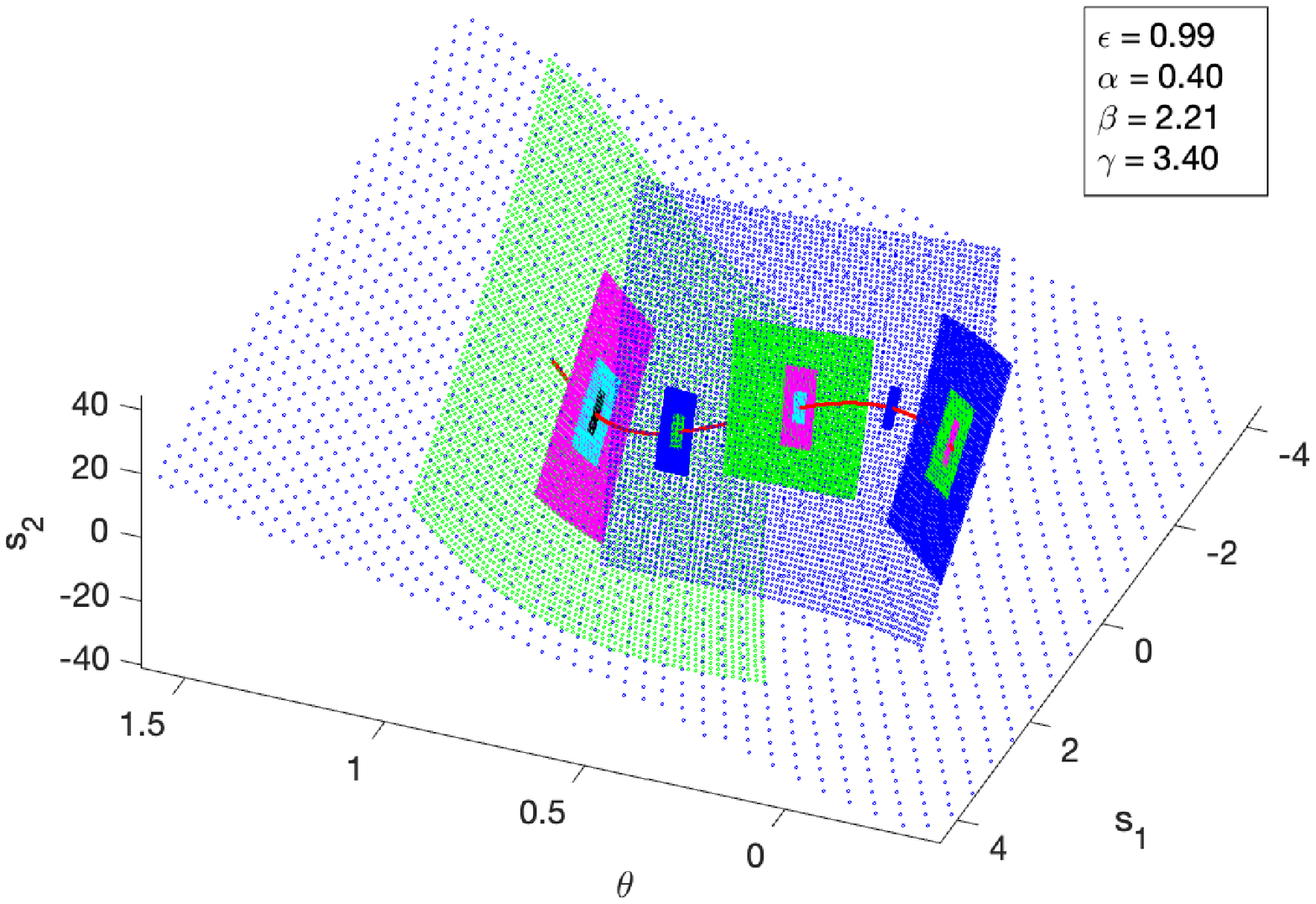}
  \caption{unstable}
  \label{fig:unstable3d}
  \footnotesize
  \emph{The isochrons expand along both the eigen-directions.}
\end{figure}

\begin{figure}[!htb]
  \includegraphics[width=0.75\textwidth]{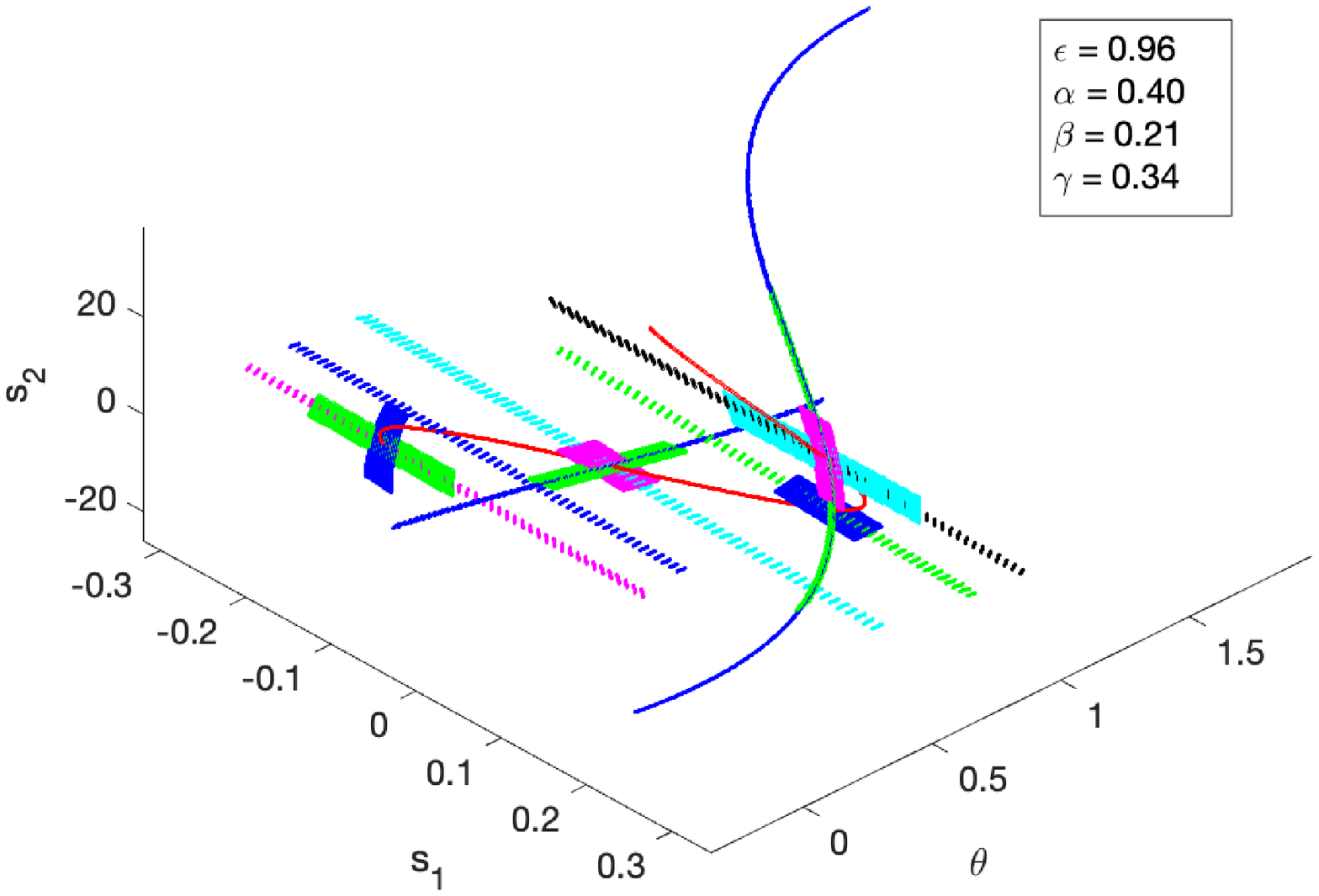}
  \caption{saddle}
  \label{fig:saddle3d}
  \footnotesize
  \emph{The isochrons contract along the stable direction and expand along the unstable direction.}
\end{figure}

Following the same continuation method as in Section~\ref{sec_continuation}, the
invariant circle and the corresponding isochrons are computed for all the
stable, unstable, and saddle choice of the parameters. For the same color choice
as in Figure~\ref{fig:cic1}, we still start with the blue isochron, which is
mapped to
the green isochron, followed by magenta, cyan, and black, correspondingly (see
Figure~\ref{fig:stable3d}, \ref{fig:unstable3d} and \ref{fig:saddle3d}). 


\appendix
\section{A Proposal for a Parallel Implementation of the Algorithm: a Digression}
\label{appendix}

Although the parallelization is not implemented in this paper, we remark that
Algorithm~\ref{algorithm} can be parallelized using a
2-dimensional mesh of processors as the communication network. In this
Appendix, we propose
a draft framework when multiple processors are available. 

This
subsection is  a digression from the main argument of the paper
and we have not implemented the algorithms discussed here. It can 
be omitted without affecting the rest. 

A two-dimensional communication network (with no wraparound link) is an
embedding network 
in the
way that the processors are arranged to a matrix, where each processor can
communicate with their nearest neighbors in the same row or column. This type of
embedding is commonly used in matrix manipulation.

\subsection{Function Representation in 2-dim Network}

In our algorithm, the functions we are dealing with are either of Type-1
(see Section~\ref{function_representation}): i.e. $a(\theta), \lambda(\theta)$,
which are stored as a 1-D array,
or of Type-2: i.e. $W(\theta, s)$, which is stored as a matrix of size $N \times
(L + 1)$. We can block distribute Type-1 functions
onto the first row of the grid, and block distribute Type-2 functions on the 2-D
grid.

More precisely, suppose the 2-D network has size $P_N \times P_L$, and $N, L$ are
integer multiples of $P_N, P_L$, respectively, then each processor on the first
row stores $\frac{N}{P_N}$ points for functions of Type-1, and each processor
stores $\frac{N}{P_N} \times \frac{L}{P_L}$ points for functions of Type-2.

Storing functions blockwise in this way allows the basic functional operations
(for example, summation, subtractions, Type-1 function multiplication) to be
processed within each processor in a pointwise manner without communication. The
delicate part is the operations when communications between processors are
needed. In our algorithm, such operations include function evaluation (for
composition) and Taylor polynomial multiplication.

\subsection{Cubic Spline Interpolation in Parallel}

Interpolations are required in our algorithm whenever composition between
functions or derivatives of functions are needed. Because of the way functions
are stored, cubic spline
interpolations need to be performed on each row of the communication network.

The computation for the coefficients of cubic splines is essentially achieved by
solving a linear system, where the matrix involved is only a tridiagonal matrix.
Thus the goal is to solve this tridiagonal system on a 1-D network.

Following \cite{MPK16}, such system can be reduced to a single variable equation
recursively in $\log_3 N$ steps, which, when implemented parallelly, has time
complexity $\mathcal{O}((\tau + \mu) \log P_N)$ (for communication) and
$\mathcal{O}(\log \frac{N}{P_N})$ (for computation), where
$\tau$ is the latency time, $\mu$ is the inverse of the bandwidth, and $m$ is the
size of the message communicating between processors (in this case is just a
constant). Pipelining techniques \cite{MM11} can be used to speedup this computation.

To perform composition, we need to first use
the hypercubic AllGather technique
to have all the processors in the same row gather all the pieces of the spline
(with communication time complexity $\mathcal{O}(\tau \log P_N + \mu N)$). For
example, when
computing $W(a(\theta), \lambda(\theta)s)$, one need to first BroadCast every
piece of $a(\theta)$ from the first row to the rest rows along the corresponding
columns (with communication time complexity $\mathcal{O}((\tau + \mu) \log
P_N)$), compute $\lambda^i(\theta)$ for row $i$ using parallel prefix sum
technique (with communication time $\mathcal{O}(\tau \log P_L + \mu
\frac{N}{P_N} \log P_L)$ and computation time $\mathcal{O}(\frac{L
  \times N}{P_L \times P_N} + \frac{N}{P_N}\log P_L)$),
perform the interpolation for each row for the corresponding $W^{(i)}(\theta)$, $i
= 0, \dots, L$, AllGather the spline on each processor, and then evaluate on
$a(\theta)$,
and then multiply it by the corresponding $\lambda^i(\theta)$ pointwisely.

\subsection{Type-2 Function Multiplication in Parallel}

In our algorithm, products between Type-2 functions are also required. Given
$g_1(\theta, s) = \sum_{i = 0}^Lg_1^{(i)}(\theta)s^i$, $g_2(\theta, s) = \sum_{i
  = 0}^Lg_2^{(i)}(\theta)s^i$, we need to compute $g_1 \times g_2$ up to the
$L$-th order. To achieve this, one can AllGather $g_2$ for each column, compute
the corresponding coefficients, and then AllReduce (with time complexity
$\mathcal{O}(\tau \log P_L) + \mu \frac{N}{P_N} \log P_L$) along each column to
update the coefficients for the product.



\end{document}